\documentclass[numbib]{imaiai}
\usepackage{algorithm,algorithmicx,algpseudocode}
\usepackage[title]{appendix}
\usepackage[section]{placeins}
\usepackage{amsmath,amsfonts,amssymb,amsthm,mathtools} 
\usepackage{pifont,stmaryrd} 
\usepackage[misc,geometry]{ifsym} 
\usepackage{graphicx,epstopdf,caption,subcaption} 
\usepackage{threeparttable,multirow} 
\usepackage[dvipsnames]{xcolor} 
\usepackage{dsfont} 
\usepackage{bm} 
\usepackage[thinlines]{easytable}
\usepackage{tabularx}

\def \P {{\mathbf P}}

\def \G {{\mathbf G}}
\def \H {{\mathbf H}}

\def \a {{\mathbf a}}

\def \R {{\mathbf R}}

\def \D {{\mathbf D}}

\def \v {{\mathbf v}}
\def \u {{\mathbf u}}
\def \x {{\mathbf x}}
\def \y {{\mathbf y}}

\def \A {{\mathbf A}}
\def \B {{\mathbf B}}

\def \Z {{\mathbf Z}}

\def \B {{\mathbf B}}

\DeclarePairedDelimiter{\abs}{\lvert}{\rvert}
\newcommand{\distas}[1]{\mathbin{\overset{#1}{\kern\z@\sim}}}%
\newcommand{\norm}[1]{\ensuremath{\left\|#1\right\|}}		

\newcommand{\cmark}{\ding{51}}%
\newcommand{\xmark}{\ding{55}}

\definecolor{blue-violet}{rgb}{0, 0, 0}

\newcommand{\revise}[1]{#1}
\renewcommand{\revise}[1]{{\color{blue-violet}#1}} 

\begin{document}

\title{Exit Time Analysis for Approximations of Gradient Descent Trajectories Around Saddle Points}

\shorttitle{Exit Time Analysis of Gradient Descent Trajectories} 
\shortauthorlist{Dixit, G\"urb\"uzbalaban, and Bajwa} 

\author{{
\sc Rishabh Dixit}\\[2pt]
Department of Electrical and Computer Engineering\\
Rutgers University--New Brunswick, NJ 08854 USA\\
{\email{Corresponding author: {\tt rishabh.dixit@rutgers.edu}}}\\[6pt]
{\sc Mert G\"urb\"uzbalaban}\\[2pt]
Department of Management Science and Information Systems\\
Department of Electrical and Computer Engineering\\
Department of Statistics\\
Rutgers University--New Brunswick, NJ 08854 USA\\
{\tt mg1366@rutgers.edu}\\[6pt]
{\sc Waheed U. Bajwa} \\[2pt]
Department of Electrical and Computer Engineering\\
Department of Statistics\\
Rutgers University--New Brunswick, NJ 08854 USA\\
{\tt waheed.bajwa@rutgers.edu}}

\maketitle

\begin{abstract}
{This paper considers the problem of understanding the exit time for trajectories of gradient-related first-order methods from saddle neighborhoods under some initial boundary conditions. Given the `flat' geometry around saddle points, first-order methods can struggle to escape these regions in a fast manner due to the small magnitudes of gradients encountered. In particular, while it is known that gradient-related first-order methods escape strict-saddle neighborhoods, existing analytic techniques do not explicitly leverage the local geometry around saddle points in order to control behavior of gradient trajectories. It is in this context that this paper puts forth a rigorous geometric analysis of the gradient-descent method around strict-saddle neighborhoods using matrix perturbation theory. In doing so, it provides a key result that can be used to generate an approximate gradient trajectory for any given initial conditions. In addition, the analysis leads to a linear exit-time solution for gradient-descent method under certain necessary initial conditions, which explicitly bring out the dependence on problem dimension, conditioning of the saddle neighborhood, and more, for a class of strict-saddle functions.}
{Exit-time analysis; Gradient descent; Non-convex optimization; Strict-saddle property.}\\
\textit{2010 Math Subject Classification}:  {90C26 ; 15Axx ; 41A58 ; 65Hxx}
\end{abstract}

\section{Introduction}
The problem of finding the convergence rate/time of gradient-related methods to a stationary point of a convex function has been studied extensively. Moreover, it has been well established that stronger conditions on function geometry yield better convergence guarantees for the class of gradient-related first-order methods. For instance, conditions like strong convexity and quadratic growth result in the so-called `linear convergence rate' to a stationary point for gradient-related first-order methods. Though there is also a class of second-order (Hessian-related) methods like the Newton method that yield super-linear convergence to stationary points of strongly convex functions, that comes at the cost of very high iteration complexity.

More recently much of the focus has shifted towards obtaining rates of convergence for gradient-related methods to stationary points of non-convex functions. To this end, there are some local geometric conditions like the Kurdyka-{\L}ojasiewicz property~\cite{Kurdyka.ADLF98,lojasiewicz1961probleme} that guarantee linear convergence rates provided the iterate is in some bounded neighborhood of the function's second-order stationary point~\cite{LiPong.FCM18}. Such guarantees, however, are hard to obtain for non-convex functions in a global sense and the linear convergence rates are often eventual, i.e., these methods usually exhibit such linear convergence only asymptotically. The main reason that restricts this speedup behavior to the asymptotic setting is the non-convex geometry that can impede fast traversal of these methods across the geometric landscape of the function. This is due to the fact that trajectories of gradient-related methods can encounter extremely flat curvature regions very near to first-order saddle points. Such regions are characterized by gradients that have very small magnitudes and it can take exponential time for the trajectory of an algorithm to traverse this extremely flat region. A natural question to ask then is whether there exist gradient-related first-order methods for which a subset of non-zero measure trajectories escape first-order saddle points of a class of non-convex functions in `linear' time.\footnote{We are slightly abusing terminology here and, in keeping with the convention of linear convergence rates in optimization literature, we are defining `linear exit time' for the trajectory of a discrete method to be one in which the trajectory escapes an $\mathcal{O}(\epsilon)$ saddle neighborhood in $\mathcal{O}(\log(\epsilon^{-1}))$ number of iterations.} The non-zero measure of such fast escaping trajectories is important since studying fast escape is only useful when the initialization set is dense in such trajectories. Section \ref{crudeanalysissection} (see Remark \ref{remarkmeasure}) in particular establishes that indeed fast saddle escape is possible from an initialization set of positive measure.

We address this question in this work by deriving an upper bound on the exit time for a certain class of gradient-descent trajectories escaping some bounded neighborhood of the first-order saddle point of a class of smooth, non-convex functions. Specifically, let $\x^*$ be a saddle point of a smooth, non-convex function $f: \mathbb{R}^n \rightarrow \mathbb{R}$ and, without loss of generality, define the bounded neighborhood around the saddle point to be an open ball of radius $\epsilon$ around $\x^*$, denoted by $\mathcal{B}_{\epsilon}(\x^*)$. Recall that the gradient at saddle point $\x^*$ is a zero vector, i.e., it is necessarily a first-order stationary point. In addition, the saddle neighborhood $\mathcal{B}_{\epsilon}(\x^*)$ exhibits certain properties that depend on Lipschitz boundedness of the function and its derivatives as well as eigenvalues of the Hessian at $\x^*$. The class of trajectories we focus on in here is assumed to have the current iterate sitting on the boundary of $\mathcal{B}_{\epsilon}(\x^*)$ and it comprises of all those trajectories of gradient descent that escape this saddle neighborhood with \emph{at least} linear rate. Note that the current iterate could have reached the boundary of $\mathcal{B}_{\epsilon}(\x^*)$ using \emph{any} gradient-related method, but that problem is not our concern. Rather, our focus here is whether there exists any gradient-descent trajectory from the current iterate that can escape $\mathcal{B}_{\epsilon}(\x^*)$ in almost linear time of order $\mathcal{O}\left(\log(\epsilon^{-1})\right)$ or better. And if such a trajectory exists, then an immediate subsequent question asks for the necessary conditions required for the existence of such gradient-descent trajectories. \revise{To answer both these questions effectively, we present a rigorous analysis of gradient-descent trajectories $\{\x_k\}$ starting at time $k=0$, when the initial iterate $\x_{0}$ sits on the boundary of the ball $\mathcal{B}_{\epsilon}(\x^*)$, till the time they exit $\mathcal{B}_{\epsilon}(\x^*)$, which we term the \emph{exit time} and denote by $K_{exit}$. It should be noted that we analyze in this work the first-order approximations of the exact trajectories, instead of the exact trajectories themselves, where the approximation error is sufficiently small. Specifically, the presence of higher-order terms ($\mathcal{O}(\epsilon^2)$ terms) in the forthcoming analysis accounts for the approximation in our analysis, while things are proved about trajectories and perturbations up to the first order in $\epsilon$.}

We conclude by noting the relevance of the exit-time results derived in this paper to the broader field of non-convex optimization. First, to the best of our knowledge, there are no works other than the ones listed in Table~\ref{table:2} that \emph{explicitly} investigate the exit times from saddle neighborhoods of the trajectories of discrete first-order methods. Rather, the focus in much of the related works discussed in Section~\ref{ssec:prior.work} is on providing rates of convergence to second-order stationary points. While such analysis necessarily implies saddle escape, this is typically accomplished through the use of noisy perturbations that allow the trajectories to move along a negative curvature direction; in particular, such approaches do not yield an explicit expression for the exit time of a trajectory from a saddle neighborhood. Second, and most importantly, the rate of convergence to a second-order stationary point is trivially a function of the time a trajectory spends near a saddle point. It therefore stands to reason that the existing convergence rates for some of the recent first-order methods can possibly be improved by identifying trajectories with linear exit time, which is the focus of this paper.

\subsection{Relation to prior work}\label{ssec:prior.work}
Convergence rates of optimization methods to the minima of convex functions have been studied for quite some time. For instance, the seminal work dealing with convergence rate analysis of gradient-related methods has been well summarized in \cite{polyak1964some}, while a recent work by \cite{nesterov2006cubic} summarizes convergence rates of Newton-type methods. These prior works rely heavily on the Lipschitz boundedness of the function along with some other form of curvature property. The works  \cite{attouch2013convergence} and \cite{bolte2007lojasiewicz} utilize the local Kurdyka--{\L}ojasiewicz property \cite{Kurdyka.ADLF98,lojasiewicz1961probleme} of a function to develop convergence guarantees and the ergodic rates using monotonicity of gradient sequences in a bounded neighborhood of the function's stationary point. However, for non-convex functions these seminal works do not analyze the exit time from a bounded neighborhood of a first-order saddle point. With the focus shifting towards characterizing the efficacy of gradient-related methods on non-convex geometries in recent years, it becomes imperative to conduct such an analysis. To the best of our knowledge, currently no work exists that analyzes (discrete) gradient-descent trajectories in the saddle neighborhood using eigenvector perturbations. Therefore, this is the first work that incorporates matrix perturbation theory to extract the local geometric information around a saddle point necessary for analyzing gradient trajectories at such small scales. As a result of the perturbation analysis, the hidden dependence of exit time on the trajectory's initialization point, conditioning of the saddle neighborhood, problem dimension, and more, is also revealed in this work (cf.~Table~\ref{table:1} in Section~\ref{exactexittimeanalysissection}).

There is a plethora of existing methods in the literature that deal with non-convex optimization problems. Within the context of this paper, we broadly classify these methods into \emph{continuous-time} Ordinary Differential Equations (ODE)-type methods/analysis and \emph{discrete-time} gradient-descent related algorithms/analysis. The latter class of methods can be further categorized into first-order and higher-order methods. Starting with the continuous-time ODE-type algorithms, we first refer to \cite{bolte2007lojasiewicz} that has developed upon the gradient flow curve analysis of non-smooth convex functions. Although this work focuses on convex problems, yet it is important in the sense that it motivates us in drawing some parallels between the discrete gradient trajectories and the continuous flow curves in our analysis of non-convex functions.

Another recent work \cite{hu2017fast} within the continuous-time setting analyzes the saddle escape problem using a stochastic ODE to characterize the rates of escape in terms of a multiplicative noise factor. Remarkably, the results in \cite{hu2017fast} give a linear rate of escape in expectation for very small stochastic noise. This work also extends these results to cascaded saddle geometries. Note that the analysis in \cite{hu2017fast} relies on an earlier important work by \cite{kifer1981exit}, which characterizes the probability distribution of the exit time of gradient curves from saddle point vicinities. The hyberbolic flow curves discussed in \cite{bolte2007lojasiewicz,hu2017fast,kifer1981exit} are the building blocks of our intuition towards analyzing discrete gradient trajectories in this work.

A Stochastic Differential Equation (SDE) approach has also been utilized in a recent work \cite{ShiSuEtAl.arxiv20} to study gradient-based (stochastic) methods for non-convex functions in the continuous-time setting. While this work also guarantees linear rates of global convergence for non-convex problems under certain assumptions, a few of which are more restrictive than our work, it does not lend itself to understanding the behavior of discrete gradient trajectories around first-order saddle points. Similarly the analysis done in \cite{murray2019revisiting} shows that fast evasion is possible for trajectories generated by normalized gradient flow from strict saddle neighborhoods of Morse functions but such an analysis is not sufficient to explain the behavior of discrete trajectories around saddle points.

Next, there exists a large collection of work analyzing discrete gradient-related methods in non-convex settings. The very basic yet most often investigated approach in these works is the Stochastic Gradient Descent (SGD) method and its variants. Such methods have been extensively studied in the literature for the purpose of escaping saddles, specifically first-order saddle points. For instance, \cite{du2017gradient,jin2017escape} provide the rates of convergence to a second order stationary point with very high probability using perturbed gradient descent, where the perturbation vector is an isotropic noise. In contrast, \revise{the work in \cite{du2017gradient} shows that---in the worst case---the time to escape cascaded saddles scales exponentially with the problem dimension}, thereby making the method impractical for highly pathological problems like optimization over jagged functions.

The work in \cite{lee2017first} provides new insights into the efficacy of gradient-descent method around strict saddle points. The authors in this work present a measure-theoretic analysis of the gradient-descent trajectories escaping strict saddle points almost surely. Their analysis uses the stable center manifold theorem in \cite{kelley1966stable} to prove that random initializations of gradient-descent trajectories in the vicinity of a strict saddle point almost never terminate into this saddle point. Note that while this is an intuitive inference, it is somewhat hard to prove for gradient flow curves around saddle points. The work  \cite{daneshmand2018escaping} also provides rates and escape guarantees under certain strong assumptions of high correlation between the negative curvature direction and a random perturbation vector. Interestingly, the convergence rate put forth in this work does not depend on the problem dimension. However due to the nature of the somewhat restrictive assumptions in \cite{daneshmand2018escaping}, the resulting method is not suited to work over a general class of non-convex problems. We also note two related recent works \cite{ErdogduMackeyEtAl.ConfNIPS18,RaginskyRakhlinEtAl.ConfCOLT17} that analyze global convergence behavior of Langevin dynamics-based variants of the SGD (and simulated annealing) for non-convex functions. Neither of these works, however, focuses on the escape behavior of trajectories around saddle neighborhoods.

There also is a sub-category of first-order methods leveraging acceleration and momentum techniques to escape saddle points. For instance, \cite{o2019behavior} uses the stable center manifold theorem to show that the heavy-ball method almost surely escapes a strict saddle neighborhood. But the rate of escape derived in this work is limited to quadratic functions; further, the ensuing analysis does not bring out the dependence on problem dimension, conditioning of the saddle neighborhood, etc. The work in \cite{reddi2017generic} provides extensions of SGD methods like the Stochastic Variance Reduced Gradient (SVRG) algorithm for escaping saddles. Recently, in works like \cite{jin2017accelerated} and \cite{xu2018first}, methods approximating the second-order information of the function (i.e., Hessian) have been employed to escape the saddles and at the same time preserve the first-order nature of the algorithm. Specifically, \cite{jin2017accelerated} shows that the acceleration step in gradient descent guarantees escape from saddle points with provably better rates; yet the rate is still worse than the linear rate. Along similar lines, the method in \cite{xu2018first} utilizes the second-order nature of the acceleration step combined with a stochastic perturbation to guarantee escape and provide escape rates.

Finally, higher-order methods are discussed in \cite{mokhtari2018escaping,paternain2019newton}, which utilize the Hessian of the function or its combinations with first-order algorithms to escape saddle neighborhoods with an impressive super linear rate while trading-off heavily with per-iteration complexity. Going even a step further, the work in \cite{anandkumar2016efficient} poses the problem with second-order saddles, thereby making higher-order methods an absolute necessity. Though these techniques optimize well over certain pathological functions like degenerate saddles or very ill-conditioned geometries, yet they suffer heavily in terms of complexity. In addition, none of these methods leverage the initial boundary condition of their methods around saddle points, which could not only influence the future trajectory but also control its exit time from some bounded neighborhood of the saddle point. This further motivates us to conduct a rigorous analysis of (approximations of) gradient-descent trajectories around saddle points for some fixed initial boundary conditions.\looseness=-1

We conclude by noting that the use of careful initial boundary conditions in order to avoid saddle points in non-convex optimization is not a fundamentally new idea. Consider, for instance, the non-convex formulation of the phase retrieval problem in \cite{candes2015phase}. A variant of the gradient descent method, termed the Wirtinger flow algorithm, can be used to solve this problem as long as the algorithm is carefully initialized along the direction of the negative curvature by means of a spectral method \cite{candes2015phase}. However, one of the implications of the results in this paper are that spectral initializations such as the one in \cite{candes2015phase}, which require costly computation of the dominant eigenvector of a matrix, are not always required for saddle escape. Rather, one might be able to escape the saddle neighborhoods in approximately linear time provided the projection of the initial gradient descent iterate along the negative curvature direction is lower bounded by a small quantity.
	
\subsection{Our contributions}
Having discussed the relevant works pertaining to the problem of characterizing the exit time of first-order methods from saddle neighborhoods, we now elaborate upon the contributions of our work.

First, none of the earlier discussed works exploit the dependence of the function gradient in saddle neighborhood on the eigenvectors of the Hessian at the saddle point. This dependence results from the eigenvector perturbations of the Hessian in the saddle neighborhood. Therefore, to our knowledge, this is the first work that utilizes the Rayleigh--Schr\"{o}dinger perturbation theory to approximate the Hessian $\nabla^2 f(\x)$ at any point $\x \in \mathcal{B}_{\epsilon}(\x^*)$. This approximate Hessian is then used to obtain the function gradient $\nabla f(\x)$ for any point $\x \in \mathcal{B}_{\epsilon}(\x^*)$.

Second, using the value of the function gradient, for any given initialization $\x_0$ and some fixed step size, we generate an approximate trajectory for the gradient-descent method inside the ball $\mathcal{B}_{\epsilon}(\x^*)$. As a consequence, we obtain the distance between the saddle point $\x^*$ and any point on the approximate trajectory inside the ball $\mathcal{B}_{\epsilon}(\x^*)$ as a function of (discrete) time. Once this distance function is known, we can estimate the exit time of the approximate trajectory from the ball $\mathcal{B}_{\epsilon}(\x^*)$. In this vein, we develop an analytical framework in this work that approximates the trajectory for gradient descent within the saddle neighborhood and establish the fact that a linear escape rate from the saddle neighborhood is possible for some approximate trajectories generated by the gradient-descent method.

Third, we utilize the initial conditions on our iterate by projecting it onto a stable and an unstable subspace of the eigenvectors of the Hessian at the saddle point. This is extremely important since the escape rate and the associated necessary conditions are heavily dependent on where the iterate or gradient trajectory started. To this end, we simply make use of the strict saddle property to split the eigenspace of the Hessian at the saddle point into orthogonal subspaces of which two are of interest, namely, the stable subspace and the unstable subspace.\footnote{There can be one more orthogonal subspace corresponding to the zero eigenvalues of the Hessian at a strict saddle point. Under the assumption of the function being a \emph{Morse function}, however, this subspace vanishes.} Taking the inner product of the iterate with these subspaces yields the respective projections. (Note that this analysis of ours can be readily adapted to obtain these projections for \emph{any} gradient-related method.) As a consequence, for any given initialization of our iterate within the saddle neighborhood, we provide the approximate iterate expression for the entire trajectory as long as it stays within this saddle neighborhood.

Finally, and most importantly, this work provides an upper bound on the exit time $K_{exit}$ for approximations of (discrete) gradient-descent trajectories that is of the order $\mathcal{O}(\log({\epsilon}^{-1}))$, where the constants inside the $\mathcal{O}(\cdot)$ term explicitly depend on the condition number, dimension, and eigenvalue gap, as detailed in Section~\ref{exactexittimeanalysissection}. Also, we develop a necessary condition on the initial iterate that is required for the existence of this exit time. It is worth noting that though the trajectory analysis developed in this work for the gradient-descent method is only approximate, we show in a follow-up work \cite{dixit2022boundary} that this approximation can only have a maximum relative error of order $\mathcal{O}(\log^2(\epsilon^{-1})\epsilon^{3/2})$, provided the exit time $K_{exit}$ is at most of the order $ \mathcal{O} (\log (\epsilon^{-1}))$. Therefore our approximate analysis of the gradient-descent trajectories and their time of exit from the saddle neighborhood can be readily adapted to develop efficient algorithms for escaping first-order saddle points at a linear rate. One such algorithm has already been developed in \cite{dixit2022boundary}, which extends the boundary conditions developed in this work for the linear exit time gradient trajectories and escapes saddle neighborhoods in linear time. The algorithm is designed to check the initial boundary conditions, after which it decides to either keep traversing along the same gradient trajectory or switch to a higher-order method for one iteration. To get a detailed understanding of this extension of our current work, we refer the reader to \cite{dixit2022boundary}.

\begin{table}[t]
\fontsize{50}{12}\selectfont
\centering
\resizebox{1.05\textwidth}{!}
{%
\hspace{-70pt}
\renewcommand{\arraystretch}{12}
\begin{tabular}{||c c c c c c||}
 \hline
 \textbf{Reference} & \textbf{Nature of analysis} & \textbf{Dynamical system analyzed} &
 \textbf{Function class} &
 \textbf{Exit time bound} & \textbf{Necessary initial conditions}
 \\ [0.5ex]
 \hline\hline
 \cite{lee2016gradient} & Asymptotic & Gradient descent method & $\mathcal{C}^2$ functions  & \xmark  & \xmark
   \\
 \hline
 \cite{o2019behavior} & Asymptotic & Heavy ball method & $\mathcal{C}^2$    functions & \xmark  & \xmark
   \\
 \hline
 \cite{o2019behavior} &  Non-asymptotic &  General accelerated methods, & Quadratics ($\langle\x, \A\x \rangle$) & $\mathcal{O}(\log(\frac{1}{\Delta}))
   $ iterations from the &  $\norm{\pi_{\mathcal{E}_{US}}(\x_0-\x^*)} \geq \Delta$  \\ & & Gradient descent method & &  unit ball $\mathcal{B}_{1}(\x^*)$ &\\
  \hline
  \cite{murray2019revisiting} &  Non-asymptotic &  Normalized gradient flow & $\mathcal{C}^2$  Morse functions & $\mathcal{O}(\epsilon)
   $ exit time from a & $\x_0 \neq \x^*$ \\ & & & & small neighborhood $\mathcal{B}_{\epsilon}(\x^*)$ & \\
 \hline
  \cite{hu2017fast} &  Non-asymptotic &  SDE-based gradient flow & $\mathcal{C}^2$ Morse functions & $\mathcal{O}(\log(\frac{1}{\tau}))
   $ mean exit time from some & \xmark \\ & & & & open neighborhood; $\tau$ is the scale & \\ & & & &  of random perturbation &\\
 \hline
  \cite{paternain2019newton} &  Non-asymptotic &  Newton-based method & $\mathcal{C}^2$ functions  & $\mathcal{O}(\log(\frac{1}{\Delta}))
   $ iterations from some &  $ \norm{\pi_{\mathcal{E}_{US}}(\nabla f(\x_0))} \geq \Delta $\\ & & & & open neighborhood  & \\
 \hline
\textbf{This work} &  Non-asymptotic &  Gradient descent method & Locally $\mathcal{C}^{\omega}$ Morse functions, & $\mathcal{O}(\log(\frac{1}{\epsilon}))
   $ iterations from the  & $\norm{\pi_{\mathcal{E}_{US}}(\x_0-\x^*)}^2 \geq \Delta > \Omega(\epsilon) $ \\
  & & & $\mathcal{C}^2$    Morse functions &  ball $\mathcal{B}_{\epsilon}(\x^*)$ &    \\
 \hline
\end{tabular}%
\renewcommand{\arraystretch}{1}
}
\caption{\small\revise{Summary of the similarities and differences between this work and some related prior works.}}
\label{table:2}
\end{table}

We conclude with Table \ref{table:2}, which highlights the similarities and differences between this work and other prior works that explicitly investigate the problem of characterizing the exit time from saddle neighborhoods. The \textit{asymptotic} analyses in this table refer to works that provide measure-theoretic results in terms of the non-convergence of trajectories to a strict saddle point, whereas the \textit{non-asymptotic} works deal with the analysis of trajectories exiting local saddle neighborhoods. The function classes $\mathcal{C}^2$ and $\mathcal{C}^{\omega}$ in the table represent twice continuously differentiable functions and analytic functions, respectively, while the class of quadratics ($\langle \x, \A \x \rangle$) represents functions with constant Hessian. The class of Morse functions is defined in Assumption~\textbf{A4} in the next section. The map $\pi_{\mathcal{E}_{US}}(.)$ is the projection map onto the unstable subspace $\mathcal{E}_{US}$ of the Hessian $\nabla^2 f(\x^*)$, where this subspace will be formally defined in Lemma \ref{lem2}. Notice that the references \cite{murray2019revisiting, hu2017fast} provide exit times from a strict saddle neighborhood for the class of $ \mathcal{C}^2$ functions but analyze continuous time dynamical systems, whereas this work provides the exit time analysis for the gradient descent method, which is a discrete dynamical system. Similarly the work \cite{o2019behavior} develops escape rates for discrete dynamical systems like gradient descent and the heavy ball method but restricts itself to the class of quadratic functions. The only work that develops escape rates for a discrete dynamical system on the class of $\mathcal{C}^2$ functions is \cite{paternain2019newton} but that analysis is for a second-order Newton based method whereas we provide an exit time bound for a first-order method.

\subsection{Notation}  All vectors are in bold lower-case letters, all matrices are in bold upper-case letters, $\mathbf{0}$ is the $n$-dimensional null vector, $\mathbf{I}$ represents the $n \times n$ identity matrix, and $\langle\cdot, \cdot\rangle $ represents the inner product of two vectors. In addition, unless otherwise stated, all vector norms $\norm{\cdot}$ are $\ell_2$ norms, while the matrix norm $\|\cdot\|_2$ denotes the operator norm. Also, for any matrix expressed as $\Z+\mathcal{O}(c)$ where $c$ is some scalar, the matrix-valued perturbation term $\mathcal{O}(c)$ is with respect to the Frobenius norm. Further, the symbol $ (\cdot)^T$ is the transpose operator, \revise{the symbols $\mathcal{O}$, $\Omega$, and $\Theta$ represent the Big-O, Big-Omega, and Big-Theta notation, respectively}, and $W(\cdot)$ is the Lambert $W$ function~\cite{corless1996lambertw}. Throughout the paper, $t$ represents the continuous-time index, while $k,K$ are used for the discrete time. Next, $\gtrapprox$ and $\lessapprox$ mean `approximately greater than' and `approximately less than', respectively. Finally, the operator $\mathbf{dist}(\cdot,\cdot)$ returns the distance between two sets, $\mathbf{diam}(\cdot)$ returns the diameter of a set, \revise{and all the eigenvectors in this work are normalized to be unit vectors}.
	
\section{Problem formulation}
\begin{figure}
 \centering
 \includegraphics[width=0.4\textwidth]{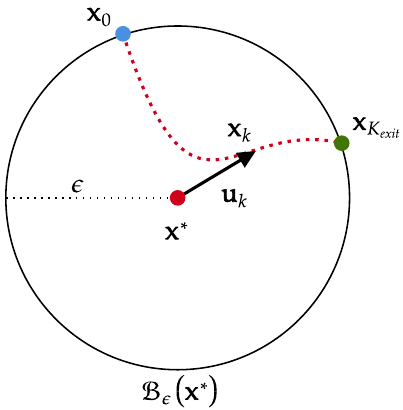}
 \caption{The radial vector evolution in a saddle neighborhood for a function defined on $\mathbb{R}^2$.}
 \label{exitballfig}
\end{figure}

Consider a non-convex smooth function $f(\cdot)$ that has strict first-order saddle points in its geometry. By strict first-order saddle points, we mean that the Hessian of function $f(\cdot)$ at these points has at least one negative eigenvalue, i.e., the function has negative curvature. Next, consider some neighborhood around a given saddle point. Formally, let $\x^*$ be some first-order strict saddle point of $f(\cdot)$ and let $\mathcal{B}_{\epsilon}(\x^*)$ be an open ball around $\x^*$, where $\epsilon$ is sufficiently small. We then generate a sequence of iterates $\x_{k}$ from a gradient-related method on the function $f(\cdot)$, where we call the vector $\u_k = \x_k -\x^*$ inside the ball $\mathcal{B}_{\epsilon}(\x^*)$ the \textbf{radial vector} (see Figure \ref{exitballfig}). Also, it is assumed that the initial iterate $\x_{0} \in \mathcal{\bar{B}}_{\epsilon}(\x^*) \backslash \mathcal{B}_{\epsilon}(\x^*)$, where $\mathcal{\bar{B}}_{\epsilon}(\x^*)$ is the closure of set $\mathcal{B}_{\epsilon}(\x^*)$. With this initial boundary condition, we are interested in analyzing the behavior of our gradient-related sequence $\x_{k}$ in the vicinity of saddle point $\x^*$. More importantly, we are interested in finding some $K_{exit}$ for which the subsequence $\{\x_{k}\}_{k>K_{exit}}$ lies outside $\mathcal{B}_{\epsilon}(\x^*)$ and establishing that $K_{exit} \leq \mathcal{O}(\log({\epsilon}^{-1}))$. Finally, we have to obtain any necessary conditions on $\x_{0}$ that are required for the existence of this `linear' exit time $K_{exit}$.

\subsection{Assumptions}
Having briefly stated the problem, we formally state the set of assumptions that are required for this problem to be addressed in this work.

\begin{figure}
    \centering
    \begin{tabular}{ccc}
         \includegraphics[height=0.8in]{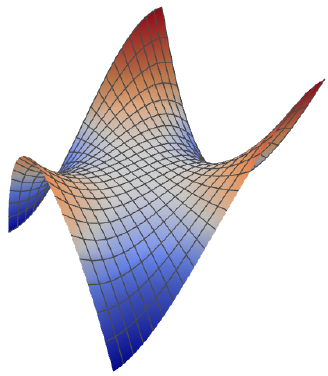} & \includegraphics[height=0.8in]{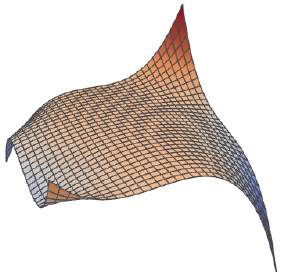} & \includegraphics[height=0.8in]{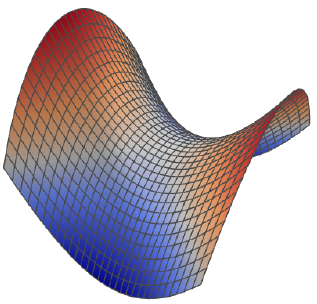}\\
         {Non-strict saddle} & {Degenerate strict saddle} & {\bf Morse function strict saddle}
    \end{tabular}
    \caption{Possible cases of saddle points where the first figure corresponds to a monkey saddle, the second figure is a strict saddle with non-invertible Hessian at the saddle point and the third figure is strict saddle with invertible Hessian at the saddle point.}
    \label{saddlegeo}
\end{figure}

\begin{itemize}
    \item[] \textbf{A1.} \textit{The function $f:\mathbb{R}^n \to \mathbb{R}$ is globally $\mathcal{C}^2$, i.e., twice continuously differentiable, and locally $\mathcal{C}^{\omega}$ in sufficiently large neighbodhoods of its saddle points, i.e., all the derivatives of this function are continuous around saddle points and the function $f(\cdot)$ also admits Taylor series expansion in these neighborhoods.}
    \item[] \textbf{A2.} \textit{The gradient of the function $f(\cdot)$ is $L-$Lipschitz continuous:} $ \norm{\nabla f(\x) - \nabla f(\y)} \leq L \norm{\x - \y}$.
    \item[] \textbf{A3.} \textit{The Hessian of the function $f(\cdot)$ is $M-$Lipschitz continuous:} $\norm{\nabla^{2} f(\x) - \nabla^{2} f(\y)}_2 \leq M \norm{\x - \y}$.		
    \item[] \textbf{A4.} \textit{The function $f(\cdot)$ has only well-conditioned first-order stationary points, i.e., no eigenvalue of the function's Hessian is close to zero at these points (see Figure \ref{saddlegeo}). Formally, if $\x^{*}$ is the first-order stationary point for $f(\cdot)$, then we have}
		\begin{align}
		\nabla f(\x^{*}) &= \mathbf{0},  \ \text{and}\nonumber \\
		\min_i \abs{\lambda_i(\nabla^{2}f(\x^*))} &> \beta, \nonumber
		\end{align}
		\textit{where $\lambda_i(\nabla^{2}f(\x^*))$ denotes the $i^{th}$ eigenvalue of the matrix $ \nabla^{2}f(\x^*)$ and $\beta > 0$. Note that such a function is termed a Morse function.}
\end{itemize}

\revise{We now make a few remarks concerning these assumptions as well as their implications.} Notice that Assumption \textbf{A1} requires $f(\cdot)$ to be locally real analytic, which may seem too restrictive to some readers since the theory of non-convex optimization is often developed around only the assumption that $f \in \mathcal{C}^2$ with Lipschitz-continuous Hessian. It is worth reminding the reader, however, that many practical non-convex problems such as quadratic programs, low-rank matrix completion, phase retrieval, etc., with appropriate smooth regularizers satisfy this assumption of real analyticity around the saddle neighborhoods; \revise{see, e.g., the formulations discussed in \cite{ma2020implicit, chen2019gradient}}. Similarly, the loss functions in deep neural networks with analytic activation functions also satisfy Assumption \textbf{A1} \revise{under certain mild conditions~\cite{kurochkin2021neural}}. It is also worth noting here that Assumption \textbf{A1} enables highly precise estimates of the exit time and the initial boundary condition, which is something that does not happen when dealing with purely $\mathcal{C}^2$ functions; see Section~\ref{crudeanalysissection} for further discussion on this topic. \revise{Next, Assumptions \textbf{A2} and \textbf{A3} are satisfied locally around any saddle point since any locally analytic function is locally $\mathcal{C}^{\infty}$ smooth and therefore is gradient and Hessian Lipschitz continuous in some compact neighborhood of the saddle point.}\looseness=-1

\revise{Lastly, the problem formulation in this work assumes the class of Morse functions (Assumption~\textbf{A4}), i.e., functions whose Hessians are invertible at their critical points. Since Morse functions can only have isolated critical points \cite{matsumoto2002introduction}, the insights from this work are not \emph{directly} applicable to non-convex optimization problems with connected saddle points. While this may appear to be a limitation of this work, Morse functions are an important tool in the study of general non-convex optimization problems since they are dense in the class of $\mathcal{C}^2$ functions \cite{matsumoto2002introduction}. It is therefore no surprise that they are routinely invoked in the non-convex optimization literature (see, e.g., \cite{mei2018landscape,mokhtari2019efficient,yang2021fast}), while neural networks with smooth activation functions are also known to be Morse functions under certain mild assumptions \cite{kurochkin2021neural}. Additionally, since connected saddle points for smooth functions generally arise only when their Hessian at the critical points has one or more zero eigenvalues, one could always add a quadratic regularization term with a sufficiently small constant to any smooth function so as to make the Hessian of the function invertible at its critical points and thus transform the function into a Morse function. As an example, we have circumvented the problem of connected saddle points within the low-rank matrix factorization problem in our follow-up work \cite{dixit2022boundary} by adding a regularization term that makes the objective function a Morse function.}

\revise{Assumption~\textbf{A4} also implies the following two propositions, both of which will be routinely invoked as part of the forthcoming analysis.}
\begin{proposition}\label{morseprop}
\revise{Under Assumption~\textbf{A4}, the function $f(\cdot)$ has only first-order saddle points in its geometry. Moreover, these first-order saddle points are strict saddle, i.e., for any first-order saddle point $\x^*$, there exists at least one eigenvalue $\lambda_i$ of $\nabla^2 f (\x^*)$ that satisfies $\lambda_i(\nabla^{2}f(\x^*)) < - \beta$.}
\end{proposition}
\begin{proof}
\revise{For any $ \mathcal{C}^m$-smooth function $f(\cdot)$ with $m \geq 2$, if $\x^*$ is its second- or higher-order saddle point then it must necessarily satisfy $ \nabla f(\x^*) = \mathbf{0}$ and $\nabla^2 f(\x^*) \succeq \mathbf{0}$, where at least one of the eigenvalues of $ \nabla^2 f(\x^*)$ is $0$. But this is not possible in our case because of Assumption~\textbf{A4}. The fact that an eigenvalue $\lambda_i$ exists such that $\lambda_i(\nabla^{2}f(\x^*)) < - \beta$ is also a direct consequence of Assumption~\textbf{A4}.}
\end{proof}

\begin{proposition}\label{eigenprop}
\revise{Under Assumption \textbf{A4}, for any sufficiently small $\epsilon$ where $\epsilon \ll \beta$, we can group the eigenvalues of the Hessian $ \nabla^{2}f(\x^*)$ at any strict saddle point $\x^*$ into $m$ disjoint sets $\{\mathcal{G}_1, \mathcal{G}_2, \dots, \mathcal{G}_m \}$ with $2\leq m\leq n$ based on the level of degeneracy of eigenvalues (closeness to one another) such that for some $\delta = \Omega(\epsilon^{1-a})$ where $a \in (0,1]$, we have the following conditions:
\begin{align}
			\mathbf{dist}(\mathcal{G}_p,\mathcal{G}_q) &\geq \delta \hspace{0.2cm} \forall \hspace{0.2cm} \mathcal{G}_p,\mathcal{G}_q \hspace{0.2cm} \text{s.t.} \hspace{0.2cm} p \neq q, \ \text{and}
		\\
		\max_{p}\{\mathbf{diam} (\mathcal{G}_p)\}  &= \mathcal{O}(\epsilon^{1-a}).			
\end{align}}
\end{proposition}
\begin{proof}
\revise{From Assumption \textbf{A4}, the eigenvalues of the Hessian $ \nabla^{2}f(\x^*)$ at any strict saddle point $\x^*$ can always be separated into two distinct groups, one consisting of positive eigenvalues and the other comprising negative eigenvalues. By this construction, the distance between these groups will be at least $2 \beta$. Since $\epsilon \ll \beta$, we get a $\delta = 2 \beta$ for this construction which satisfies the constraint $\delta = \Omega(1)$. Next, we check whether the diameter of these two groups is larger than $\Theta(\epsilon^{1-a})$; if yes then we split that particular group into two more groups at the first eigenvalue where the consecutive eigenvalue gap within that group exceeds $\Theta(\epsilon^{1-a})$. This eigenvalue gap becomes our new $\delta$ and by construction it will satisfy the constraint $\delta = \Omega(\epsilon^{1-a})$ for some $a > 0$ since $\delta > \Theta(\epsilon^{1-a})$. Repeating this process recursively, we would have constructed the disjoint sets $\{\mathcal{G}_1, \mathcal{G}_2, \dots, \mathcal{G}_m \}$ with $2\leq m\leq n$. Since $n$ is finite, this process will terminate in finite steps (maximum $n-1$ steps) and therefore after the final splitting, we will obtain $\delta = \Omega(\epsilon^{1-a})$ for some $a \in (0,1]$ such that $\max_{p}\{\mathbf{diam} (\mathcal{G}_p)\}  = \mathcal{O}(\epsilon^{1-a}) $.}
\end{proof}

\revise{Proposition \ref{eigenprop} describes a fundamental property of any $\mathcal{C}^2$ function that arises due to the algebraic multiplicity / (approximate) degeneracy of the eigenvalues of its Hessian at the saddle points. Note that, as a consequence of the strict-saddle property (Assumption~\textbf{A4} / Proposition~\ref{morseprop}) and Proposition~\ref{eigenprop}, we get the following necessary condition:}
\begin{align}
	\revise{\beta \geq \frac{\delta}{2}.}
\end{align}

\section{Gradient trajectories and their approximations around strict saddle point}
In this section, we analyze the behavior of the gradient descent algorithm in the vicinity of our strict saddle point, i.e., the region given by the set of points contained in $\mathcal{B}_{\epsilon}(\x^*)$. It has been already established that gradient descent converges to minimizers and almost never ends up terminating into a strict saddle point \cite{lee2016gradient}. However, the geometric structure of the region $\mathcal{B}_{\epsilon}(\x^*)$ has not been utilized completely in prior works when it comes to developing rates of escape (possibly linear). Although linear rates of divergence from a strict saddle point are provided in \cite{o2019behavior} for the Nesterov accelerated gradient method, their analysis is reserved only for quadratic functions. Intuitively, for saddle neighborhoods with sufficient curvature magnitude $\beta$ (\revise{Assumption \textbf{A4}, Proposition \ref{morseprop}}), there should exist some gradient trajectories that escape the saddle neighborhood $\mathcal{B}_{\epsilon}(\x^*)$ with linear rate every time. Moreover, these trajectories should have some dependence on their initialization $\x_{0}$. To support this intuition of a linear escape rate, we first need an understanding of the behavior of gradient flow curves in the saddle point neighborhood, following which parallels can be drawn between flow curves and gradient trajectories.

We start by formally defining the gradient descent update and the corresponding flow curve equation. For a constant step size, the gradient descent method is given by
\begin{align}
	\x_{k+1} = \x_{k} - \alpha \nabla f(\x_k),  \label{gd}
\end{align}
where $\alpha$ is the step size and we require that $\alpha \leq \frac{1}{L}$.

Next, the corresponding gradient flow curve is defined. If the step size $\alpha$ in \eqref{gd} is taken to $0$, the discrete iterate equation in index $k$ of gradient descent can be transformed into a continuous-time ODE in $t$ given by
\begin{align}
	 \frac{d \x(t)}{dt} & = - \nabla f(\x(t)), \label{ode0}
\end{align}
which is the gradient flow equation in the limit of $\alpha \to 0$ \cite{bolte2010characterizations}. Note that although $\norm{\x_{k+1} -\x_{k}}$ is $\mathcal{O}(\epsilon)$ here since both $\x_{k}$ and $\x_{k+1}$ lie inside $\mathcal{B}_{\epsilon}(\x^*)$, we still require that $\alpha \to 0$ to transform the discrete iterate update into a continuous-time ODE.

We now state the following lemma about the gradient norm $\norm{\nabla f(\x)}$ when $\x \in \mathcal{B}_{\epsilon}(\x^*)$.
\begin{lemma}\label{lem1}
For every point $\x \in \mathcal{B}_{\epsilon}(\x^*)$, the gradient $\nabla f(\x)$ will have $\mathcal{O}(\epsilon)$ magnitude.
\end{lemma}
\begin{proof}
This can be verified using Assumption \textbf{A2}:
\begin{align}
    \norm{\nabla f(\x) - \nabla f(\x^*)} & \leq L \norm{\x - \x^*}
			 \leq L \epsilon.
\end{align}
\end{proof}

This lemma is of importance since it will help us in characterizing the gradients in the ball $ \mathcal{B}_{\epsilon}(\x^*)$ in terms of the Hessian $\nabla^2 f(\x^*) $ at the saddle point, from which we will develop approximations of gradient trajectories around the saddle point.

\subsection{Intuition behind the linear time of escape} \label{odeanalysisflowcurve}
From the ODE analysis of flow curves for gradient-related methods such as those in \cite{hu2017fast,kifer1981exit}, it can be readily inferred that the gradient flow curves show hyperbolic behavior in the vicinity of saddle points. Since the discrete gradient method \eqref{gd} is the Euler discretization of the gradient flow curve ODE \eqref{ode0}, the geometric behavior of these two equations should be similar to one another with a deviation between them not more than of order $\mathcal{O}(\alpha)$ when the step size $\alpha$ is sufficiently small.\footnote{The actual deviation between the gradient flow curve and the gradient descent method after $k$ iterations tends to be on the order of $\mathcal{O}(k \alpha)$ for a fixed $\epsilon$. However, the factor $k$ can be suppressed provided the trajectories generated by the two methods do not have large exit times.} Therefore a crude analysis of flow curves should be sufficient to make approximate deductions for the discrete gradient method.

Concretely, we first define a time-varying vector $\u(t)$ that points to our iterate $\x(t)$ from the first-order strict saddle point $\x^*$. By this definition, we have that
\begin{align}
    \u(t) &= \x(t) - \x^*
    \implies   \frac{d \u(t)}{dt}  =\frac{d \x(t)}{dt}. \label{ode1}
\end{align}
Now, computing the norm squared of $\u(t)$, differentiating it with respect to $t$ and using \eqref{ode0}, we get
\begin{align}
    \norm{\u(t)}^2 & = \norm{\x(t)-\x^*}^2  \\
    \implies  \frac{d  \norm{\u(t)}^2}{dt}  & = 2\langle (\x(t) - \x^*),  - \nabla f(\x(t)) \rangle. \label{odeflowcurve}
\end{align}
Next, let the gradient flow curve enter $\mathcal{B}_{\epsilon}(\x^*)$ ball at time $t=0$ and exit this ball at time $t=T$. Geometrically, the inner product defined in \eqref{odeflowcurve} is negative at the entry point of the ball $\mathcal{B}_{\epsilon}(\x^*)$ (i.e., vectors $(\x(0) - \x^*)$ and $- \nabla f(\x(0))$ form an obtuse angle), becomes equal to $0$ at some point $\x_{critical}$ inside this ball and is positive at the exit point (i.e., vectors $(\x(T) - \x^*)$ and $- \nabla f(\x(T))$ form an acute angle).

Using Taylor's expansion around $\x^*$ along the direction $\x(t)-\x^*$, we can write $\nabla f(\x(t))$ in the following manner:
\begin{align}
    \nabla f(\x(t)) = \nabla f(\x^*) + \int_{p=0}^{p=1}\nabla^2 f(\x^* + p\u(t)) \u(t) dp. \label{taylorl1}
\end{align}
If $\norm{\u(t)}$ is sufficiently small or is of order $\mathcal{O}(\epsilon)$, we can approximate $\nabla^2 f(\x^* + p\u(t)) \approx \nabla^2 f(\x^*) $. After substituting this approximation in \eqref{taylorl1} we obtain
\begin{align}
    \nabla f(\x(t)) \approx \nabla^2 f(\x^*)\u(t). \label{taylorapproxl1}
\end{align}
Using this result in \eqref{odeflowcurve} yields
\begin{align}
   	 \frac{d  \norm{\u(t)}^2}{dt}  & = 2\langle (\x(t) - \x^*),  - \nabla f(\x(t)) \rangle \approx -2 \langle \u(t), \nabla^2 f(\x^*)\u(t) \rangle. \label{ode3}
\end{align}
Also using \eqref{ode0} and \eqref{taylorapproxl1} we get that
\begin{align}
   \frac{d\u(t)}{dt}  & =    \frac{d\x(t)}{dt} \approx -\nabla^2 f(\x^*)\u(t). \label{ode2}
\end{align}
Now consider the case where Assumptions \textbf{A1} to \revise{\textbf{A4}} are satisfied. Since the eigenvalues of $\nabla^2 f(\x^*)$ are both positive and negative, the approximate ODE \eqref{ode2} will have the following  solution:
\begin{align}
  \u(t) = \sum_{i=1}^{n} c_i \v_i(0) e^{-\lambda_{i}(0)t},
\end{align}
where $(\lambda_{i}(0), \v_i(0))$ represents the $i^{th}$ eigenvalue--eigenvector pair for the Hessian $\nabla^2 f(\x^*)$ and $c_i$ are non-negative constants that depend on the initialization $\u(0)$. \revise{(Here, the non-negativity of $c_i$'s can be assumed without loss of generality because the sign of the eigenvectors can be chosen arbitrarily.)}

From this equation it is clear that we have a solution that is exponential in $t$. Moreover from the approximate ODE \eqref{ode3}, it is evident that a hyperbolic curve is generated with an exponential rate of change. Therefore, for any initialization, i.e., for any choice of constants $c_i$, $\norm{\u(t)}^2$ \revise{eventually increases at an exponential rate}, thereby giving a linear escape rate for $\x(t)$ from the region $\mathcal{B}_{\epsilon}(\x^*)$ provided $c_i \neq 0$ corresponding to at least one of the negative eigenvalues.

However, the approximation $\nabla^2 f(\x^* + p\u(t)) \approx \nabla^2 f(\x^*) $ fails to capture the first-order perturbation terms in the Hessian $\nabla^2 f(\x^* + p\u(t))$. Given a sufficiently small saddle neighborhood $\mathcal{B}_{\epsilon}(\x^*)$, for any $\x \in \mathcal{B}_{\epsilon}(\x^*)$, the eigenvalues and eigenvectors of the Hessian $\nabla^2 f(\x)$ can have $\mathcal{O}(\epsilon)$ variations with respect to those of the Hessian $\nabla^2 f(\x^*)$. Taking this $\mathcal{O}(\epsilon)$ perturbation into account complicates the gradient flow curve analysis inside the $\mathcal{B}_{\epsilon}(\x^*)$ ball,\footnote{This is formally taken into account in subsequent sections using matrix perturbation theory.} which otherwise is straightforward from \eqref{ode2}. Moreover, for all practical purposes, we cannot take our step size $\alpha \to 0$ for the sake of using ODE analysis. Choosing arbitrarily small step sizes causes the number of iterations needed to escape from the ball $\mathcal{B}_{\epsilon}(\x^*)$ increase to infinity. Therefore a discrete gradient trajectory analysis using matrix perturbation theory becomes an absolute necessity to obtain trajectories (or approximate trajectories) with linear exit time from $\mathcal{B}_{\epsilon}(\x^*)$.

\subsection{Warm-up: Rudimentary analysis of the exit time for discrete gradient trajectories}\label{crudeanalysissection}
The intuition developed as part of the ODE-based analysis \emph{suggests} linear time escape of discrete gradient trajectories from strict-saddle neighborhoods. We now present a rudimentary analysis of the gradient descent method that uses elementary facts about first-order methods, as opposed to matrix perturbation theory, to derive a bound on the exit time of the gradient descent method from the saddle neighborhood $\mathcal{B}_{\epsilon}(\x^*)$. The purpose of this analysis is twofold. First, it shows that (discrete) gradient-descent trajectories can indeed escape strict-saddle neighborhoods in linear time. Second, it highlights the limitations of existing analytical techniques in deriving linear escape rates for discrete gradient trajectories, thereby motivating the need for the matrix perturbation-based analysis of gradient trajectories in the next section for derivation of a linear escape rate from $\mathcal{B}_{\epsilon}(\x^*)$. 

Note that the analysis in this section requires only a relaxed version of Assumption \textbf{A1} on the function $f(\cdot)$, namely, it is twice continuously differentiable: $f \in \mathcal{C}^2$. But the remaining assumptions (\revise{Assumptions~\textbf{A2}--\textbf{A4}}) stay the same. Now consider the following that follows from the gradient descent iteration: 
\begin{align}
    \x_{k+1} - \x^* &= \x_{k} - \x^* - \alpha \nabla f(\x_k) \\
   &= \x_{k} - \x^* - \alpha \nabla^2 f(\x^*)(\x_k -\x^*)  - \alpha (\nabla f(\x_k)-\nabla^2 f(\x^*)(\x_k -\x^*)) \\
   &= (\mathbf{I}- \alpha \nabla^2 f(\x^*))(\x_k -\x^*) - \alpha r(\x_k), \label{crudeanalysis1}
\end{align}
where $r(\x_k) = (\nabla f(\x_k)-\nabla^2 f(\x^*)(\x_k -\x^*))$. Using the Hessian Lipschitz continuity of $f(\cdot)$ and the fact that $\nabla f(\x_k) = \nabla f(\x^*) + \int_{p=0}^{p=1}\nabla^2 f(\x^* + p(\x_k -\x^*)) (\x_k -\x^*) dp $ since $f \in \mathcal{C}^2$, we get that
\begin{align}
    \norm{r(\x_k)} &= \norm{\int_{p=0}^{p=1}\nabla^2 f(\x^* + p(\x_k -\x^*)) (\x_k -\x^*) dp-\nabla^2 f(\x^*)(\x_k -\x^*)} \\
& \leq \bigg(\int_{p=0}^{p=1} \norm{\nabla^2 f(\x^* + p(\x_k -\x^*))  -\nabla^2 f(\x^*)}dp \bigg)\norm{(\x_k -\x^*)} \\
& \leq \frac{M\norm{(\x_k -\x^*)}^2}{2}.
\end{align}
Thus, $\norm{r(\x_k)} \leq \frac{M\epsilon^2}{2}$ whenever $ \x_k \in \mathcal{B}_{\epsilon}(\x^*)$. Inducting \eqref{crudeanalysis1} up to $k=0$ yields:
\begin{align}
    \x_{k+1} - \x^* = (\mathbf{I}- \alpha \nabla^2 f(\x^*))^{k+1}(\x_0 -\x^*) - \alpha \sum\limits_{i=0}^{k}  (\mathbf{I}- \alpha \nabla^2 f(\x^*))^{k-i}r(\x_i). \label{crudeanalysis2}
\end{align}

Next, in order to analyze the worst case bounds on the exit time, assume that the unstable subspace of $\nabla^2 f(\x^*)$ has dimension $1$, i.e., $\lambda_j > 0   $ for all $j \in \{1,2,\ldots,n-1\}$ and $ \lambda_n < 0 $, where $ \lambda_j$ is the $j^{th}$ eigenvalue of $\nabla^2 f(\x^*)$. Also let $ \v_n$ be 
\revise{an eigenvector of $\nabla^2 f(\x^*)$ of unit norm} corresponding to the eigenvalue $\lambda_n$, where $ \lambda_n < - \beta $ from Assumption \revise{\textbf{A4}}. Since divergence can happen only from the unstable subspace, our assumption on $\nabla^2 f(\x^*)$ will leave only a single direction of escape, i.e. along $\v_n$, for the gradient trajectories. Moreover since both $\v_n$ and $-\v_n$ will be the eigenvectors of $\nabla^2 f(\x^*)$, hence without loss of generality let us assume that $ \langle\v_n, (\x_0 - \x^*)\rangle \geq 0 $, where $\x_0 \in \bar{\mathcal{B}}_{\epsilon}(\x^*) \backslash \mathcal{B}_{\epsilon}(\x^*)$, and we are required to find the exit time $K_{exit}$ that satisfies
\begin{align}
    K_{exit} = \inf_{k>0}\{k \vert \norm{\x_k - \x^*} > \epsilon\}.
    \label{eq-Kexit-1}
\end{align}
\revise{As we show in Lemma \ref{lemma-exit-times-equivalent} in Appendix~\ref{app:equivalence}, this} is equivalent to the following condition:
\begin{align}
     K_{exit} = \inf_{k>0}\{k \vert  \langle\v_n, (\x_k - \x^*)\rangle > \gamma_k \epsilon\}, \label{exittimecrudeanalysis}
\end{align}
where $ \gamma_k = \frac{ \langle\v_n, (\x_k - \x^*)\rangle}{\norm{\x_k - \x^*}}$ and \revise{we have assumed for the sake of the crude analysis that} $\gamma_k \in (0,1]$ \revise{for every $k$}. Now, taking the inner product of $\v_n$ with $(\x_{k+1} - \x^*)$ in \eqref{crudeanalysis1}, 
 and using the Hessian Lipschitz continuity and $\norm{r(\x_k)} \leq \frac{M \epsilon^2}{2} $, we get:
\begin{align}
   \langle\v_n, \x_{k+1} - \x^* \rangle &= \langle\v_n,(\mathbf{I}- \alpha \nabla^2 f(\x^*))(\x_k -\x^*)\rangle - \alpha   \langle\v_n,r(\x_k)\rangle \geq (1+ \alpha \beta)\langle\v_n, \x_{k} - \x^* \rangle - \frac{\alpha M \epsilon^2}{2}, \label{crudeanalysis3}
\end{align}
where we have used the substitution $ (\mathbf{I}- \alpha \nabla^2 f(\x^*)) = \sum\limits_{j=1}^{n}(1- \alpha \lambda_j)\v_j \v_j^T$. To show divergence from $\x^*$, it then suffices to show that for some $ \rho \in (0,1)$ we have
\begin{align}
    (1+ \alpha \beta)\langle\v_n, \x_{k} - \x^* \rangle - \frac{\alpha M \epsilon^2}{2} \geq  (1+ \rho\alpha \beta)\langle\v_n, \x_{k} - \x^* \rangle \label{crudeanalysis4}
\end{align}
hold for all $k$, \revise{which will then imply that} $ \langle\v_n, \x_{k} - \x^* \rangle$ is strictly monotonically increasing with $k$.\footnote{In general, we do not need the monotonicity condition for all $k$ but only after a sufficiently large $k$ that is smaller than the exit time. Such trajectories will also have linear exit times as proved in a subsequent counterexample.} Further simplifying \eqref{crudeanalysis4} we get the condition
\begin{align}
     \beta(1-\rho)\langle\v_n, \x_{k} - \x^* \rangle \geq \frac{ M \epsilon^2}{2}, \label{crudeanalysis4.2}
\end{align}
which should hold for all $k$. A sufficient boundary condition for this inequality to hold is:
\begin{align}
    \langle\v_n, \x_{0} - \x^* \rangle \geq \frac{ M \epsilon^2}{2\beta(1-\rho)}. \label{crudesufficientcondition}
\end{align}
Now if the boundary condition \eqref{crudesufficientcondition} holds then from \eqref{crudeanalysis3} and \eqref{crudeanalysis4} we have:
\begin{align}
    \langle\v_n, \x_{k} - \x^* \rangle \geq (1+\rho \alpha \beta)^k\langle\v_n, \x_{0} - \x^* \rangle \geq  (1+\rho \alpha \beta)^k\frac{ M \epsilon^2}{2\beta(1-\rho)}.
\end{align}
Then using \eqref{exittimecrudeanalysis}, exit from the ball $\mathcal{B}_{\epsilon}(\x^*)$ can be guaranteed by setting the following condition:
\begin{align}
    \langle\v_n, \x_{k} - \x^* \rangle  \geq  (1+\rho \alpha \beta)^k\frac{ M \epsilon^2}{2\beta(1-\rho)} &> \gamma_k \epsilon \\
    \Longleftrightarrow \quad k & \geq \frac{\log(\frac{2\gamma_k\beta(1-\rho)}{M \epsilon})}{\log (1+\rho \alpha \beta)},
\end{align}
which implies $K_{exit} \leq \frac{\log\left(\frac{2\gamma_{K_{exit}}\beta(1-\rho)}{M \epsilon}\right)}{\log (1+\rho \alpha \beta)}$ as long as the sufficient condition \eqref{crudesufficientcondition} is satisfied.

The preceding rudimentary analysis guarantees a linear exit time bound for the gradient descent method under the sufficient boundary condition \eqref{crudesufficientcondition}. But the resulting exit time bound is loose due to its dependency on the unknown factors $\gamma_{K_{exit}}$ and $\rho$, where $\gamma_{K_{exit}}$ could be arbitrarily small and the presence of $\rho$ in the boundary condition makes this analysis more restrictive than the matrix perturbation-based analysis presented in Section \ref{exactexittimeanalysissection}. Also, the exit time analysis in this section does not bring out the dependence of boundary conditions and exit time bound on the problem dimension, conditioning of the neighborhood, and spectral gap, etc. Such dependencies are captured in the analysis of Section \ref{exactexittimeanalysissection} and Table \ref{table:1} in Section \ref{ssec:comparison} summarizes the corresponding differences between the two analytical approaches. More importantly this analysis guarantees a linear exit time bound only for those trajectories starting at $\x_0$ that satisfy the monotonicity property implied by \eqref{crudeanalysis3} and \eqref{crudeanalysis4}. That is, it does not capture the trajectories for which $ \langle\v_n, \x_{k} - \x^* \rangle $ does not increase monotonically with $k$. A simple counterexample to the need for this monotonicity property for derivation of a linear exit time bound can be easily constructed. We refer the reader to Appendix \ref{Appendix E} for one such counterexample. This implies there exist gradient trajectories that can exit in a linear time while violating the monotonicity condition, thereby illustrating that the rudimentary exit time analysis does not capture all the trajectories with linear exit times.\looseness=-1

\begin{remark} \label{remarkmeasure}
Note that the sufficient condition of $  \langle\v_n, \x_{0} - \x^* \rangle \geq \frac{ M \epsilon^2}{2\beta(1-\rho)}$ from \eqref{crudesufficientcondition} guarantees linear exit time gradient trajectories. Moreover this condition makes sure that such trajectories do not have zero measure since the set of initialization given by $ \{ \x_0 \hspace{0.1cm} \vert \hspace{0.1cm}  \langle\v_n, \x_{0} - \x^* \rangle \geq \frac{ M \epsilon^2}{2\beta(1-\rho)}\}$ has positive measure for sufficiently small $\epsilon$.
\end{remark}

In summary, to analyze the complete set of gradient trajectories around the saddle point that escape in linear time and develop a precise exit time bound we need more than the class of twice-differentiable functions; hence the need to work with analytic functions.\footnote{In order to get a highly precise bound on exit time, we need the best possible first-order approximations of gradient trajectories, which can only be obtained for analytic functions. Therefore even the class of $\mathcal{C}^{\infty}$ functions is not sufficient for our analysis; see also the discussion in Remark~\ref{remark:Davis.Kahan} in Section \ref{degenmat} in this regard.} Note that many optimization and learning problems, such as quadratic functions and deep neural networks with smooth activation functions,  satisfy real analyticity in some neighborhoods of stationary points, if not over the entire domain.

\subsection{\revise{An informal statement of the main result}}%
\revise{In this section, we provide an informal statement of the main result of this paper as well as a brief discussion of the implications of this result.}

\begin{theorem}[Informal Main Result]\label{informalthm}
\revise{Under Assumptions \textbf{A1}--\textbf{A4}, the approximate trajectories of the gradient descent method with step size $\alpha = \frac{1}{L}$, when initialized on the boundary of some $\epsilon$ neighborhood of a strict saddle point $\x^*$ of $f(\cdot)$, where $\epsilon < \min\bigg\{ \frac{2 \beta}{M}, \Omega\bigg(\frac{\delta}{n^2}\bigg)\bigg\}$ and $\epsilon \ll 1$, can exit this neighborhood in approximately linear time, i.e., $K_{exit} \lessapprox \mathcal{O}\bigg(\log\bigg( \frac{\delta}{\epsilon n} \bigg)\bigg)$, where $K_{exit}$ is the exit time for the approximate trajectory, $n$ is the problem dimension and $\delta$ is the eigen gap from Proposition \ref{eigenprop}. However, this linear exit time bound holds only if the initial radial vector $\u_0= \x_0 -\x^*$ is not orthogonal to the unstable subspace of $\nabla^2 f(\x^*)$ and subtends some non-zero angle with the stable subspace of $\nabla^2 f(\x^*)$. In particular, the cosine square of the angle between the initial radial vector and the unstable subspace of $\nabla^2 f(\x^*)$ must be at least of the order $\Omega\bigg(\frac{\epsilon n}{\delta}\bigg)$, where this cosine square is referred to as the unstable subspace projection value.}
\end{theorem}
\revise{A formal statement of this result, which includes precise characterizations of the approximate trajectory, exit time, and the bounds on $\epsilon$ as well as the necessary initial unstable subspace projections, is provided in Theorem~\ref{theorem2}. We also refer the reader to Figure \ref{familyofcurves} for a concrete intuition of the angle between the initial radial vector and the unstable subspace of $\nabla^2 f(\x^*)$ as well as its relation to the unstable subspace projection.}

\revise{We now briefly summarize the implications of this main result, while additional discussion is provided after Theorem~\ref{theorem2}. For a function $f(\cdot)$ satisfying Assumptions \textbf{A1}--\textbf{A4}, let the gradient descent method with step size $\alpha = \frac{1}{L}$ be initialized on the boundary of some $\epsilon$ neighborhood of a strict saddle point $\x^*$ of $f(\cdot)$ such that the initial radial vector $\u_0= \x_0 -\x^*$ subtends some angle with the unstable subspace of $\nabla^2 f(\x^*)$ that is not equal to $\frac{\pi}{2}$. Then we have the following statements:
\begin{itemize}
    \item [\textbf{S1.}] There exists some lower bound on the cosine square of this angle (termed as the `sufficient condition') for which the approximate trajectories of the gradient descent method will exit the saddle neighborhood in linear time.
    \item [\textbf{S2.}] Also, there exists a strict lower bound on the cosine square of this angle (termed as the `necessary condition') that is of the order $\Omega\bigg(\frac{\epsilon n}{\delta}\bigg)$. If the cosine square of the angle between the initial radial vector and the unstable subspace of $\nabla^2 f(\x^*)$ is smaller than $\Theta\bigg(\frac{\epsilon n}{\delta}\bigg)$, the approximate trajectories of the gradient descent method can never exit the saddle neighborhood in linear time.
\end{itemize}
This work rigorously establishes Statement \textbf{S2} and also shows that Statement \textbf{S1} is not vacuous (cf.~Section ~\ref{nonvacuousclaim} in Appendix~\ref{Appendix D}). Note that a rigorous characterization of the lower bound in Statement \textbf{S1} requires a more sophisticated proof machinery, which has been pursued in our follow-up work \cite{dixit2022boundary}.}
\begin{remark}
\revise{A fast exit time in terms of the scaling with $\frac{1}{\epsilon}$ in and of itself might not preclude the gradient descent method from converging super slowly in the worst case. The carefully constructed function with cascaded saddles in \cite{du2017gradient}, in particular, is a prime example of this behavior, as the gradient descent method takes an exponentially---{in dimension $n$}---large time in the worst case to escape the cascaded saddles and converge to a local minimum for this function. However, the particular class of functions within the family of Morse functions being considered in this work excludes the construction in  \cite{du2017gradient}. Going further, we have established in a follow-up work \cite{dixit2022boundary} that the time to escape cascaded saddles and reach a second-order stationary point for functions in this class does not scale exponentially in the dimension for a simple variant of the gradient descent method.}
\end{remark}

\subsection{Brief overview of results and proof sketch for the linear exit time bound}
Our matrix perturbation-based analysis utilizes the standard gradient-descent method \eqref{gd} in the saddle neighborhood $\mathcal{B}_{\epsilon}(\x^*)$. Since we are interested in developing analysis suited only for the region $\mathcal{B}_{\epsilon}(\x^*)$, we assume that initially our iterate $\x_{0}$ sits on the boundary of $\mathcal{\bar{B}}_{\epsilon}(\x^*)$. We then follow the given sequence of steps in order to obtain linear exit time bound for approximations of gradient descent trajectories around a saddle point.

\begin{enumerate}
\item Starting with Lemma \ref{lem2} we show that the region $\mathcal{B}_{\epsilon}(\x^*)$ around the strict saddle point $\x^*$ is comprised of a stable and an unstable subspace, which are orthogonal to one another.
\item Next, for any $\x \in \mathcal{B}_{\epsilon}(\x^*)$ we write $ \nabla f(\x)$ in terms of the radial vector $\u = \x-\x^*$ as $\nabla f(\x) =  \bigg(\displaystyle\int_{p=0}^{p=1}\nabla^2 f(\x^* + p\u)  dp\bigg) \u$.
\item Then in Lemma \ref{lem4} using matrix perturbation theory we express the Hessian $\nabla^2 f(\x)$ at $\x = \x^* + p \u$, where $\x \in \mathcal{B}_{\epsilon}(\x^*)$, $p \in [0,1]$, and $\norm{\u} \leq \epsilon$ in terms of a perturbation of $\nabla^2 f(\x^*)$, as $$\nabla^2 f(\x^* + p\u) = \nabla^2 f(\x^*) + \D(\x),$$ with the perturbation matrix $\D(\x)$ bounded as $$ \norm{\D(\x)}\leq M p \epsilon.$$
\item We iterate the Gradient descent method in terms of the radial vector $\u_k$ as follows:
    \begin{align*}
   \u_{k+1} &=\x_{k} - \x^* - \alpha \nabla f(\x_k) = \bigg(\mathbf{I} - \alpha \int_{0}^{1} \nabla^2 f(\x^* + p \u_k)dp \bigg) \u_k  \\
  \implies \u_{k+1}  &= \bigg(\mathbf{I} - \alpha \nabla^2 f(\x^*) -\underbrace{\alpha \int_{0}^{1} \D(\x^* + p \u_k)dp}_{\R(\u_k) =\mathcal{O}(\epsilon)} \bigg) \u_k
    \end{align*}
    where $ \norm{\R(\u_k)} = \norm{ {\alpha \int_{0}^{1} \D(\x^* + p \u_k)dp}} = \mathcal{O}(\epsilon)$ from the last step. Using this radial vector update in Lemma \ref{lem5}, we induct the above recursion up to initialization $\u_0$ and obtain the exact trajectory expression: $$ \u_{K+1}  = \Pi_{k=0}^{K}\bigg(\mathbf{I} - \alpha \nabla^2 f(\x^*) - \R(\u_k)\bigg) \u_0.$$
\item In Lemma \ref{lem6} we expand the product of the $K+1$ non-commuting matrices from the last step up to first order as follows: $$\tilde{\u}_{K+1} := \Pi_{k=0}^{K} \A_k \u_0- \sum\limits_{r=0}^{K}(\Pi_{k=r+1}^{K}\A_r\R(\u_r)\Pi_{k=0}^{r-1}\A_r)\u_0,$$ where $\tilde{\u}_{K+1} \approx {\u}_{K+1}$ and $\A_k :=  \mathbf{I} - \alpha \nabla^2 f(\x^*) $ for all $k$ \revise{in the case of gradient descent}. This is the most  crucial step in the analysis since we obtain the approximate trajectory $\{\tilde{\u}_K\}$ in this step.\footnote{\revise{Even though $\A_k$ is constant for the gradient-descent iteration \eqref{gd}, we have purposefully not removed its subscript $k$ since it may not be constant for a general dynamical system. Consider, for instance, gradient descent with variable step size $\alpha_k$ instead of constant step size $\alpha$ and we then have $\A_k = \mathbf{I} - \alpha_k \nabla^2 f(\x^*)$. Hence, with the subscript $k$ intact, the expression for the approximate trajectory $\{\tilde{\u}_K\}$ can be easily adapted to a general class of first-order methods.}}
\item The approximate trajectory $\{\tilde{\u}_K\}$ obtained above cannot be {uniquely} determined since it is a function of the eigenvalues of the Hessian $\nabla^2{f(\x^*)}$, which are known only up to an interval. Therefore in Lemma \ref{lem7} we obtain a parametrized family of approximate trajectories for a fixed $\u_0$, denoted by $\{\tilde{\u}_K^{\tau}\}$, where the parameter $\tau \in \mathbb{R}$ varies with variations in the eigenvalues of the Hessian $\nabla^2{f(\x^*)}$. Next, we construct the \textit{minimal} approximate trajectory from this family, defined as one that stays closest to $\x^*$ for each $K$ and show that this \textit{minimal} approximate trajectory has the maximum exit time among all approximate trajectories.
\item In Theorem \ref{theorem1} we obtain the closed form expression of the normalized radial distance for the \textit{minimal} approximate trajectory given by $\Psi(K)$ where $ \epsilon^2 \Psi(K) \leq \inf_{\tau} \norm{\tilde{\u}_K^{\tau}}^2 < \epsilon^2$.
\item Finally in Theorem \ref{theorem2} we obtain the smallest upper bound on $K$ of the order $\mathcal{O}(\log(\epsilon^{-1}))$ that satisfies the condition $\Psi(K) > 1$ which will imply $\epsilon^2 < \epsilon^2 \Psi(K) \leq \inf_{\tau} \norm{\tilde{\u}_K^{\tau}}^2  $. This condition gives the linear exit time bound from the saddle neighborhood. We then derive any necessary conditions on $\x_0$ for guaranteeing this linear exit time.
\end{enumerate}

Before formally beginning our analysis of discrete gradient trajectories, we state the following lemma that will be utilized frequently in our analysis.  	
\begin{lemma}\label{lem2}
  For any point $\x \in \mathcal{B}_{\epsilon}(\x^*)$, the vector $\u$ given by $\u = \x - \x^*$ belongs to a vector space $\mathcal{E}$ that is comprised of a stable subspace $\mathcal{E}_{S}$ (subspace corresponding to contraction dynamics) and an unstable subspace $\mathcal{E}_{US}$ (subspace corresponding to expansive dynamics). Formally, this can be written as
  \begin{align}
  	\mathcal{E} = \mathcal{E}_{S}  \bigoplus \mathcal{E}_{US}, \nonumber
  \end{align}
  where $\bigoplus$ denotes the direct sum of two spaces.
\end{lemma}
\begin{proof}
The eigenvalues of the Hessian $\nabla^2 f(\x^*)$ are both positive and negative. Without loss of generality, these can be classified into two sets of stable and unstable eigenvalues with the stable set comprising positive eigenvalues and the unstable set having negative eigenvalues. Then the corresponding subspaces can be written as
\begin{align}
  \mathcal{E}_{S} &= span\{\v_i | \lambda_{i}(\nabla^2 f(\x^*)) > 0\}, \ \text{and} \\
  		\mathcal{E}_{US} &= span\{\v_j  | \lambda_{j}(\nabla^2 f(\x^*)) < 0\},
\end{align}
where $\lambda_{i}(\nabla^2 f(\x^*)), \v_i$ represent the $i^{th}$ eigenvalue-eigenvector pair. Since these subspaces are orthogonal and span the complete space $\mathcal{E}=\mathbb{R}^n$, any vector $\u = \x - \x^*$ is spanned by these subspaces. Next, we define the two index sets $ 	\mathcal{N}_{S} = \{i | \lambda_{i}(\nabla^2 f(\x^*)) > 0\} $ and $\mathcal{N}_{US} =\{ j | \lambda_{j}(\nabla^2 f(\x^*)) < 0\}$ for the two subspaces. Since these subspaces are orthogonal, their index sets are disjoint.
\end{proof}

\subsection{Analysis of discrete gradient trajectories using matrix perturbation theory}\label{exactexittimeanalysissection}
Now that we have established all the necessary preliminaries, we can move on to develop approximate bounds on the escape time from the region $\mathcal{B}_{\epsilon}(\x^*)$ for gradient descent. From here onward we restrict ourselves to discrete time iterates denoted by subscripts $k$ and the entire analysis is carried out in discrete time. Also, we assume that Assumptions \textbf{A1} to \revise{\textbf{A4}} hold along with the additional condition of $m=n$ in \revise{Proposition \ref{eigenprop}}, i.e., there are no degenerate eigenvalues. Section \ref{degenmat} after Lemma \ref{lem4} discusses the analysis for the degenerate eigenvalues, i.e., the case when $m \neq n$ in \revise{Proposition \ref{eigenprop}}. In there, we show that the analysis for the degenerate case is very straightforward and easy to extend from the non-degenerate analysis. \revise{It should also be noted that instead of analyzing exact trajectories, we analyze from here onward the first-order approximations of the exact trajectories, where the approximation error is sufficiently small. The presence of the higher-order terms ($\mathcal{O}(\epsilon^2)$ terms) in the forthcoming analysis accounts for the approximation in our analysis, and things are proved about trajectories and perturbations up to the first order in $\epsilon$.} To summarize our next set of steps, we begin with a lemma that characterizes the approximate Hessian behavior in the region $\mathcal{B}_{\epsilon}(\x^*)$, followed by a lemma that expresses $\x_{k}$ for any $k\geq0$ approximately in terms of $\x_{0}$ and a theorem that characterizes an approximate lower bound on the distance of $\x_{k}$ from $\x^*$.

\begin{lemma}\label{lem4}
Let $r_j(\u)$ be a function of the vector $\u$ defined as $r_j(\u) = \norm{\bigg(\frac{d^j }{dw^j}\nabla^2 f(\x^*+w\u)\bigg\vert_{w=0}\bigg)}_2$ and $\epsilon>0$ be a constant that satisfies the necessary condition of $\epsilon < \inf_{{\norm{\u}=1}}\bigg(\limsup_{j\to \infty} \sqrt[j]{\frac{r_j(\u)}{j!}}\bigg)^{-1}$. Then for any $\x_k\in\mathcal{B}_{\epsilon}(\x^*)$ such that $\x_{k} = \x^* + p \u_{k}$ with $0<p\leq1$, the Hessian $\nabla^2 f(\x_k)$ is given by
\begin{align}
    \nabla^2 f(\x_k) &= \nabla^2 f(\x^*) + p \norm{\u_{k}} \H(\hat{\u}_k) + \mathcal{O}(\epsilon^2),
\end{align}
where $\hat{\u}_k= \frac{\u_k}{\norm{\u_k}}$ and we have that
\begin{align}
    \H(\hat{\u}_k) &=  \sum_{i=1}^{n} \bigg(\langle \v_i(0), \H(\hat{\u}_k)\v_i(0)\rangle \v_i(0)\v_i(0)^{T} + \lambda_i(0)(\sum_{l \neq i} \frac{\langle \v_l(0), \H(\hat{\u}_k)\v_i(0)\rangle}{\lambda_{i}(0)-\lambda_{l}(0)}\v_l(0))\v_i(0)^{T} \nonumber \\ &+ \lambda_i(0)\v_i(0)(\sum_{l \neq i} \frac{\langle \v_l(0), \H(\hat{\u}_k)\v_i(0)\rangle}{\lambda_{i}(0)-\lambda_{l}(0)}\v_l(0))^{T}\bigg)   \label{perturbation}
\end{align}      	
    with $\lambda_{i}(0), \v_i(0)$ being the $i^{th}$ eigenvalue--eigenvector pair of the Hessian $\nabla^2 f(\x^*)$.
\end{lemma}

The proof of this lemma is given in Appendix \ref{Appendix A}. \revise{Note that the expression for $\H(\hat{\u}_k)$ in the lemma statement is more of a property rather than a definition, where $\norm{\H(\hat{\u}_k)}_2 $ is bounded. However, it may not be the case that $\H(\hat{\u}_k) = \mathcal{O}(\epsilon)$. In particular, we have the following bound from inequality \eqref{reviewreporteqnew} in Appendix \ref{Appendix B}:
\begin{align}
 \norm{\H(\hat{\u}_k)}_2 &\leq  M + \mathcal{O}(\epsilon),
\end{align}
which suggests that $\norm{\H(\hat{\u}_k)}_2 $ could even be a constant-order term; see Appendix \ref{Appendix B} for further details.}
\begin{remark}
The condition $\epsilon < \inf_{{\norm{\u}=1}}\bigg(\limsup_{j\to \infty} \sqrt[j]{\frac{r_j(\u)}{j!}}\bigg)^{-1}$ is necessary but may not be sufficient to guarantee this lemma's result. Since evaluating the radius of convergence for an expansion generated by the Rayleigh--Schr\"{o}dinger perturbation analysis is beyond the scope of this work, we only put forth this necessary condition here.
\end{remark}

\begin{remark}
Note that the quantity $\inf_{{\norm{\u}=1}}\bigg(\limsup_{j\to \infty} \sqrt[j]{\frac{r_j(\u)}{j!}}\bigg)^{-1} $ is exactly equal to the radius of convergence for the Taylor series expansion of the matrix $\nabla^2 f(\x^* + w \u)$ about $w>0$ for all $\{\u : \norm{\u}_2 = 1\}$, which is strictly positive due to the analytic nature of $f(\cdot)$. A proof of this claim is given in Appendix~\ref{Appendix A}.
\end{remark}

\subsubsection{Statement about the generality of Lemma \ref{lem4}}\label{degenmat}
It should be noted that while obtaining \eqref{perturbation}, we assumed a minimum gap of $\delta$ between any two eigenvalues of the Hessian $\nabla^2 f(\x^*)$. However, we can have many groups of equal or almost similar eigenvalues from \revise{Proposition \ref{eigenprop}}; this creates singular terms in the coefficient denominators of first-order eigenvector corrections in \eqref{perturbation}. This can be solved easily from the degenerate matrix perturbation theory, which extends the results of Rayleigh--Schr\"{o}dinger theory. From that we obtain the following new first-order correction term in place of \eqref{subsectioncorrectionterm} in the proof of the lemma for the $i^{th}$ eigenvector $\tilde{\v}_i(w)$:
\begin{align}
  	\frac{d}{dw}(\tilde{\v}_i(w))\bigg\vert_{w=0} & = \sum_{l \not \in \mathcal{G}_p} \frac{\langle \tilde{\v}_l(0), \H(\hat{\u}_k)\tilde{\v}_i(0)\rangle}{\lambda_{i}(0)-\lambda_{l}(0)}\tilde{\v}_l(0), \label{degenerate}
\end{align}
where the corresponding $i^{th}$ unperturbed eigenvalue $\lambda_{i}(0)$ belongs to the set $\mathcal{G}_{p}$. Also note that we have a new basis of eigenvectors $\tilde{\v}_i$ instead of ${\v}_i$, which resolves the degeneracy issue within the groups of similar eigenvalues. This change of basis can always be done since there are infinitely many solutions to the eigenvectors belonging to the degenerate subspaces. More importantly, we are never required to compute these eigenvectors explicitly in our analysis. To get a detailed understanding of the degenerate matrix perturbation theory, the reader can refer to \cite{MBT,MBT1}.

Therefore for the case with degenerate eigenvalue sets, the analysis will remain the same, but with fewer first-order perturbation terms (\eqref{degenerate} has $n-\abs{\mathcal{G}_{p}}$ orthogonal terms in the summation instead of the $n-1$ orthogonal terms that appear in \eqref{subsectioncorrectionterm}). Now, these fewer $\mathcal{O}(\epsilon)$ terms in \eqref{degenerate} will result in weaker first-order perturbations on the distance $\norm{\x_{k}-\x^*}$ when compared to that from \eqref{subsectioncorrectionterm}. In a subsequent lemma (Lemma \ref{lem7}), it will be established that the worst-case trajectory is obtained when the first-order perturbation terms are used to minimize $\norm{\x_{k}-\x^*}$ for every $k$. This worst-case trajectory stays inside the ball $\mathcal{B}_{\epsilon}(\x^*)$ for the maximum number of iterations. For the case of degenerate eigenvalues, fewer first-order terms from \eqref{degenerate} means a weaker perturbation effect over $\norm{\x_{k}-\x^*}$, which implies that $\norm{\x_{k}-\x^*}$ cannot be minimized completely. This is in contrast to the case of \eqref{subsectioncorrectionterm} which has more first-order terms ($n-1$) and hence a stronger perturbation effect over $\norm{\x_{k}-\x^*}$. Now, a stronger perturbation can be used to contain the worst-case trajectory inside $\mathcal{B}_{\epsilon}(\x^*)$ for a longer duration (part of the proof for Lemma \ref{lem7}). As a consequence, the worst-case trajectory from the non-degenerate case will have a larger exit time compared to that of the degenerate case. Therefore, we are not required to perform the analysis for the degenerate case because the worst-case performance in terms of exit time is captured in the current analysis for the non-degenerate case.

\begin{remark}\label{remark:Davis.Kahan}
It is worth noting here that the exit time analysis in this work could have been carried out using the Davis--Kahan theorem \cite{davis1970rotation}. Such an analysis would have required the function $f(\cdot)$ to only be $\mathcal{C}^2$, as opposed to analytic, but it would have necessitated the eigensubspaces of the Hessian $\nabla^2 f(\x^*)$ to be non-degenerate. However, non-degeneracy of the eigensubspaces is a much stronger assumption in many real-world problems than the analyticity assumption of the function $f(\cdot)$, which is needed for use of the degenerate matrix perturbation theory in our analysis.
\end{remark}

We now move on to the lemmas that express $\x_{K} \in \mathcal{B}_{\epsilon}(\x^*)$ for any $K\geq0$ approximately in terms of $\x_{0}$ provided $K$ and $\epsilon$ satisfy certain necessary conditions.
\begin{lemma}\label{lem5}
Given an initialization of the radial vector $\u_{0}$ and $\epsilon < \inf_{{\norm{\u}=1}}\bigg(\limsup_{j\to \infty} \sqrt[j]{\frac{r_j(\u)}{j!}}\bigg)^{-1}$, at any iteration $K$ the radial vector $\u_{K}$ is given by the product
\begin{align}
	\u_{K} = \prod_{k=0}^{K-1} \bigg[\A_k + \epsilon \P_k   \bigg]\u_{0},
\end{align}
where  $\epsilon \P_k = \B_k+\mathcal{O}(\epsilon^2)$, $\B_k = \mathcal{O}(\epsilon)$ for $\x_{k} \in \mathcal{B}_{\epsilon}(\x^*)$ and $\A_k,\B_k$ are given by the following equations:
\begin{align}
	\A_k &= \sum_{i \in \mathcal{N}_{S}}  c^{s}_i(k)\v_i(0)\v_i(0)^{T} + \sum_{j \in \mathcal{N}_{US}}  c^{us}_j(k)\v_j(0)\v_j(0)^{T} \\
		\B_k &= \sum_{i=1}^{n}  \sum_{l \neq i} \bigg(d_{l,i}(k)\v_l(0)\v_i(0)^{T}   +d_{i,l}(k)\v_i(0)\v_l(0)^{T}\bigg).
\end{align}
The coefficient terms $c^{s}_i(k)$, $c^{us}_j(k)$, $d_{i,l}(k)$ and $d_{l,i}(k)$ are as follows:	
\begin{align}
	c^{s}_i(k) &= \bigg(1-\alpha \lambda_{i}(0)-\alpha \frac{\norm{\u_{k}}}{2}\langle \v_i(0), \H(\hat{\u}_k)\v_i(0)\rangle\bigg),   \\
		c^{us}_j(k) &= \bigg(1-\alpha \lambda_{j}(0)-\alpha \frac{\norm{\u_{k}}}{2}\langle \v_j(0), \H(\hat{\u}_k)\v_j(0)\rangle\bigg), \ \text{and}   \\
		d_{i,l}(k) &= d_{l,i}(k) =  \frac{\langle \v_l(0), \H(\hat{\u}_k)\v_i(0)\rangle\lambda_i(0)\alpha \norm{\u_{k}}}{2(\lambda_{l}(0)-\lambda_{i}(0))}.
\end{align}
\revise{Further, suppose $\nu_n \leq \dots  \leq \nu_1$ are the absolute values of the eigenvalues of the matrix $\prod_{k=0}^{K-1} \bigg[\A_k + \epsilon \P_k \bigg]$ and we have that $ \sup_{0\leq k \leq K-1}\norm{\A_k}_2 = \norm{\A}_2$, $ \sup_{0\leq k \leq K-1}\norm{\A_k^{-1}}_2 = \norm{\A^{-1}}_2$ and $\sup_{0\leq k \leq K-1}\norm{\P_k}_2 = \norm{\P}_2$ for some matrices $\A$ and $\P$. Then} for $\epsilon < \min \bigg\{\inf_{{\norm{\u}=1}}\bigg(\limsup_{j\to \infty} \sqrt[j]{\frac{r_j(\u)}{j!}}\bigg)^{-1},\frac{\norm{\A^{-1}}_2^{-1}}{\norm{\P}_2}\bigg\}$ and $K\epsilon \ll 1$, the following condition holds provided $\A_k$ is non-singular for all $k$:
\begin{align}
    \norm{\A^{-1}}_2^{-K} \bigg( 1- K\epsilon \frac{\norm{\P}_2}{\norm{\A^{-1}}^{-1}_2} - \mathcal{O}\bigg((K\epsilon)^2\bigg) \bigg) \leq \nu_n \leq \dots \leq \nu_1  \leq  \norm{\A}_2^K \bigg( 1+ K\epsilon \frac{\norm{\P}_2}{\norm{\A}_2} + \mathcal{O}\bigg((K\epsilon)^2\bigg) \bigg).
\end{align}
\end{lemma}

The proof of this lemma is given in Appendix \ref{Appendix B}. This lemma states that the radial vector $\u_{K}$ evolves linearly at every iteration $K$, where the transition matrix from the initial state $\u_{0}$ to the state $\u_{K}$ is given by $\prod_{k=0}^{K-1} \bigg[\A_k + \epsilon \P_k   \bigg]$. This lemma also states that the absolute value of the eigenvalues of this transition matrix are bounded between terms that are expressed up to $K\epsilon$ precision if $K\epsilon \ll 1$ and $\epsilon$ is upper bounded by the value provided in the lemma. This result is extremely useful in establishing that the matrix product given by $\prod_{k=0}^{K-1} \bigg[\A_k + \epsilon \P_k   \bigg]$ can be computed explicitly up to $K\epsilon$ precision without trading off much on the accuracy of the radial vector $\u_{K}$.

\begin{remark}
Notice that the matrix product $\prod_{k=0}^{K-1} \bigg[\A_k + \epsilon \P_k   \bigg]$ in this lemma is hard to compute where  expansion of this product will generate $K$ terms. The hardness lies in the fact that the higher order terms in $\epsilon$ appearing in the expansion do not simplify due to the fact that matrices $\P_k$ do not commute. Beyond first order the expansion of this matrix product cannot be simplified with ease. Therefore Lemma \ref{lem5} is of utmost importance in the sense that it provides the conditions under which the the tail error generated by the first order approximation $\prod_{k=0}^{K-1} \bigg[\A_k + \epsilon \P_k   \bigg] \approx \prod_{k=0}^{K-1} \A_k + \epsilon \sum\limits_{r=0}^{K}(\Pi_{k=r+1}^{K}\A_r\P_r\Pi_{k=0}^{r-1}\A_r)$ remains bounded.
\end{remark}

\begin{lemma}\label{lem6}
Given an initialization of the radial vector $\u_{0}$, at any iteration $K$ such that $K = \mathcal{O}\bigg(\frac{1}{\epsilon}\bigg)$ and $\epsilon < \min \bigg\{\inf_{{\norm{\u}=1}}\bigg(\limsup_{j\to \infty} \sqrt[j]{\frac{r_j(\u)}{j!}}\bigg)^{-1},\frac{2 \delta(1-\alpha L)}{\alpha M ( 2 L n^2   +\delta)}  + \mathcal{O}(\epsilon^2)\bigg\}$ when $\alpha \in \bigg(0,\frac{1}{L}-\mathcal{O}(\epsilon)\bigg]$ or $ \epsilon < \min \bigg\{\inf_{{\norm{\u}=1}}\bigg(\limsup_{j\to \infty} \sqrt[j]{\frac{r_j(\u)}{j!}}\bigg)^{-1},\frac{2L\delta}{M(2Ln^2 -\delta )} + \mathcal{O}(\epsilon^2)\bigg\} $ when $\alpha \in \bigg(\frac{1}{L}-\mathcal{O}(\epsilon),\frac{1}{L}\bigg]$, the radial vector $\u_{K}$ can be approximately given as
\begin{align}
    \u_{K}\approx \tilde{\u}_K =\epsilon\sum_{i \in \mathcal{N}_{S}}  \bigg(&\prod_{k=0}^{K-1}c^{s}_i(k){\theta}^{s}_i + \sum_{l \in \mathcal{N}_{S}} \sum_{r=0}^{K-1} \prod_{k=0}^{r-1} c^{s}_i(k) d_{i,l}(r)  \prod_{k=r+1}^{K-1} c^{s}_l(k){\theta}^{s}_l \nonumber \\ & + \sum_{l \in \mathcal{N}_{US}} \sum_{r=0}^{K-1} \prod_{k=0}^{r-1} c^{s}_i(k) d_{i,l}(r)  \prod_{k=r+1}^{K-1} c^{us}_l(k){\theta}^{us}_l  \bigg) \v_i(0) + \nonumber \\
    \epsilon \sum_{j \in \mathcal{N}_{US}}  \bigg(&\prod_{k=0}^{K-1}c^{us}_j(k){\theta}^{us}_j + \sum_{l \in \mathcal{N}_{S}} \sum_{r=0}^{K-1} \prod_{k=0}^{r-1} c^{us}_j(k) d_{j,l}(r)  \prod_{k=r+1}^{K-1} c^{s}_l(k){\theta}^{s}_l \nonumber \\ & + \sum_{l \in \mathcal{N}_{US}} \sum_{r=0}^{K-1} \prod_{k=0}^{r-1} c^{us}_j(k) d_{j,l}(r)  \prod_{k=r+1}^{K-1} c^{us}_l(k){\theta}^{us}_l  \bigg)  \v_j(0), \label{explicitfunc1}
\end{align}
where $\epsilon {\theta}^{s}_{i} = \langle \u_{0},\v_i(0)\rangle$, $\epsilon {\theta}^{us}_{j} = \langle \u_{0},\v_j(0)\rangle$ and we have that
\begin{align}
    \u_{0} &= \epsilon\sum_{i \in \mathcal{N}_{S} } {\theta}^{s}_{i} \v_i(0) + \epsilon\sum_{j \in \mathcal{N}_{US} } {\theta}^{us}_{j} \v_j(0)
\end{align}
with ${\theta}^{s}_{i} \geq 0$, ${\theta}^{us}_{j}\geq 0$ for all $i,j$. The coefficient terms $c^{s}_i(k)$, $c^{us}_j(k)$, $d_{i,l}(k)$, $d_{l,i}(k)$ are the same as in Lemma \ref{lem5}.
\end{lemma}

The proof of this lemma is given in Appendix \ref{Appendix B}. The approximation $\tilde{\u}_K$ in this lemma for the radial vector $\u_K$ is generated by explicitly computing the matrix product $\prod_{k=0}^{K-1} \bigg[\A_k + \epsilon \P_k   \bigg]$ from Lemma \ref{lem5} up to first order in $\epsilon$. \revise{Also note that the non-negativity of $\theta_i^s$ and $\theta_j^{us}$ here can be assumed without loss of generality.}
	
\begin{remark}
The conditions $\epsilon < \min \bigg\{\inf_{{\norm{\u}=1}}\bigg(\limsup_{j\to \infty} \sqrt[j]{\frac{r_j(\u)}{j!}}\bigg)^{-1},\frac{2 \delta(1-\alpha L)}{\alpha M ( 2 L n^2   +\delta)}  + \mathcal{O}(\epsilon^2)\bigg\}$ when $\alpha \in \bigg(0,\frac{1}{L}-\mathcal{O}(\epsilon)\bigg]$ or $ \epsilon < \min \bigg\{\inf_{{\norm{\u}=1}}\bigg(\limsup_{j\to \infty} \sqrt[j]{\frac{r_j(\u)}{j!}}\bigg)^{-1},\frac{2L\delta}{M(2Ln^2 -\delta )}+ \mathcal{O}(\epsilon^2) \bigg\} $ when we have $\alpha \in \bigg(\frac{1}{L}-\mathcal{O}(\epsilon),\frac{1}{L}\bigg]$ are necessary but may not be sufficient due to unavailability of the radius of convergence from the Rayleigh--Schr\"{o}dinger perturbation analysis. Also note that here $r_j(\u)$ has the same definition as in Lemma \ref{lem4}.
\end{remark}
		
In words, this lemma states that the radial vector $\u_{K}$ can be expressed by explicitly computing the matrix product $\prod_{k=0}^{K-1} \bigg[\A_k + \epsilon \P_k   \bigg]$ from Lemma \ref{lem5} to $K\epsilon$ precision provided $K\epsilon \ll 1$ and $\epsilon$ is bounded above. This approximate solution represented by $ \tilde{\u}_K$ generates the trajectory $ \{\tilde{\u}_K\}_{K=1}^{K_{exit}}$, which we refer to as the $\epsilon$-precision trajectory.
\begin{remark}
Notice that from \eqref{explicitfunc1} we obtain a closed form expression for the $\epsilon$ precision trajectory inside $\mathcal{B}_{\epsilon}(\x^*)$ for some initialization $\u_0$. However the solution is not unique due to the fact that the coefficients $c^{s}_i(k), c^{us}_j(k), d_{l,i}(k) $ from Lemma \ref{lem5} are known only up to an interval. This is due to the fact that the eigenvalues $\lambda_i(0), \lambda_j(0)$ are known up to an interval. Hence we will obtain a \textbf{family of $\epsilon$ precision trajectories} from the expression of $\tilde{\u}_K$. The next lemma provides a handle on the exit times for such a family of approximate trajectories.
\end{remark}

\begin{figure}
    \centering
    \includegraphics[width=0.85\textwidth]{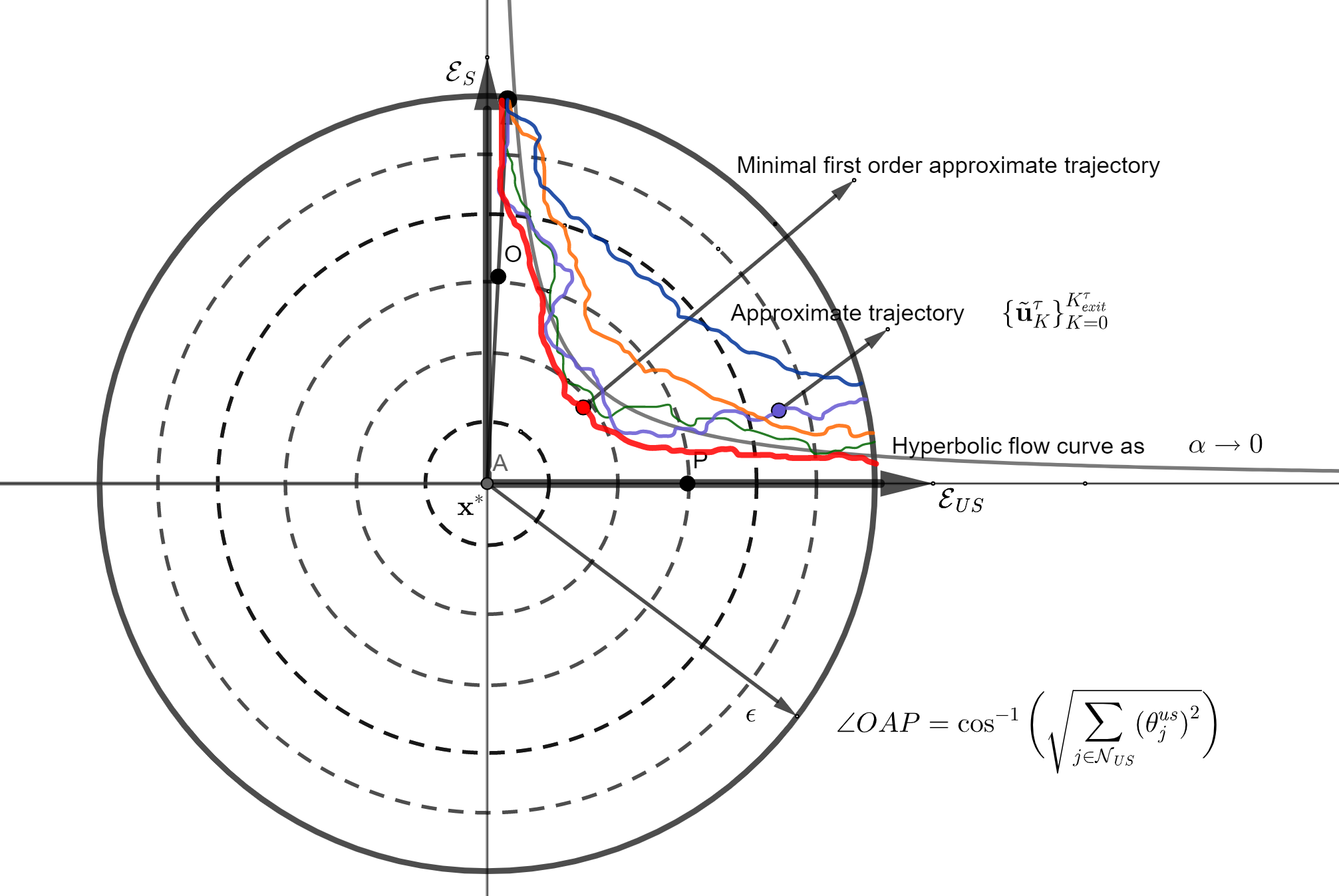}
    \caption{A 2-D representation of the approximate trajectories, where every approximate trajectory has its own exit time and the \textit{minimal} approximate trajectory is the one which has the largest exit time. The initial radial vector subtends a very large angle of $\angle OAP $ ($0 \ll \angle OAP < \frac{\pi}{2}$) from the unstable subspace, where the initial unstable projection is given by $ {\sum_{j \in \mathcal{N}_{US}} ({\theta}^{us}_j)^2}$.}
    \label{familyofcurves}
\end{figure}

\begin{lemma}\label{lem7}
Let $S_{\epsilon}=\bigg\{  \{\tilde{\u}_K^{\tau}\}_{K=1}^{K_{exit}^{\tau}}  \bigg| \u_{0} \bigg\}$ be the set of all possible $\tau$-parameterized $\epsilon$-precision trajectories generated by the approximate equation \eqref{explicitfunc1} in Lemma \ref{lem6}, where the parameter $\tau \in \mathbb{R}$ varies with variations in the sequence $\bigg \{\{c^{s}_i(k), c^{us}_j(k), d_{l,i}(k)\}_{k=0}^{K-1} \bigg \}_{K=1}^{K_{exit}}$. Let $K_{exit}^{\tau}$ be the exit time of the $\tau$-parameterized trajectory $\{\tilde{\u}_K^{\tau}\}_{K=1}^{K_{exit}^{\tau}}$ from the ball $\mathcal{B}_{\epsilon}(\x^*)$ where we have that
\begin{align}
	 K_{exit}^{\tau} & = \inf_{K\geq 1} \bigg\{K \hspace{0.2cm}\bigg| \hspace{0.2cm} \norm{\tilde{\u}_K^{\tau}}^2 > \epsilon^2 \bigg\}.
\end{align}
Formally, $\tilde{\u}_K^{\tau}$ is a possible solution to the equation \eqref{explicitfunc1} in $\tilde{\u}_K$, where $1\leq K\leq K_{exit}^{\tau}$ and $\tilde{\u}_K$ varies with variations in the sequence $\{c^{s}_i(k), c^{us}_j(k), d_{l,i}(k)\}_{k=0}^{K-1} $.
	
Let $K^{\iota}$ be the exit time of the infimum over all possible $\tau$-parameterized trajectories, where infimum is taken with respect to the squared radial distance $\norm{\tilde{\u}_K^{\tau}}^2$. This $K^{\iota}$ can be defined as
\begin{align}
	K^{\iota} & = \inf_{K\geq 1} \bigg\{K \hspace{0.2cm}\bigg| \hspace{0.2cm} \inf_{\tau}\bigg\{\norm{\tilde{\u}_K^{\tau}}^2 \bigg \}> \epsilon^2 \bigg\}.
\end{align}
Then we have the following condition:
\begin{align}
	K^{\iota} &\geq \sup_{\tau} \bigg\{K_{exit}^{\tau} \bigg\}  = \sup_{\tau} \inf_{K\geq 1} \bigg\{K \hspace{0.2cm}\bigg| \hspace{0.2cm} \norm{\tilde{\u}_K^{\tau}}^2 > \epsilon^2 \bigg\}.
\end{align}
\end{lemma}

The proof of this lemma is given in Appendix \ref{Appendix C}. This particular lemma states an important result about the exit time $K^{\iota}$ of the  trajectory generated by selecting that approximate vector $\tilde{\u}^{\tau}_K$ from all possible $\tau$ that has the minimum radial distance from $\x^*$ at each $K$. It claims that this minimal trajectory has the maximum exit time from $\mathcal{B}_{\epsilon}(\x^*)$. Though seemingly trivial, this result is extremely important in proving the worst-case exit time for trajectories with linear escape rates. A representation of the family of approximate trajectories along with the constructed minimal approximate trajectory is provided in Figure \ref{familyofcurves}.

\begin{theorem}\label{theorem1}
For every value of the parameter $\tau$, there exists a lower bound on the squared radial distance $\norm{\tilde{\u}_K^{\tau}}^2$ for all $K$ in the range $1 \leq K \leq \sup_{\tau} \bigg\{K_{exit}^{\tau} \bigg\}  $ provided $K\epsilon \ll 1$. Moreover, this lower bound can be expressed using a function of $K$ called the trajectory function $\Psi(K)$. Formally, for $1 \leq K < \sup_{\tau} \bigg\{K_{exit}^{\tau} \bigg\}  $ we have that
\begin{align}
	   \epsilon^2 \geq\inf_{\tau}\norm{\tilde{\u}_K^{\tau}}^2 > &  \epsilon^2 \Psi(K) ,
\end{align}
where the trajectory function $\Psi(K)$ is defined as follows:
\begin{align}
    \hspace{-0.2cm}\Psi(K) =& \bigg(c_1^{2K} -2Kc_2^{2K-1}b_1 - b_2c_3^Kc_2^K - b_2c_3^{2K}\bigg)\sum_{i \in \mathcal{N}_{S}}({\theta}^{s}_i)^2 + \bigg( c_4^{2K} - 2Kc_3^{2K-1}b_1- b_2c_3^Kc_2^K -b_2c_3^{2K}\bigg)\sum_{j \in \mathcal{N}_{US}}({\theta}^{us}_j)^2
\end{align}
with $c_1 = \bigg(1 - \alpha L - \frac{\alpha \epsilon M}{2}- \mathcal{O}(\epsilon^2) \bigg)$, $c_2 =\bigg(1 - \alpha \beta + \frac{\alpha \epsilon M}{2}+\mathcal{O}(\epsilon^2)\bigg) $, $c_3 = \bigg(1 + \alpha L + \frac{\alpha \epsilon M}{2}+\mathcal{O}(\epsilon^2)\bigg)$, $c_4 =\bigg(1 + \alpha \beta - \frac{\alpha \epsilon M}{2}-\mathcal{O}(\epsilon^2)\bigg)$, $b_1 = \bigg(\frac{\alpha\epsilon M L n}{2 \delta}+\mathcal{O}(\epsilon^2)\bigg)$, $b_2 = \frac{\bigg(\frac{\alpha\epsilon M L n}{2 \delta}+\mathcal{O}(\epsilon^2)\bigg)\bigg(1+\mathcal{O}(K\epsilon)\bigg)}{ \bigg( \alpha L + \alpha\beta + \mathcal{O}(\epsilon^2)\bigg)}$ \revise{and $\delta$ is defined in Proposition \ref{eigenprop}}.

We also require that $\epsilon < \min \bigg\{\inf_{{\norm{\u}=1}}\bigg(\limsup_{j\to \infty} \sqrt[j]{\frac{r_j(\u)}{j!}}\bigg)^{-1},\frac{2 \delta(1-\alpha L)}{\alpha M ( 2 L n^2   +\delta)}  + \mathcal{O}(\epsilon^2)\bigg\}$ when $\alpha \in \bigg(0,\frac{1}{L}-\mathcal{O}(\epsilon)\bigg]$, while $ \epsilon < \min \bigg\{\inf_{{\norm{\u}=1}}\bigg(\limsup_{j\to \infty} \sqrt[j]{\frac{r_j(\u)}{j!}}\bigg)^{-1},\frac{2L\delta}{M(2Ln^2 -\delta )}+ \mathcal{O}(\epsilon^2) \bigg\} $ when we have $\alpha \in \bigg(\frac{1}{L}-\mathcal{O}(\epsilon),\frac{1}{L}\bigg]$.
\end{theorem}

The proof of this theorem is given in Appendix \ref{Appendix C}. Theorem \ref{theorem1} states that for a given initialization $\u_{0}$, all the possible $\epsilon$-precision trajectories generated have their radial distance from $\x^*$ lower bounded using some function $\Psi(K)$. Now this $\Psi(K)$ can be used to determine $K^{\iota}$ and hence $K_{exit}$ for any choice of the step size $\alpha$.

\begin{remark}\label{remark-trajectory-func}
Notice that the trajectory function $ \Psi(K)$ corresponds to the minimal approximate trajectory $\inf_{\tau} \norm{\tilde{\u}_K^{\tau}}$. Now $K^{\iota}$ can be obtained by solving the condition $K^{\iota} = \inf_{\substack{K\geq 1\\ {K = \mathcal{O}(\log(\epsilon^{-1}))}}}\{K \hspace{0.1cm} \vert \hspace{0.1cm}\Psi(K) > 1\}$. The condition $ K = \mathcal{O}(\log(\epsilon^{-1}))$ ensures linear time solutions, which are the only solutions of interest to the problem. Then from Lemma \ref{lem7} we will have $K_{exit} < K^{\iota} = \mathcal{O}(\log(\epsilon^{-1}))$. \revise{It is worth reminding the reader here that the $\mathcal{O}(\epsilon^2)$ terms in the theorem statement account for the approximation in analysis and things are proved about
trajectories and perturbations up to first order in $\epsilon$.}
\end{remark}

\revise{Observe that in the expression for the trajectory function $\Psi(K)$, the term accompanying $\sum_{i \in \mathcal{N}_{S}}({\theta}^{s}_i)^2 $ is $\bigg(c_1^{2K} -2Kc_2^{2K-1}b_1 - b_2c_3^Kc_2^K - b_2c_3^{2K}\bigg)$, which is a decreasing function of $K$ since $c_1<1$. Moreover, the rate of decrease of the term $\bigg(c_1^{2K} -2Kc_2^{2K-1}b_1 - b_2c_3^Kc_2^K - b_2c_3^{2K}\bigg)$ for small values of $K$ is governed by $ c_1$ and not by $c_3$, where $c_3>1$ due to the fact that $b_1, b_2 $ are of order $\mathcal{O}(\epsilon)$ and so $-2Kc_2^{2K-1}b_1$, $-b_2c_3^{2K} $ will not decrease as fast as $ c_1^{2K}$ for small enough $K$ since we assumed $K\epsilon \ll 1$. Next, by a similar argument the term accompanying $\sum_{j \in \mathcal{N}_{US}}({\theta}^{us}_j)^2 $ given by $\bigg( c_4^{2K} - 2Kc_3^{2K-1}b_1- b_2c_3^Kc_2^K -b_2c_3^{2K}\bigg) $ is an increasing function of $K$ for $K\epsilon \ll 1$ since $c_4>1$ and so $c_4^{2K} $ dominates the term $ 2Kc_3^{2K-1}b_1+ b_2c_3^Kc_2^K +b_2c_3^{2K} $. Also notice that $\Psi(K)<1$ at $K=0$ since $\sum_{i \in \mathcal{N}_{S}}({\theta}^{s}_i)^2 + \sum_{j \in \mathcal{N}_{US}}({\theta}^{us}_j)^2 = 1 $. Therefore, provided the initial unstable subspace projection $\sum_{j \in \mathcal{N}_{US}}({\theta}^{us}_j)^2 $ is not too small, the trajectory function $\Psi(K)$ first increases for small $K$, where $K\epsilon \ll 1$, and then decreases to $-\infty$. Then for some small $K$ if $\Psi(K)>1$, we are guaranteed that the minimal approximate trajectory $\inf_{\tau} \norm{\tilde{\u}_K^{\tau}}$ escapes $\mathcal{B}_{\epsilon}(\x^*)$. Section~\ref{psik_simul} simulates the evolution of the trajectory function $\Psi(K)$ on the phase retrieval problem, which corroborates this theoretical understanding.}

Before moving on to the next theorem, we introduce the notion of conditioning of a function. The condition number  at the stationary point of a non-convex function is given by the ratio of the largest absolute eigenvalue to the smallest absolute eigenvalue of the Hessian of the function at that point. Also, a function is called perfectly conditioned if the condition number is equal to $1$. In the current problem setting, the condition number of the function $f(\cdot)$ at the saddle point $\x^*$ is given by $\frac{L}{\beta}$. Now, the function $f(\cdot)$ is well-conditioned if the condition number $\frac{L}{\beta}$ is not arbitrarily large or equivalently $\frac{\beta}{L}$ is bounded away from $0$.

\begin{theorem}\label{theorem2}
For the gradient update equation with the step size $\alpha = \frac{1}{L}$, \revise{there exists a minimum projection $\Delta$} of the radial vector initialization $\u_{0}$ on the unstable subspace $\mathcal{E}_{US}$ such that \revise{whenever $ \sum_{j \in \mathcal{N}_{US}} ({\theta}^{us}_j)^2 \geq  \Delta$, where} $\u_{0}+\x^* \in \mathcal{\bar{B}}_{\epsilon}(\x^*) \backslash \mathcal{B}_{\epsilon}(\x^*)$, $\u_{0} = \epsilon\sum_{i \in \mathcal{N}_{S} } {\theta}^{s}_{i} \v_i(0) + \epsilon\sum_{j \in \mathcal{N}_{US} } {\theta}^{us}_{j} \v_j(0)$, \revise{the} $\epsilon$-precision trajectories $\{\tilde{\u}_{K}\}_{K=1}^{K_{exit}}$ \revise{can exit $\mathcal{B}_{\epsilon}(\x^*)$ in  linear time.} Moreover their exit time $K_{exit}$ from the ball $\mathcal{B}_{\epsilon}(\x^*)$ is approximately upper bounded as follows:
\begin{align}
K_{exit} < K^{\iota} & \lessapprox \frac{\log \bigg(\bigg(2 + \frac{ \epsilon M}{2L}\bigg) \log\bigg(\frac{2 + \frac{ \epsilon M}{2L}}{1 + \frac{\beta}{L} - \frac{\epsilon M}{2L}}\bigg)\frac{2 \delta}{\epsilon M n}\bigg)}{2\log\bigg(\frac{2 + \frac{ \epsilon M}{2L}}{1 + \frac{\beta}{L} - \frac{\epsilon M}{2L}}\bigg)},
\end{align}
where $ \epsilon < \min \bigg\{\inf_{{\norm{\u}=1}}\bigg(\limsup_{j\to \infty} \sqrt[j]{\frac{r_j(\u)}{j!}}\bigg)^{-1},\frac{2L\delta}{M(2Ln^2 -\delta )} + \mathcal{O}(\epsilon^2), \revise{\frac{2 \beta}{M}}\bigg\} $ \revise{and we must necessarily have that} $\Delta >\epsilon \frac{ M Ln}{ \delta (L+ \beta )}$ \revise{with $\delta$ defined in Proposition \ref{eigenprop}}.
\end{theorem}

The proof of this theorem is given in Appendix \ref{Appendix D}. \revise{In terms of the order notation, we have $K_{exit} \lessapprox \mathcal{O}\bigg(\log\bigg( \frac{\delta}{\epsilon n} \bigg)\bigg)$ and the initial unstable subspace projection satisfies $\sum_{j \in \mathcal{N}_{US}} ({\theta}^{us}_j)^2 \geq\Delta >\Omega\bigg(\frac{\epsilon n}{\delta}\bigg)$.}
\begin{remark}
This theorem guarantees the existence of $\epsilon$-precision trajectories with linear exit time and gives an upper bound on their exit time $K_{exit}$ from the ball $\mathcal{B}_{\epsilon}(\x^*)$. However, the sufficient conditions that guarantee the existence of this exit time \revise{$K_{exit} \lessapprox \mathcal{O}\bigg(\log\bigg(\frac{1}{\epsilon}\bigg)\bigg)$} depend on the quantity $\Delta$. Note that the condition $  \Delta >\epsilon \frac{ M Ln}{ \delta (L+ \beta )}$ is necessary for the existence of order $\mathcal{O}\bigg(\log\bigg(\frac{1}{\epsilon}\bigg)\bigg)$ solution of $K^{\iota}$ but not sufficient. Since this work only deals with the existence of linear exit time solutions, we refrain from developing tighter lower bounds on $\Delta$. Obtaining such sufficient conditions requires a more rigorous analysis of the trajectory function $\Psi(K)$, which is beyond the scope of current work. \revise{In particular, our followup work \cite{dixit2022boundary} derives one such sufficient condition.}
\end{remark}

\begin{remark}
\revise{Observe that the bound on the exit time from Theorem \ref{theorem2} depends on quantities like the Lipschitz parameters, condition number, problem dimension and the eigen gap. However for structured problems such as those in \cite{chen2019gradient}, one can leverage the specialized function geometry and obtain rates of convergence independent of these parameters. But in the absence of any other assumption on the function, and since we are dealing with a much larger function class, i.e., the class of Morse functions, these parameters become necessary to evaluate the escape rates. In order to better understand the utility of these local Lipschitz parameters in the derivation of our results for the general (as opposed to the specialized) non-convex functions, observe that the local Hessian Lipschitz parameter $M$ is required to bound $ \norm{ \textbf{H}(\hat{\u}_k) }_2$, where $\textbf{H}(\hat{\u}_k) $ is used to determine the Hessian at any point $\x_k \in \mathcal{B}_{\epsilon}(\x^*)$ from Lemma \ref{lem4}. Next, the local gradient Lipschitz parameter $L$ controls the coefficient terms $c^{s}_i(k)$, $c^{us}_j(k)$, $d_{i,l}(k)$ from Lemmas \ref{lem5} and \ref{lem6}, where these terms depend on the eigenvalues of $\nabla^2 f(\x^*)$, the difference between these eigenvalues, and the matrix $\textbf{H}(\hat{\u}_k)$, which comes from Lemma \ref{lem4}. Since these coefficient terms determine the expression for the approximate gradient trajectory in Lemma \ref{lem6}, one cannot generate a closed-form expression of the approximate gradient trajectory in the absence of the gradient Lipschitz parameter. Finally, the minimal approximate trajectory function from Theorem \ref{theorem1} relies on the precise bounds for these coefficients. Without the gradient Lipschitz parameter, the eigenvalues of $\nabla^2 f(\x^*)$ cannot be bounded and similarly without the Hessian Lipschitz parameter one cannot obtain an upper bound on $ \norm{ \textbf{H}(\hat{\u}_k) }_2$.}
\end{remark}

Theorem \ref{theorem2} guarantees a linear exit time bound from the ball $\mathcal{B}_{\epsilon}(\x^*)$ for $\epsilon$-precision trajectories under some necessary initial conditions on $\x_{0}$. The necessary condition of $ \Delta >\epsilon \frac{ M Ln}{ \delta (L+ \beta )}$ requires that the initial radial vector is not aligned too much with the stable subspace of the Hessian $\nabla^2 f(\x^*)$ and has some \revise{order $\Omega(\epsilon)$} alignment with the unstable subspace so as to facilitate the linear time escape. It should be noted that this necessary condition of $\Delta >\epsilon \frac{ M Ln}{ \delta (L+ \beta )}$ is sufficient to claim that these gradient descent trajectories \revise{for $\alpha < \frac{1}{L}$ will almost surely not} terminate into the strict saddle point $\x^*$ from the following {Lemma \ref{lem3}}.

\begin{lemma}\label{lem3}
The discrete gradient trajectories \revise{for $\alpha < \frac{1}{L}$} ending into the first-order strict saddle point $\x^*$ have zero Lebesgue measure with respect to the space $\mathcal{E}$ and are referred to as trivial trajectories. This result can be established using the stable center manifold theorem from \cite{shub2013global}.
\end{lemma}

We refer the reader to \cite{lee2016gradient} for a proof of this lemma. \revise{Note that the assumption on the step size $\alpha < \frac{1}{L}$ in Lemma \ref{lem3} is necessary since the zero measure result can only be developed when the map $G : \x_k \mapsto \x_{k+1}$, where $\x_{k+1} = \x_k - \alpha \nabla f(\x_k) =: G(\x_k)$, is a diffeomorphism (or is at least locally bi-Lipschitz). A crucial step in \cite{lee2016gradient} where this diffeomorphism property is utilized involves pulling back measure zero sets under the diffeomorphism $G$ to again get measure zero sets. However for the case of $\alpha = \frac{1}{L}$, the map $G : \x_k \mapsto \x_{k+1}$ fails to be a diffeomorphism (or even locally bi-Lipschitz); see details in \cite{lee2016gradient}.}\looseness=-1

\revise{We also note that the condition of minimal non-zero projection of the initial point on the unstable subspace of $\nabla^2 f(\x^*)$ from Theorem \ref{theorem2}, given by the bound $\Delta >\epsilon \frac{ M Ln}{ \delta (L+ \beta )}$, is tight. Moreover, this necessary condition does not contradict any existing results regarding the almost sure non-convergence of randomly initialized gradient descent to strict saddle points. Further, recall that the gradient descent method may get stuck at the saddle point for a particular set of initializations. In Theorem~\ref{theorem2}, however, we provide a condition on the initialization that ensures its exclusion from such a set. This condition, which is one of the major contributions of this work, requires the projection of the initial point on the unstable subspace of the Hessian $\nabla^2 f(\x^*)$ at the saddle point $\x^*$ to be at least on the order of $\Omega(\epsilon)$. Take, for instance, a specific example of the strict saddle Morse function $f(x,y) = x^2 - y^2$ with the initialization scheme of $(x_0, 0)$ for any $x_0 \in \mathbb{R}$. Under this given initialization scheme, the gradient descent method will eventually get stuck at the origin, which is a strict saddle point. However, since the initialization point completely lies in the stable subspace of $\nabla^2 f(0,0)$, which is $\textsf{span}\{(1,0)\}$, it has a null projection on the unstable subspace of $\nabla^2 f(0,0)$, which is $\textsf{span}\{(0,1)\}$. Therefore, this example violates the minimal projection condition of Theorem~\ref{theorem2} and does not affect the validity of our claims.}

\subsection{Comparison with the exit time bound from Section \ref{crudeanalysissection} }\label{ssec:comparison}

\begin{table}[h!]
\centering
\resizebox{1.1\columnwidth}{!}{
\renewcommand{\arraystretch}{2}
\hspace{-1cm}
\begin{tabular}{||c c c ||}
 \hline
 \textbf{Assumptions / Techniques / Metrics} & \textbf{Exit Time Analysis from Section \ref{crudeanalysissection}} & \textbf{Exit Time Analysis from Section \ref{exactexittimeanalysissection}}  \\ [0.5ex]
 \hline\hline
 Function class & $\mathcal{C}^2$ Morse functions & locally $\mathcal{C}^{\omega}$ Morse functions  \\
 \hline
 Proof techniques & Sequential monotonicity of  & Matrix perturbation theory and  \\
 & the unstable subspace projection & approximation theory\\
 \hline
 Key metrics & Saddle neighborhood's radius $\epsilon$,  & Saddle neighborhood's radius $\epsilon$, \\
  &unknown factors $\gamma_{K_{exit}}, \rho$ & dimension $n$ and eigenvalue gap $\delta$ \\
    \hline
    Closed-form expression for the trajectory /  & \xmark  & \cmark \\
    approximate trajectory inside $\mathcal{B}_{\epsilon}(\x^*)$ & & \\
    \hline
    Constraints on the set of trajectories /  & Gradient trajectories for which $\langle \v_n, \x_k - \x^* \rangle$ & No constraints\\
    approximate trajectories analyzed &  increases monotonically with $k$ & \\
    \hline
Linear exit time bound & $ \mathcal{O}\bigg( \log\bigg(\frac{\gamma_{K_{exit}}(1-\rho)}{ \epsilon}\bigg)\bigg) $ & $\mathcal{O}\bigg(\log\bigg( \frac{\delta}{\epsilon n} \bigg)\bigg)$  \\
\hline
 Nature of the exit time bound & Exact & Approximate  \\
 \hline
 Initial boundary conditions & $\langle\v_0, \x_{0} - \x^* \rangle \geq \Omega\bigg(\frac{ \epsilon^2}{1-\rho}\bigg)$  & $\sum_{j \in \mathcal{N}_{US}} ({\theta}^{us}_j)^2 \geq \Delta >\Omega\bigg( \frac{\epsilon n}{ \delta }\bigg)$  \\
 \hline
  Bounds on $\epsilon$ & \xmark &  \cmark\\
 [1ex]
 \hline
\end{tabular}
\renewcommand{\arraystretch}{1}}
\caption{\revise{Comparison of the exit time analyses that follow from existing analytical techniques (Section \ref{crudeanalysissection}) and the novel matrix perturbation-based analytical approach of Section \ref{exactexittimeanalysissection}.}}
\label{table:1}
\end{table}

It can be seen from Theorem \ref{theorem2} that the exit time bound for the approximate trajectory and the necessary initial condition using the matrix perturbation-based analysis depend on quantities like the inverse of the condition number $ \frac{\beta}{L}$, minimum eigenvalue gap $\delta$, function's dimension $n$ and the size of the saddle neighborhood $\epsilon$. In contrast, the rudimentary analysis in Section \ref{crudeanalysissection} does not bring out the dependence of the exit time bound and the initial boundary condition on these key problem parameters.  
Moreover, the analysis developed in Section \ref{crudeanalysissection} leaves more open questions by introducing unknown parameters like $\rho$ and $\gamma_{K_{exit}}$, where $\gamma_{K_{exit}}$ could be arbitrarily small and the presence of $\rho$ in the boundary condition makes the analysis from Section \ref{crudeanalysissection} more restrictive than the analysis presented in Section \ref{exactexittimeanalysissection} where matrix perturbation theory is used. The main reason for this difference between the results of Section \ref{crudeanalysissection} and those of Theorem \ref{theorem2} is that, by restricting the class of functions from $\mathcal{C}^2$ to real analytic, we are able to develop tight approximations to discrete trajectories using the matrix perturbation theory that lead to precise expressions for the exit time bound and the initial boundary condition that depend on the key problem parameters. These differences between the rudimentary analytical approach of Section \ref{crudeanalysissection} and the matrix perturbation-based approach of Section \ref{exactexittimeanalysissection} are also summarized in Table \ref{table:1}. 
Notice that there is a cross (\xmark) marked against the `Closed form expression for the trajectory' in Table \ref{table:1} in the column corresponding to the analysis of Section~\ref{crudeanalysissection}. This is because although \eqref{crudeanalysis2} provides an expression for the trajectory inside the ball $\mathcal{B}_{\epsilon}(\x^*)$, its exact closed form cannot be determined due to the fact that we only have information on $\norm{r(\x_k)}$ in Section \ref{crudeanalysissection}. In contrast, the same $r(\x_k)$ is known up to first-order precision in Section \ref{exactexittimeanalysissection} and therefore a closed-form expression for the $\epsilon$-precision trajectory is available from Lemma \ref{lem6}.

\section{Numerical results}
To support the theoretical framework developed in this work and showcase the effectiveness of gradient trajectories with large initial unstable projections in escaping from strict saddle neighborhoods, we evaluate the performance of the gradient descent method on the phase retrieval problem \cite{candes2015phase}. Briefly, the phase retrieval problem formulation is given by
\begin{align}
    \min_{\x \in \mathbb{R}^n} f(\x) = \frac{1}{4m} \sum\limits_{j=1}^{m}\bigg[\langle\a_j,\x \rangle^2 - y_j \bigg]^2, \label{phaseretrieval}
\end{align}
where the $y_j$'s are known observations and the $\a_j$'s are independent and identically distributed (i.i.d.) random vectors whose entries are generated from a normal distribution. \revise{Note that the variable '$m$' here in \eqref{phaseretrieval} should not be confused with the number of eigenvalue groups '$m$' defined in proposition \ref{eigenprop}.} The formulation in \eqref{phaseretrieval} is the least-squares problem reformulation for the Short-Time Fourier Transform (STFT) of the actual phase retrieval problem (see \cite{jaganathan2016stft}). Moreover, the above least-squares reformulation of the original phase retrieval problem can also be found in recent works like \cite{ma2020implicit}, which highlight the efficacy of simple gradient descent method on structured non-convex functions. Clearly, the function in \eqref{phaseretrieval} satisfies Assumption \revise{\textbf{A1} and also Assumptions \textbf{A2} and \textbf{A3} locally in every compact set.}\looseness=-1

\begin{figure}[h]
\centering
\begin{tabular}{cc}
   \hspace{-0.9cm} \includegraphics[width=3in]{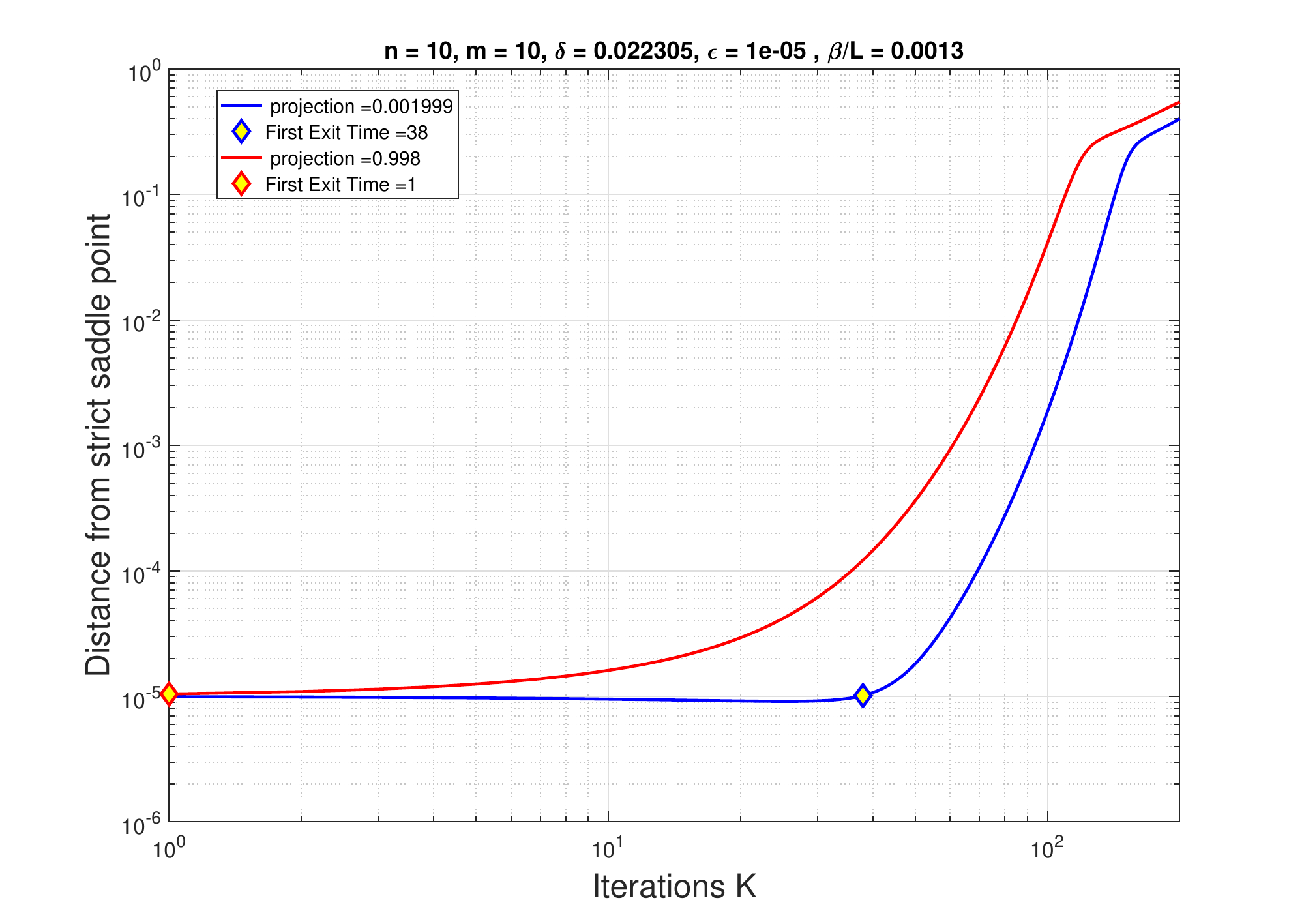} & \includegraphics[width=3in]{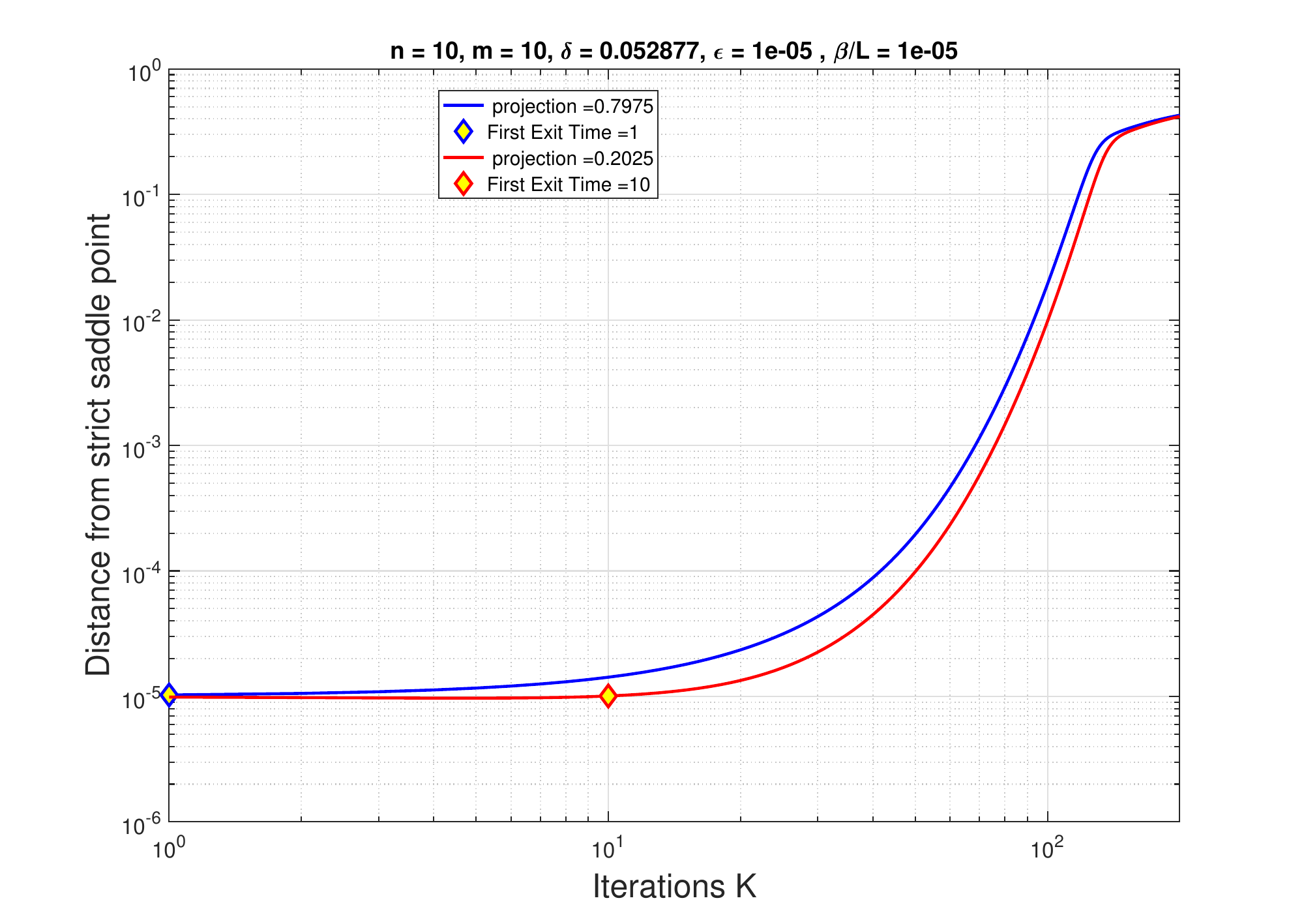} \\
\end{tabular}
\begin{tabular}{cc}
 (a) & (b) \\
   \hspace{-0.9cm} \includegraphics[width=3in]{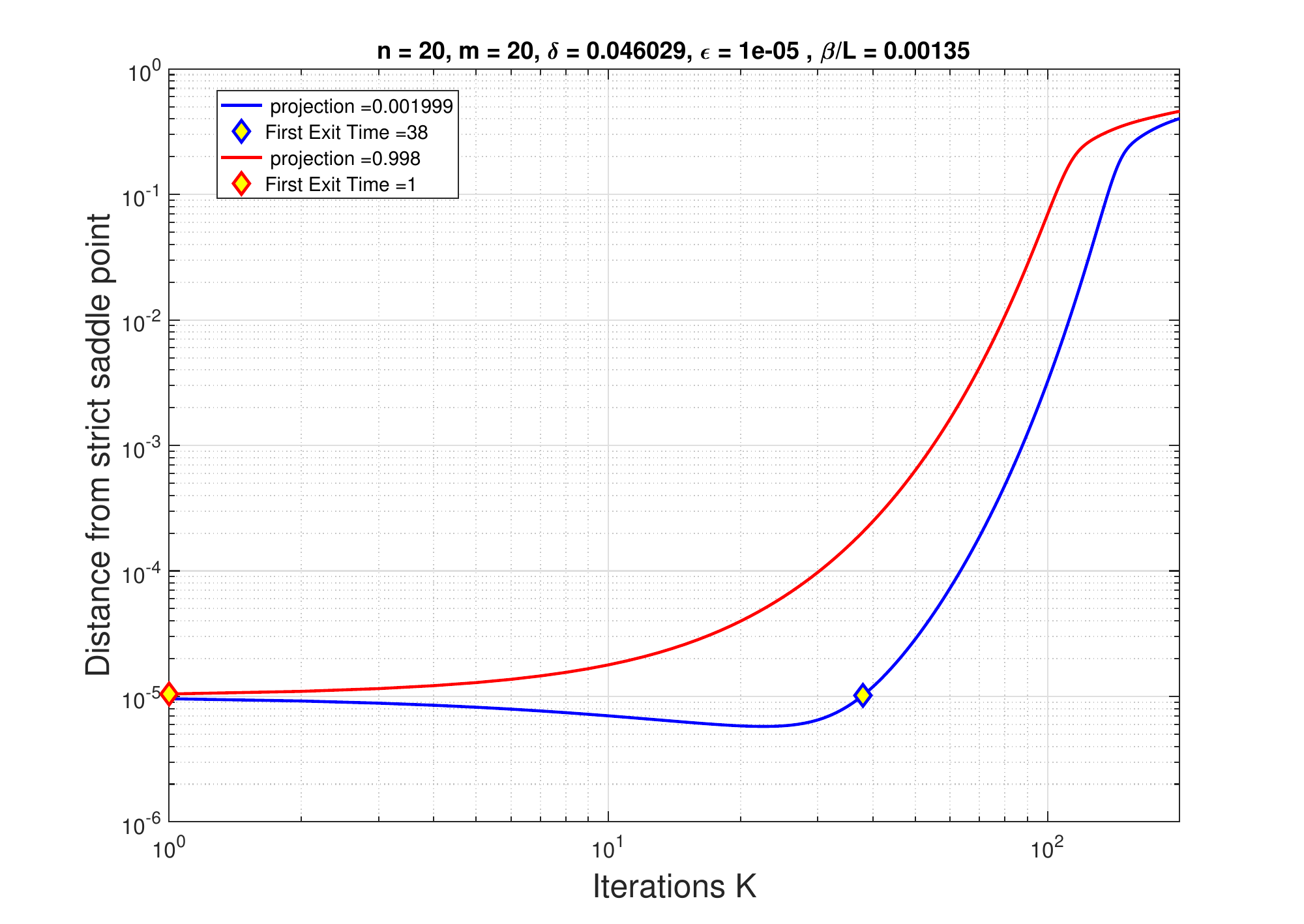} & \includegraphics[width=3in]{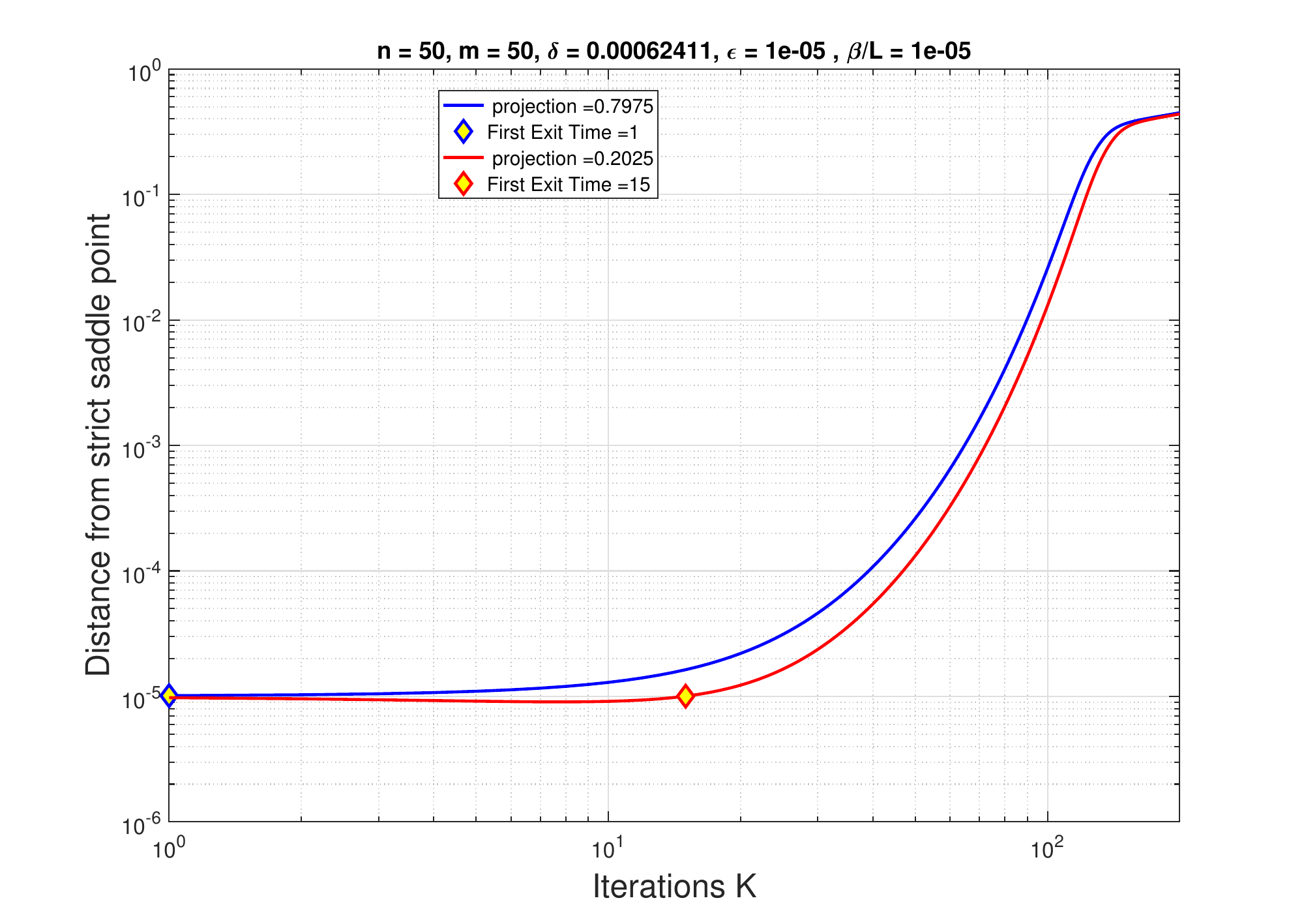} \\
    (c) & (d)
\end{tabular}
\caption{\revise{Simulating gradient trajectories on the phase retrieval problem with $\alpha = 0.1/L$ under certain initial unstable projections for various values of $m$, $n$ and $\epsilon$.}}
\label{fig11}
\end{figure}

\begin{figure}[h]
\centering
\begin{tabular}{cc}
   \hspace{-0.9cm}  \includegraphics[width=3in]{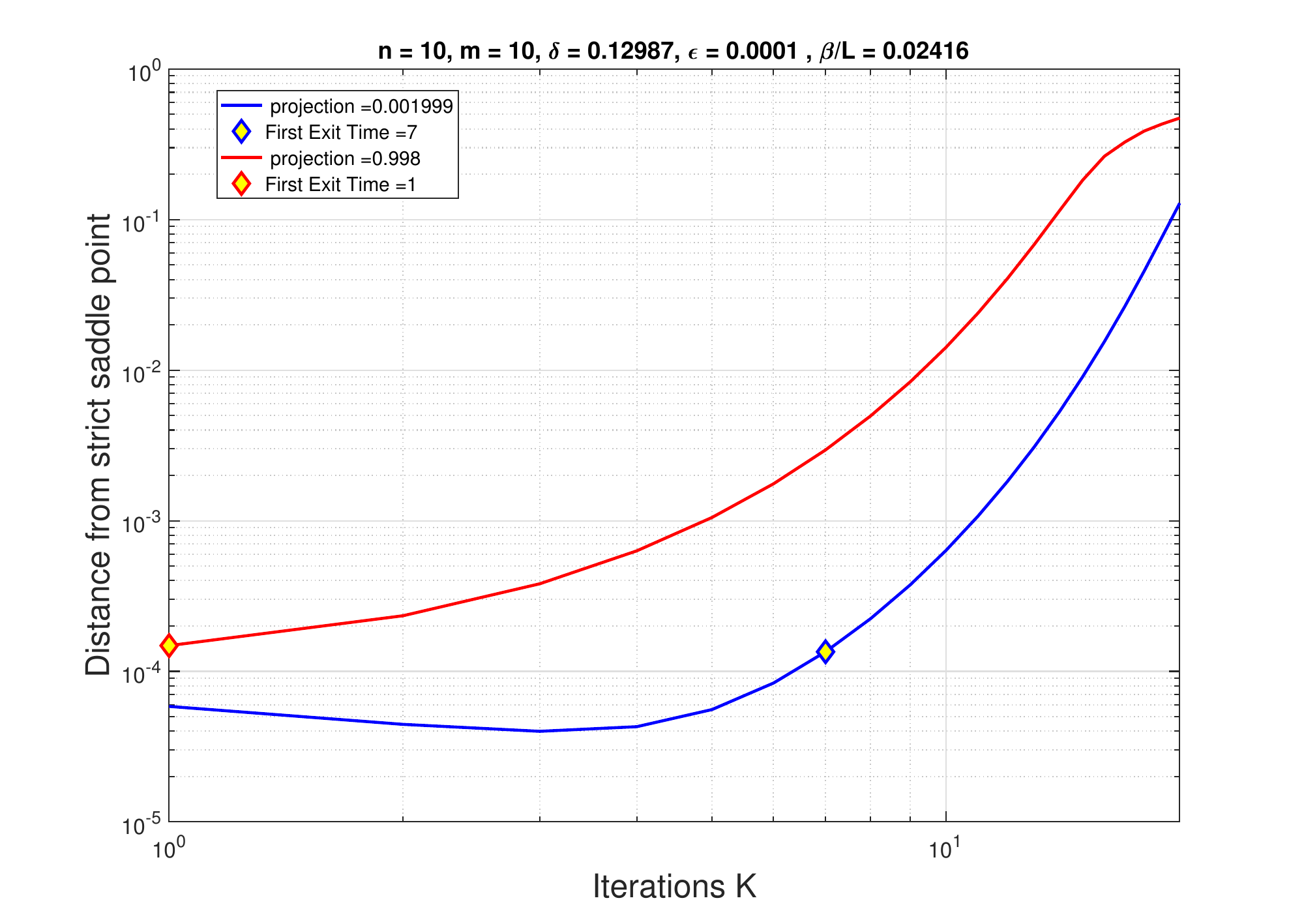} & \includegraphics[width=3in]{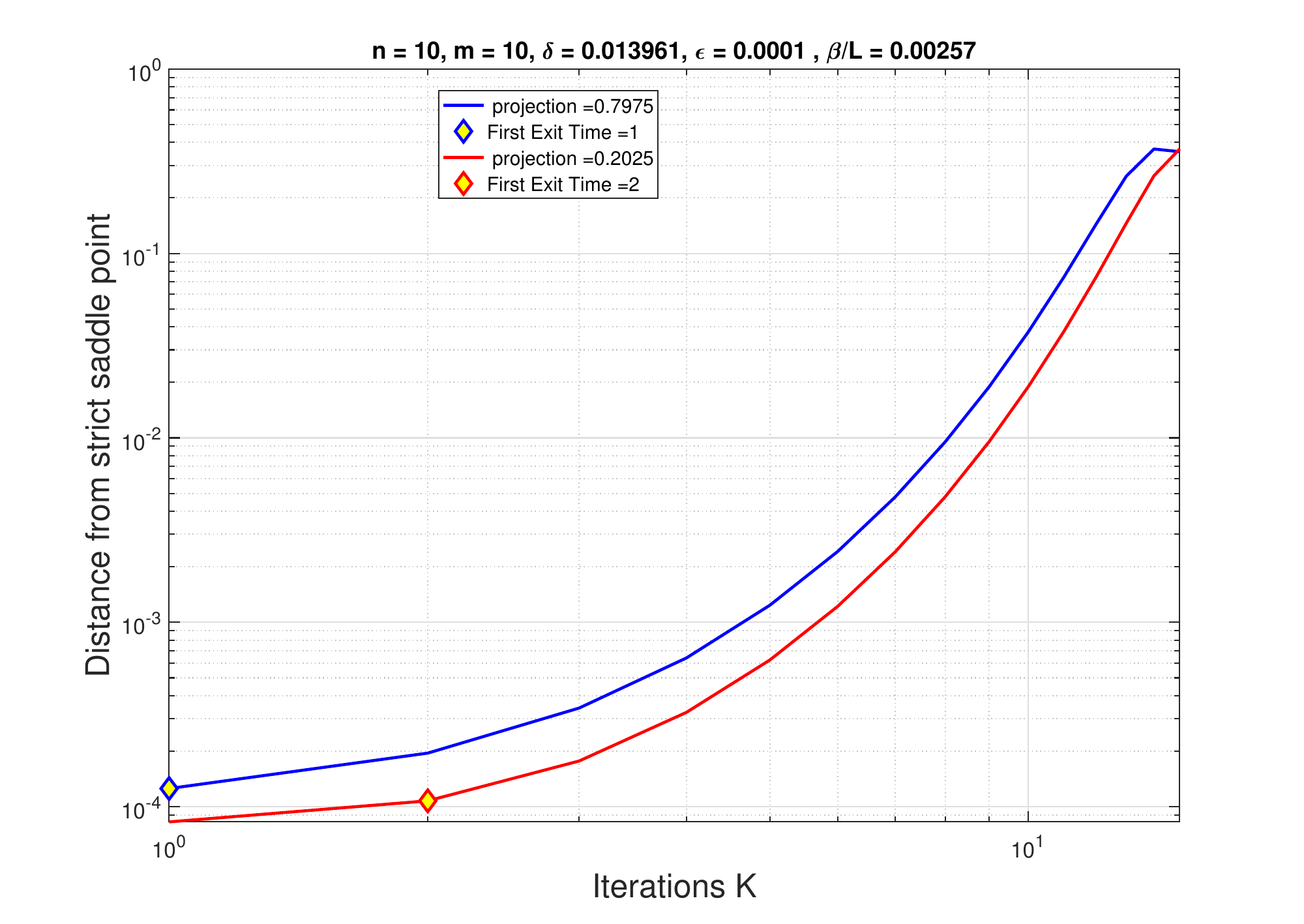} \\
    (a) & (b) \\
   \hspace{-0.9cm}  \includegraphics[width=3in]{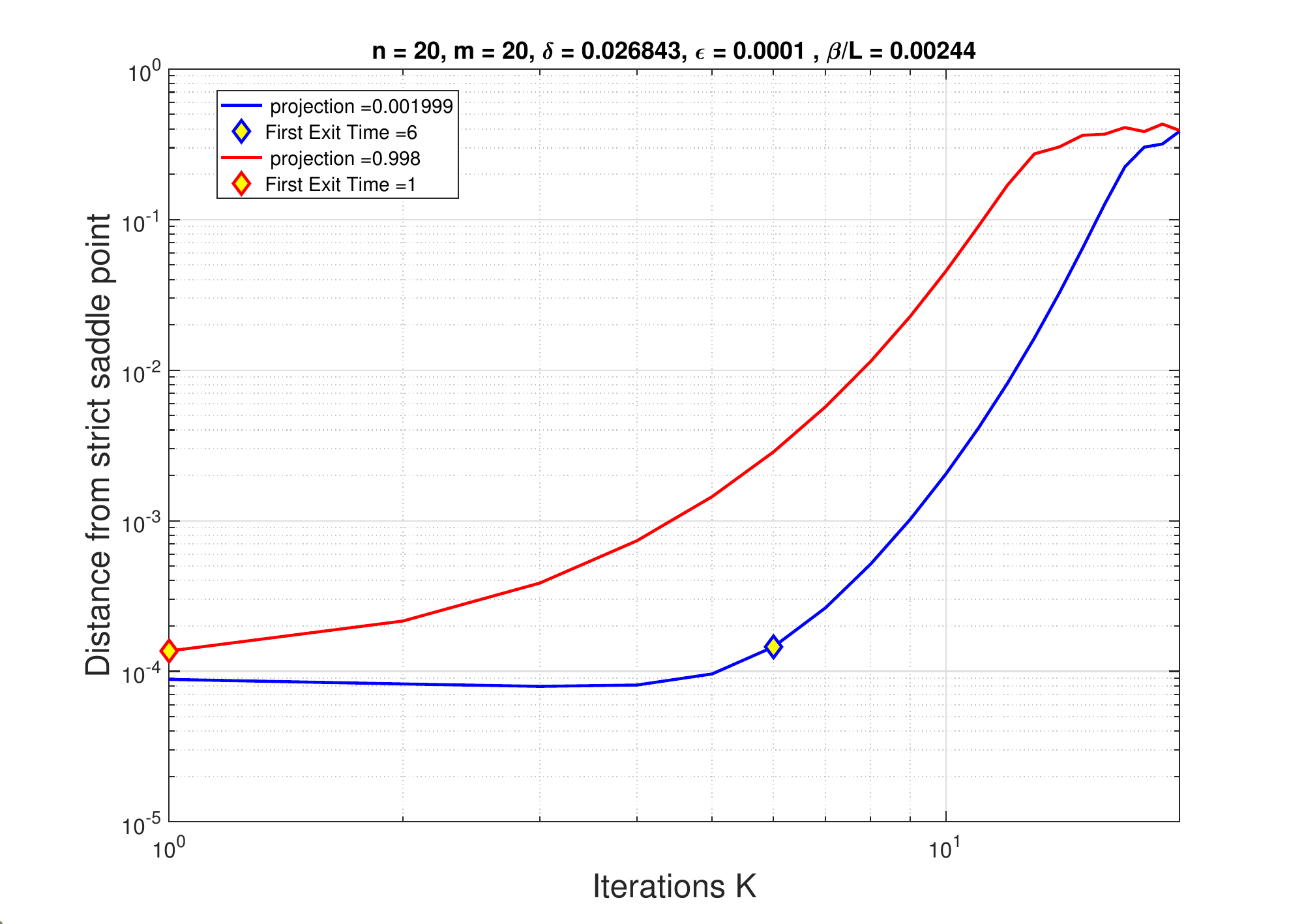} & \includegraphics[width=3in]{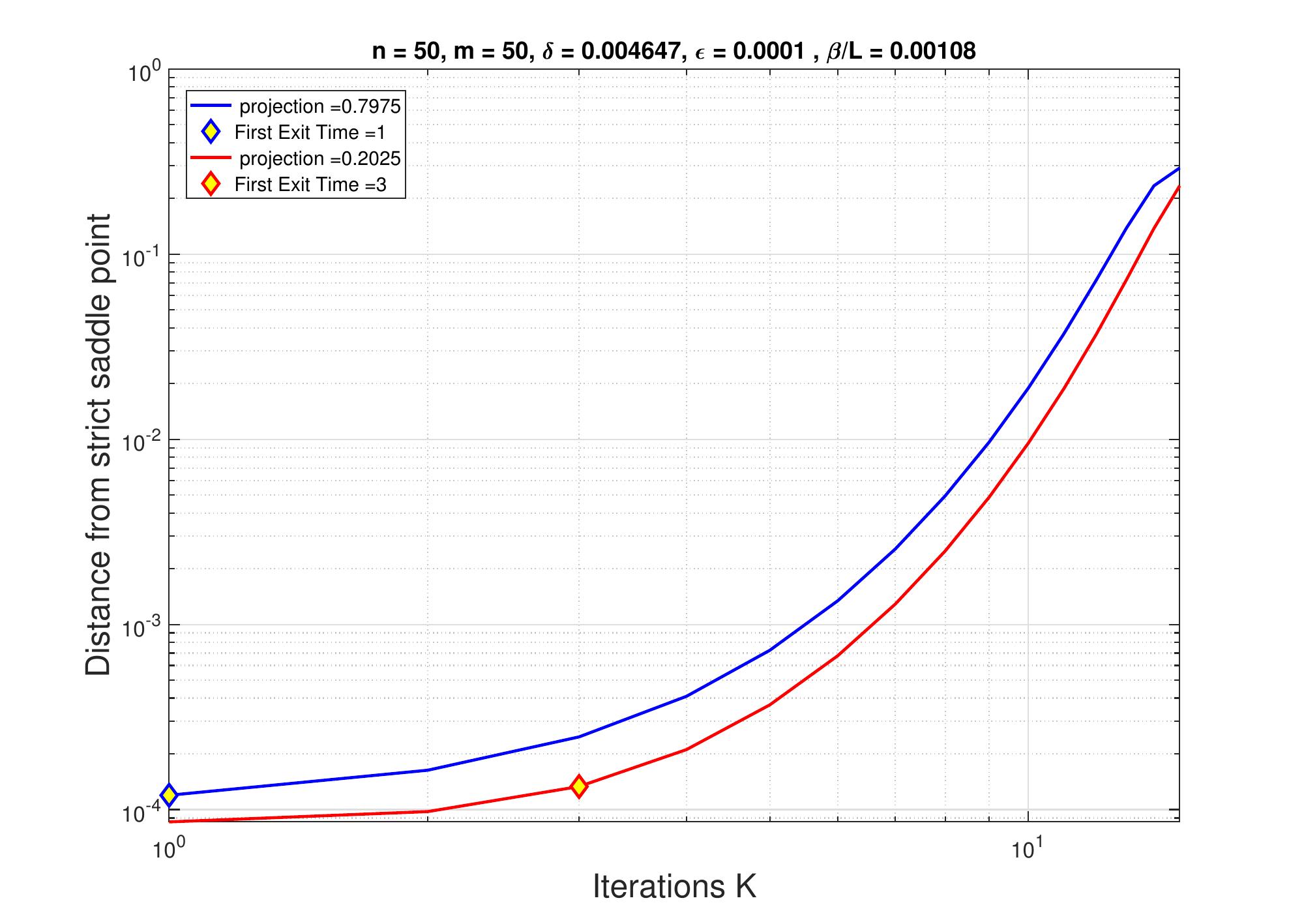} \\
    (c) & (d)
\end{tabular}
\caption{\revise{Simulating gradient trajectories on the phase retrieval problem with $\alpha = 1/L$ under certain initial unstable projections for various values of $m$, $n$ and $\epsilon$.}}
\label{fig12}
\end{figure}

In the simulations, we set $y_j = 1$ for $ 1 \leq j \leq \left \lfloor{\frac{m}{2}}\right \rfloor $ and $y_j = -1$ otherwise. Also, for the sake of simplicity we always set $m=n$ so that the system of equations $y_j = \langle \a_j, \x\rangle^2$ is neither under determined nor over determined and the Hessian of the function $f(\cdot)$ is full rank. The i.i.d.\ nature of the $\a_j$'s thus implies that the parameter $\frac{\beta}{L}$ is not too small \revise{and therefore Assumption \textbf{A4} gets satisfied}. The closed-form expressions for the gradient and the Hessian of the function in \eqref{phaseretrieval} are, respectively, as follows:
\begin{align}
    \nabla f(\x) &= \frac{1}{m} \sum\limits_{j= 1}^m \bigg(\langle\a_j,\x \rangle^2 - y_j \bigg) \langle\a_j,\x \rangle\a_j, \quad \text{and} \\
    \nabla^2 f(\x) &=  \frac{1}{m} \sum\limits_{j= 1}^m \bigg(3\langle\a_j,\x \rangle^2 - y_j \bigg) \a_j\a_j^T. \label{eigenphase}
\end{align}
For the particular choice of $y_j'$s it is observed that $\x^*= \mathbf{0} $ is a strict saddle point. We now initialize the gradient descent method in the $\epsilon$-neighborhood of $\x^*$ and examine the exit-time behavior of its trajectories for different values of $n,m,\epsilon,$ and the `projection' of the initial iterate on the unstable subspace, which corresponds to the quantity $ \sum_{j \in \mathcal{N}_{US}} (\theta_j^{us})^2$. The results are reported in Figure~\ref{fig11} for the step size of $\alpha = 0.1/L$ and in Figure~\ref{fig12} for the step size of $\alpha = 1/L$, with $L$ being the largest eigenvalue of $ \nabla^2 f(\mathbf{\x^*})$. Note that each subplot in both of the figures corresponds to different random $\a_j$'s. In order to highlight the dependence of the exit time on the unstable projection, we compare two different initializations of the gradient descent method for the same set of problem parameters in terms of the radial distance of the respective generated trajectories from the saddle point. Also the "first exit time" (the iteration when the gradient trajectory exits $\mathcal{B}_{\epsilon}(\x^*)$ for the first time) from the saddle neighborhood for the two trajectories are marked on each of the curves in colors matching with their respective radial distance curves.

It is evident from the two figures that, as suggested by the theoretical developments in this paper, a larger initial unstable subspace projection results in a faster exit time. More importantly, Figure \ref{fig12} corroborates our findings from Theorem \ref{theorem2} that for the step size of $\frac{1}{L}$, even with very small initial unstable subspace projections, i.e., $  \sum_{j \in \mathcal{N}_{US}} (\theta_j^{us})^2 = \mathcal{O}(\epsilon)$ such as those in Figure \ref{fig12}(a) and Figure \ref{fig12}(b), faster exit times are possible. Such conclusion does not necessarily hold for small step size, as in Figure \ref{fig11}(a) and Figure \ref{fig11}(b), where small initial unstable subspace projections yield relatively larger exit times.

\revise{We next illustrate the dependence of the exit time estimate on the dimension $n$ and eigen gap $\delta$. We first develop a numerical setup to showcase the dependence on $\delta$. To give an idea of our experimental setup, below is a step-by-step methodology used to perform simulations}:
\begin{enumerate}
\item \revise{Suppose the Hessian of function $f(\cdot)$ for the phase retrieval problem \eqref{phaseretrieval} has three distinct groups of eigenvalues,\footnote{\revise{We can introduce more groups of eigenvalues but refrain from doing so for the sake of simplicity.}} where the eigenvalues within any group are identical such that one group has eigenvalues equal to the gradient Lipschitz constant $L$ (as before, $L$ is the largest eigenvalue of $\nabla^2 f(\mathbf{\x^*})$, where $\x^*$ is the strict saddle point), the other group has eigenvalues equal to $ -\beta$, and the third group is placed on the eigenvalue spectrum so that it is at a $\delta$ distance from one of these groups. Further, suppose the third eigenvalue group has eigenvalues $-\beta + \delta$ where $(L+\beta)/2 > \delta > 2 \beta$. This construction preserves the parameters $L, \beta$ from Assumptions \textbf{A2}, \textbf{A4} for the function $f(\cdot)$ as the eigen gap $\delta$ is varied. Though the Hessian Lipschitz parameter $M$ for the function $f(\cdot)$ may not be preserved by this construction,\footnote{\revise{The Hessian Lipschitz parameter $M$ may change but will remain bounded in every compact set and therefore will be upper bounded by a constant term in the ball $\mathcal{B}_{\epsilon}(\x^*)$. Also, $M$ will remain constant with respect to the dimension $n$ since $n$ is fixed here.}} yet this setup is still able to control a given maximum number of parameters, i.e., $L, \beta$ and the problem dimension $n$.}
\item \revise{Next, we set $m=n=100$ in the phase retrieval problem \eqref{phaseretrieval}, where $\x \in \mathbb{R}^n$, $\a_j$'s are taken to be the canonical basis of $ \mathbb{R}^n$, and the eigen gap $\delta$ varies in the range $[0.15,2.13]$. Using the setup described in the previous bullet point, we then set the $y_j$'s as follows:}
    \begin{align}
         \revise{y_j =\begin{cases}
           \frac{m}{20}  \hspace{0.1cm}\textit{ ;} \hspace{0.1cm}  &1 \leq j \leq \left \lfloor{\frac{m}{3}}\right \rfloor  \\
              \frac{m}{20} - m \delta  \hspace{0.1cm}\textit{ ;} \hspace{0.1cm} & \left \lfloor{\frac{m}{3}}\right \rfloor+1 \leq j \leq 2\left \lfloor{\frac{m}{3}}\right \rfloor \hspace{0.1cm} \label{setup1} \\
             -5 m \hspace{0.1cm}\textit{ ;} \hspace{0.1cm} & \textit{otherwise}.
        \end{cases}}
    \end{align}
\item \revise{Since the $\a_j$'s are orthonormal, it can be readily checked using \eqref{eigenphase} that the eigenvalues of $\nabla^2 f(\x)$ at $\x = \mathbf{0}$ are $-y_j/m$ and we have $L = 5$, $\beta =  1/20$ from the given choice of $y_j$'s. By the above choice of $y_j$'s, the eigenvalues belong to three distinct groups and $\x = \mathbf{0}$ is a strict saddle point. In particular, the choice $ y_j = \frac{m}{20} - m \delta $ from above corresponds to the case where the free eigenvalue group has eigenvalues equal to $ ( \delta - \frac{1}{20})$.}
\item \revise{Finally, for the eigen gap $\delta$ in the range $[0.15,2.13]$, we compute the exit time from $\epsilon$-neighborhood of the origin for different values of the initial unstable subspace projections. }
\end{enumerate}

\begin{figure}[h]
\centering
\begin{tabular}{c}
    \includegraphics[width=4in]{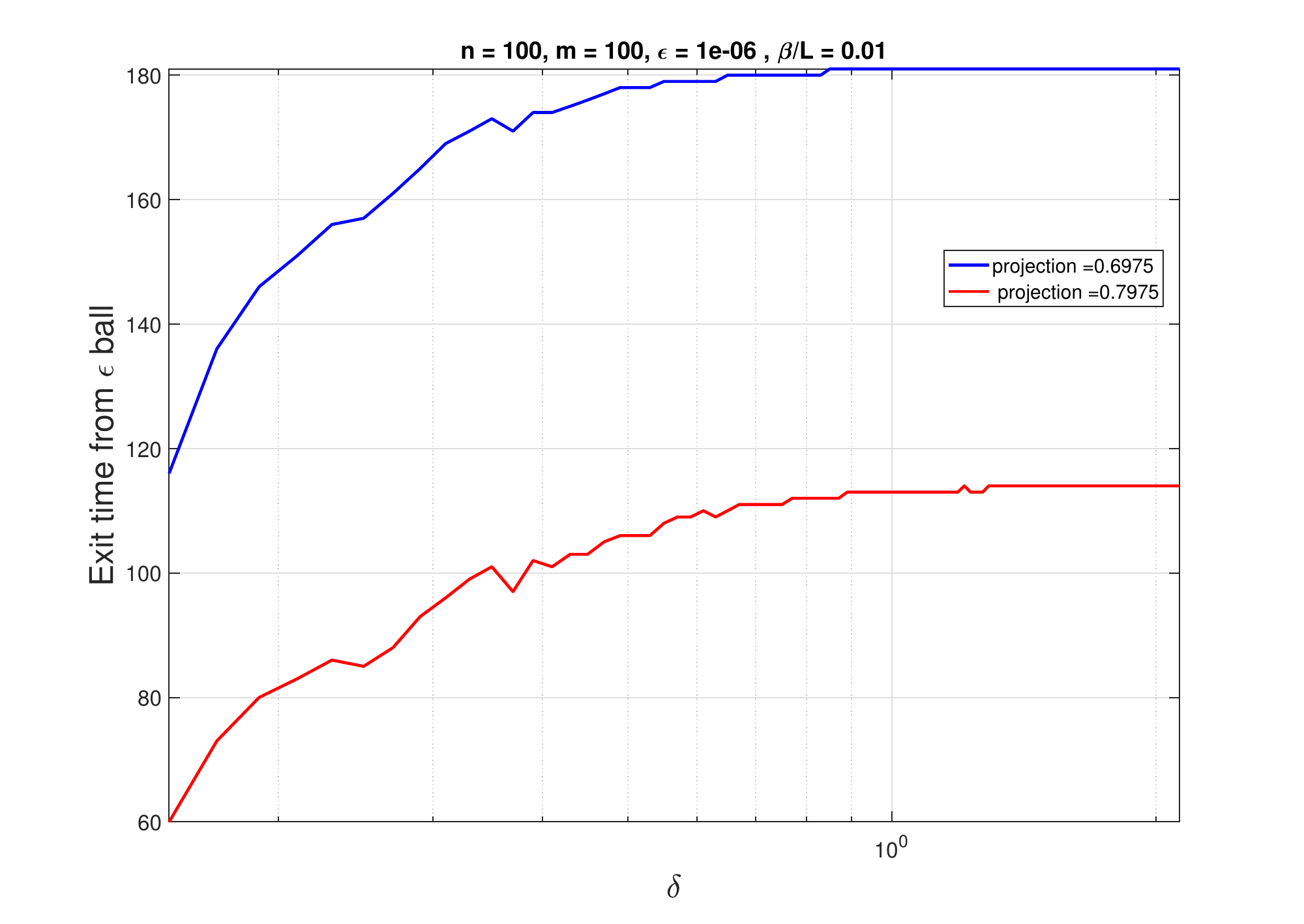}
\end{tabular}
\caption{\revise{Exit time versus the eigen gap $\delta$ (logarithmic scale) under certain initial unstable subspace projections for given values of $n$, $L$, $\beta$, and $\epsilon$.}}
\label{fig13}
\end{figure}

\revise{The results for this numerical setup are plotted in Figure \ref{fig13} for two values of the initial unstable subspace projections for $\alpha = 0.1/L$, where we have displayed the exit time versus $\delta$ on the logarithmic scale. We observe from the figure that the exit time increases with increasing eigen gap $\delta$ at least initially, which agrees with Theorem \ref{theorem2} where we have $K_{exit} \lessapprox \mathcal{O}(\log\delta) $.}

\revise{Next, we illustrate the dependence of the exit time on the problem dimension. Note that in general as the dimension $n$ increases, the gradient as well as the Hessian Lipschitz parameters ($L, M$) increase. In particular, we have {$L =\Theta(n)$, $M =\Theta(n)$} (see the discussion within Section~3 of \cite{chen2019gradient}). However, we can showcase the dependence of the exit time on the problem dimension for very particular choice of functions by keeping the gradient and Hessian Lipschitz parameters fixed with respect to the order of dimension. To this end, we modify the cost function in the phase retrieval problem \eqref{phaseretrieval} by normalizing it with dimension $n$ and rewriting \eqref{phaseretrieval} as:
\begin{align}
    \min_{\x \in \mathbb{R}^n} f(\x) = \frac{1}{4m n} \sum\limits_{j=1}^{m}\bigg[\langle\a_j,\x \rangle^2 - y_j \bigg]^2, \label{phaseretrievalnew}
\end{align}
where the normalization factor of $1/n$ helps in keeping the Hessian Lipschitz parameter independent of the dimension $n$. Note that in the earlier formulation \eqref{phaseretrieval} if we had {$M =\Theta(n)$} then in the new formulation \eqref{phaseretrievalnew} we will have {$M =\Theta(1)$}.}

\begin{figure}[h]
\centering
\begin{tabular}{c}
    \includegraphics[width=4in]{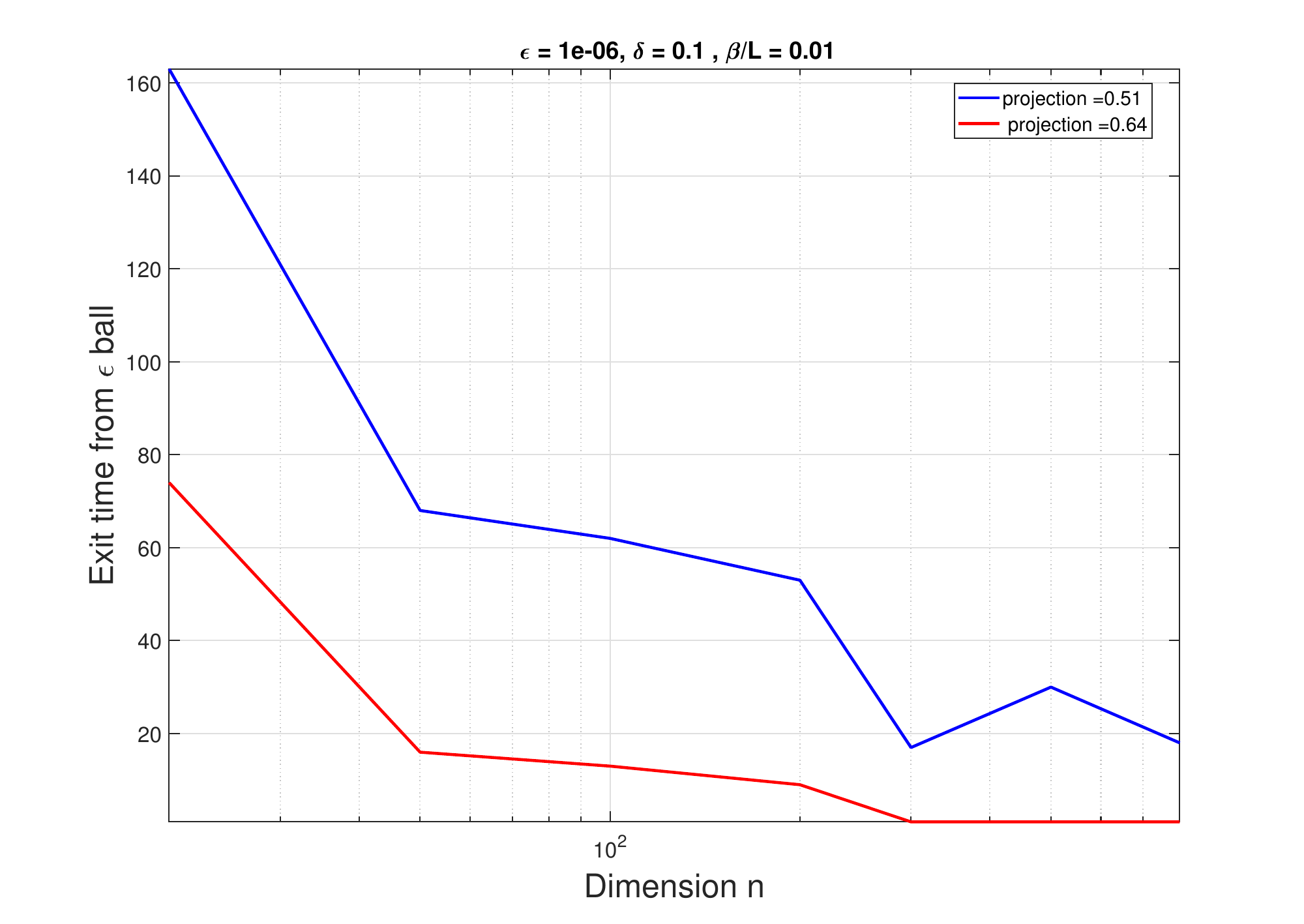}
\end{tabular}
\caption{\revise{Exit time vs dimension $n$ under certain initial unstable subspace projections for given values of $\delta$, $L$, $\beta$ and $\epsilon$.}}
\label{fig13a}
\end{figure}

\revise{Next, we once again set $m=n$ in \eqref{phaseretrievalnew}, where $\x \in \mathbb{R}^n$, and vary $n$ in the interval $ [20, 800]$. As before the $\a_j$'s are the canonical basis of $ \mathbb{R}^n$, while the eigen gap $\delta$ is fixed at $0.1$. We then set the $y_j$'s as follows:
    \begin{align}
         y_j =\begin{cases}
           \frac{m n}{20}  \hspace{0.1cm}\textit{ ;} \hspace{0.1cm}  &1 \leq j \leq \left \lfloor{\frac{m}{2}}\right \rfloor  \\
              -\frac{m n}{20}   \hspace{0.1cm}\textit{ ;} \hspace{0.1cm} & \left \lfloor{\frac{m}{2}}\right \rfloor+1 \leq j \leq 2\left \lfloor{\frac{m}{2}}\right \rfloor -1\hspace{0.1cm} \label{setup2} \\
             -5 m n \hspace{0.1cm}\textit{ ;} \hspace{0.1cm} & \textit{otherwise}.
       \end{cases}
    \end{align}
Since the $\a_j$'s are orthonormal, it can be readily checked after adapting \eqref{eigenphase} for the modified formulation \eqref{phaseretrievalnew} that the eigenvalues of $\nabla^2 f(\x)$ at $\x = \mathbf{0}$ are $\frac{-y_j}{m n}  $ and we have $L = 5$, $\beta =  1/20$, and $\delta = 0.1$ from the given choice of $y_j$'s. This construction preserves the parameters $L, \beta$ from Assumptions \textbf{A2}, \textbf{A4} and the eigen gap $\delta$ from Proposition \ref{eigenprop} for the function $f(\cdot)$ as the problem dimension $n$ is varied (the parameter $M$ from Assumption \textbf{A3} also gets independent of the dimension $n$). Finally, for $n \in [20, 800]$ we compute the exit time from the $\epsilon$-neighborhood of origin for different values of initial unstable subspace projections. The results are plotted in Figure \ref{fig13a} for two values of initial unstable subspace projections for $\alpha = 0.1/L$, where we have displayed the exit time versus dimension $n$ on the logarithmic scale. We observe that the exit time decreases with increasing {dimension $n$}, which agrees with Theorem \ref{theorem2} where we have $K_{exit} \lessapprox \mathcal{O}(\log n^{-1}) $.}

\subsection{\revise{Evolution of the trajectory function $\Psi(K)$ from Theorem \ref{theorem1} on phase retrieval problem}}\label{psik_simul}
\begin{figure}[h]
\centering
\begin{tabular}{c}
    \includegraphics[width=4in]{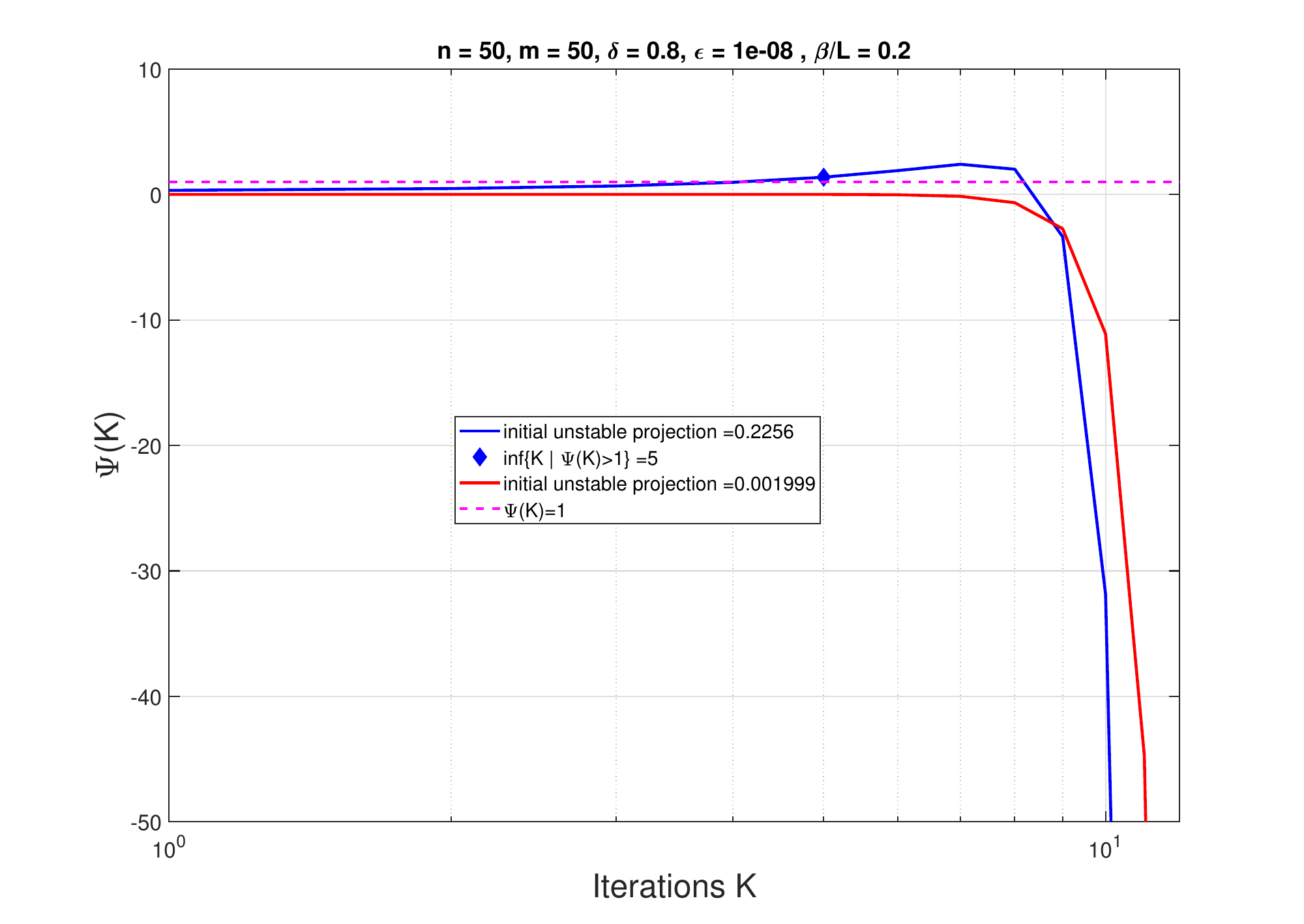}
\end{tabular}
\caption{\revise{$\Psi(K)$ vs $K$ under certain initial unstable subspace projections for given values of $n$, $L$, $\beta$, $\delta$ and $\epsilon$. The blue curve has a sufficient initial unstable subspace projection that allows it to first increase, become greater than $1$ and then decrease whereas the red curve always remains below $\Psi(K)=1$ and keeps on decreasing since it has a very small initial unstable subspace projection.}}
\label{fig16}
\end{figure}

\revise{We now illustrate that the trajectory function $\Psi(K)$ first increases to a maximum and then continuously decreases to $-\infty$ from the example of the phase retrieval problem. In particular, if the initial unstable subspace projection is not too small then there exists a non-trivial finite $K$ where $\Psi(K)>1 $, which is the upper bound on the exit time. In the phase retrieval problem \eqref{phaseretrieval} we set $m=n=20$, where $\x \in \mathbb{R}^n$, the $\a_j$'s are taken to be the canonical basis of $ \mathbb{R}^n$, the eigen gap $\delta=0.5$, $L=20 $, and $\beta = 2$. We then set the $y_j$'s as follows:
\begin{align}
         y_j =\begin{cases}
           m (\beta + \delta)  \hspace{0.1cm}\textit{ ;} \hspace{0.1cm}  & j=1  \\
             m \beta  \hspace{0.1cm}\textit{ ;} \hspace{0.1cm} & j=2 \\
             -m \beta \hspace{0.1cm}\textit{ ;} \hspace{0.1cm} & 3\leq j \leq m-1 \\
             -m L \hspace{0.1cm}\textit{ ;} \hspace{0.1cm} & \textit{otherwise}.
       \end{cases} \label{setup1_a}
\end{align}
The results are plotted in Figure \ref{fig16} for two values of initial unstable subspace projections for $\alpha = 1/L$. Clearly, the trajectory function $\Psi(K)$ first increases to a maximum and then continuously decreases to $-\infty$, which agrees with Theorem \ref{theorem1} (see also the discussion following Remark \ref{remark-trajectory-func}).}

\section{Conclusion}
This work has focused on the analysis of gradient-descent trajectories in some small neighborhoods of a strict saddle point. Using tools from matrix perturbation theory and first-order eigenvector perturbations, a proof technique has been developed that describes the behavior of gradient-descent method as a function of the local geometry around a strict saddle point. Two novel lemmas have been presented in this work that quantify the radius of a saddle neighborhood within which an approximate analysis for the gradient-descent trajectory can be developed, provided the trajectory stays inside this neighborhood for a bounded interval. Next, this work has also presented two novel theorems that quantify this approximate trajectory distance from the saddle point at every iteration and provide an exit time from the saddle neighborhood based on the initial unstable projection of the radial vector. Developing a robust algorithm that can leverage this analysis so as to efficiently escape saddle neighborhood and a rigorous analysis of the trajectory function are some of the directions that have been pursued in a follow-up paper \cite{dixit2022boundary} to this work.

\section{Data Availability Statement}
The data underlying this paper are available in the paper and in its online supplementary material.

\section*{Acknowledgments}
This work was supported in part by the National Science Foundation under grants CCF-1453073, CCF-1907658, CCF-1910110, OAC-1940209, CNS-2148104, CCF-1814888, and DMS-2053485, by the Army Research Office under grants W911NF-17-1-0546 and W911NF-21-1-0301, by the Office of Naval Research Award Number N00014-21-1-2244 and by the DARPA Lagrange Program under ONR/SPAWAR contract N660011824020. The authors would also like to thank H.\ Vincent Poor, an anonymous reader, and the reviewers for their careful reading and many helpful suggestions that have helped improve the paper.

\section*{Appendices}
\begin{appendices}
\section{\revise{On the equivalence of \eqref{eq-Kexit-1} and \eqref{exittimecrudeanalysis}}}\label{app:equivalence}
\begin{lemma}\label{lemma-exit-times-equivalent} \revise{In the setting of Section \ref{crudeanalysissection}, the exit time \eqref{eq-Kexit-1}
is equivalent to 
 \eqref{exittimecrudeanalysis}.}
\end{lemma}
\begin{proof}
\revise{First, we show that the condition $\norm{\x_k - \x^*} = 0$ does not hold for any finite $k\geq 0$. For $k=0$, this is a trivial statement as our initialization is such that $\norm{\x_k - \x^*} = \epsilon>0$. We then proceed by induction. Suppose that $ \norm{\x_{k-1} - \x^*} >0$ and $ \norm{\x_k - \x^*} = 0$ for some finite $k\geq 1$. Since $ \x_k - \x^* = \x_{k-1} - \x^* - \alpha \nabla f(\x_{k-1})$, we can write $\norm{\x_k - \x^*} = \norm{ (\mathbf{I} - \alpha \textbf{M})(\x_{k-1} - \x^*)} $ with $\textbf{M} = \int_{p=0}^1 \nabla^2 f(\x^* + p(\x_{k-1} - \x^*))dp,$ where we have used Taylor's formula (with an integral form) to represent the gradient of $f$ as an integral over the Hessian of $f$. By assumption, we have $\alpha \leq \frac{1}{L}$ and we first consider the case of $\alpha < \frac{1}{L}$ so that $\alpha\norm{\textbf{M}}_2 <1$, which implies $ \norm{(\mathbf{I} - \alpha \textbf{M})^{-1}}_2^{-1} > 0$ and we can therefore write $(\mathbf{I} - \alpha \textbf{M})^{-1} (\x_k - \x^*) = \x_{k-1} - \x^*$. Then, $\norm{\x_k - \x^*} \geq \norm{(\mathbf{I} - \alpha \textbf{M})^{-1}}_2^{-1}\norm{(\x_{k-1} - \x^*)} >0$, which leads to a contradiction. Therefore, we conclude that $\norm{\x_k - \x^*}>0$ for every $k$. Correspondingly, the quantity $\gamma_k = \frac{ \langle\v_n, (\x_k - \x^*)\rangle}{\norm{\x_k - \x^*}}$ is well-defined in the sense that its denominator cannot vanish. Here,  $\gamma_k \in [0,1]$ because
 the vectors $\v_n$ and $\frac{\x_k - \x^*}{\norm{\x_k - \x^*}}$ are both unit vectors and if the dot product is negative, we can always flip the sign of the eigenvector $\v_n$.
Note that throughout this crude analysis section, for the sake of simplicity, we assume the dot product does not vanish, i.e., $\gamma_k \neq 0$ for any $k$, because otherwise the set $\{k \vert  \langle\v_n, (\x_k - \x^*)\rangle > \gamma_k \epsilon\}$ can be empty.\footnote{\revise{This assumption would be satisfied for instance for quadratic objectives if the initialization has a non-zero component in the stable subspace of the Hessian at the saddle point; this can be verified as the solutions admit an explicit formula for every $k$ in the quadratic case.}}} 

\revise{Next, notice that by the definition of $\gamma_k$, we have $ \langle\v_n, (\x_k - \x^*)\rangle > \gamma_k \epsilon $ $ \iff$ $ \norm{\x_k - \x^*} > \epsilon$; this is because by multiplying the latter inequality with the positive scalar $\gamma_k$, we can simply obtain the former inequality. Therefore, we conclude that the sets $\{k \vert \norm{\x_k - \x^*} > \epsilon\}$ and $\{k \vert  \langle\v_n, (\x_k - \x^*)\rangle > \gamma_k \epsilon\}$ (defined in (3.19) and (3.20) respectively) are identical for $ \gamma_k\in (0,1]$ and $\alpha < \frac{1}{L}$. When $\alpha = \frac{1}{L}$, we can have $ \norm{\x_k -\x^*} = 0$ for some finite $k=K$, but since $\x^*$ is a fixed point of the gradient descent iteration, we will get $ \norm{\x_k -\x^*} = 0$ for all $k> K$, which implies $\inf_{k>0}\{k \vert \norm{\x_k - \x^*} > \epsilon\} = \infty$. Since we are looking for finite exit times in the crude analysis, we can disregard the case of $ \norm{\x_k -\x^*} = 0$ for some finite $k$ when $\alpha  = \frac{1}{L}$, and then for $ \gamma_k\in (0,1]$, we again conclude that the sets
$\{k \vert \norm{\x_k - \x^*} > \epsilon\}$ and $\{k \vert  \langle\v_n, (\x_k - \x^*)\rangle > \gamma_k \epsilon\}$ are identical. Therefore, we conclude that (3.19) and (3.20) are equivalent. This completes the proof.}
\end{proof}

\section{Proof of Lemma \ref{lem4} (Hessian perturbation)} \label{Appendix A}
    \proof{}
	From the Taylor expansion around the strict saddle point $\x^*$ along the direction $\x_k-\x^*$ we have the following:
	\begin{align}
	\nabla f(\x_k) & = \nabla f(\x^*) + \int_{p=0}^{p=1}\nabla^2 f(\x^* + p\u_k) \u_k dp \\
	\implies    \nabla f(\x_k) & = \nabla f(\x^*) + \int_{p=0}^{p=1}\nabla^2 f(\x^* + p\norm{\u_{k}}\hat{\u}_k) \u_k dp, \label{integralbound}
	\end{align}
	where $\u_{k} = \x_{k} - \x^*$ and $\{\x_{k}\}$ is the sequence of iterates generated from the gradient descent method \eqref{gd}.
	
	Note that here in the last step, we have made the substitution of $\u_k = \norm{\u_{k}} \hat{\u}_k$ and we have that $\norm{\u_k} \leq \epsilon$ since our iterate $\x_k$ lies inside the ball $\mathcal{B}_{\epsilon}(\x^*)$. $\hat{\u}_k$ represents the unit vector in the direction of $\u_k$.

Next, we start developing the term $\nabla^2 f(\x^* + p\norm{\u_{k}}\hat{\u}_k)$ using matrix perturbation theory and variational calculus. We start with introducing a matrix function $\G(\cdot) : \mathbb{R} \rightarrow \mathbb{R}^{n \times n}$ which is given by
	\begin{align}
	\G(w) = \nabla^2 f(\x^* + w\hat{\u}_k),
	\end{align}
	where $w= p \norm{\u_{k}}$, $p$ being the variable of previous integration and therefore $w = \mathcal{O}(\epsilon) $. For sufficiently small $\epsilon$ we can utilize the Taylor series expansion of $\G(w)$ around $w=0$:
	\begin{align}
	\hspace{-3cm}\G(w) &= \G(0) + w\frac{d\G}{dw}\bigg\vert_{w=0} + \frac{w^2}{2} \frac{d^2\G}{dw^2}\bigg\vert_{w=0} + \hspace{0.5cm}+ \dots \\
	\implies  \nabla^2 f(\x^* + w\hat{\u}_k) &= \underbrace{\nabla^2 f(\x^*) + w \frac{d}{dw}(\nabla^2 f(\x^* + w\hat{\u}_k))\bigg\vert_{w=0}}_{S_1} + \underbrace{\frac{w^2}{2} \frac{d^2}{dw^2}(\nabla^2 f(\x^* + w\hat{\u}_k))\bigg\vert_{w=0}  + \hspace{0.5cm}+ \dots \hspace{0.5cm}}_{R_1} .\label{taylorrem}
	\end{align}
	With $w = \mathcal{O}(\epsilon) $ and the eigenvalues of $\nabla^2 f(\x^*)$ separated by $\delta$ or more, we can get rid of all the higher-order terms in the Taylor sequence from $w^2$ onwards. It is a reasonable approximation from the Rayleigh--Schr\"{o}dinger perturbation theory (\cite{zhong2017eigenvector,MBT,MBT1}) as long as we have \revise{Proposition \ref{eigenprop}}, i.e., there are at least two eigenvalue groups of $\nabla^2 f(\x^*)$ that are not degenerate or too close to one another. This leaves us with the following first order approximation:
	\begin{align}
	\nabla^2 f(\x^* + w\hat{\u}_k) &= \nabla^2 f(\x^*) + w \frac{d}{dw}(\nabla^2 f(\x^* + w\hat{\u}_k))\bigg\vert_{w=0} + \mathcal{O}(\epsilon^2),
	\end{align}
	where we have that $S_1 = \nabla^2 f(\x^*) + w \frac{d}{dw}(\nabla^2 f(\x^* + w\hat{\u}_k))\bigg\vert_{w=0}$ and the order of the remainder term $R_1$ is $\mathcal{O}(\epsilon^2)$. This remainder term $R_1$ is easy to obtain from Taylor's Remainder theorem. Applying this theorem to \eqref{taylorrem} with the substitution $\nabla^2 f(\x^* + u\hat{\u}_k) = \G(u)$ yields
	\begin{align}
	R_1 &= \int_{0}^{w} u \frac{d^2 \G}{du^2} du \label{taylorremfirstorder} \\
	\implies \norm{R_1}_2 &= \norm{\int_{0}^{w} u \frac{d^2 \G}{du^2} du}_2 < \bigg(\int_{0}^{w} \norm{ \frac{d^2 \G}{du^2} }_2^2 du\bigg)^{\frac{1}{2}} \bigg(\int_{0}^{w} u^2 du \bigg)^{\frac{1}{2}} \leq \frac{B_2w^2}{\sqrt{3}} < \frac{B_2{\epsilon}^2}{\sqrt{3}}.
	\end{align}
	Here in the last step we applied the Cauchy-Schwarz inequality followed by an extra assumption on the spectral radius of $\frac{d^2 \G}{du^2}$ which is $\norm{\frac{d^2 \G}{du^2}}_2 \leq B_2$ for some finite positive value $B_2$. The final inequality follows from the fact that $w = p\norm{\u_{k}} < \epsilon$ where $0<p\leq1$. Hence the remainder term $R_1$ is of order $\mathcal{O}(\epsilon^2)$. Note that the condition of $\norm{\frac{d^2 \G}{du^2}}_2 \leq B_2 < \infty$ is valid for any analytic function $f(\cdot)$. Moreover, it bounds the variations of the Hessian inside the ball $\mathcal{B}_{\epsilon}(\x^*)$.

Next, using a matrix substitution of $\H(\hat{\u}_k) = \frac{d}{dw}(\nabla^2 f(\x^* + w\hat{\u}_k))|_{w=0}$, our first order Hessian approximation becomes
	\begin{align}
	\nabla^2 f(\x^* + w\hat{\u}_k) &= \nabla^2 f(\x^*) + w \H(\hat{\u}_k)  + \mathcal{O}(\epsilon^2) \label{firstorderperturbation}\\
	\implies \nabla^2 f(\x^* + p{\u}_k) &= \nabla^2 f(\x^*) + p\norm{\u_{k}} \H(\hat{\u}_k)  + \mathcal{O}(\epsilon^2). \label{hessianapprox}
	\end{align}
	
	\subsection{Rayleigh--Schr\"{o}dinger perturbation analysis}
	We can now find the matrix $\H(\hat{\u}_k)$ using the spectral theorem and the Rayleigh--Schr\"{o}dinger perturbation theory. Note that this matrix $\H(\hat{\u}_k)$ depends on the unit vector $\hat{\u}_k$.
	
	We first apply the spectral theorem on the real symmetric matrix $\nabla^2 f(\x^* + w\hat{\u}_k) $ to get the following decomposition in terms of its eigenvalues $\lambda_i(w)$ and the eigenvectors $\v_i(w)$:
	\begin{align}
	\nabla^2 f(\x^* + w\hat{\u}_k)  &= \sum_{i=1}^{n} \lambda_i(w)\v_i(w)\v_i(w)^{T}.
	\end{align}
	Now, differentiating this decomposition w.r.t.\ the variable $w$ and obtaining its value at the point $w=0$ we get
	\begin{align}
	\frac{d}{dw}(\nabla^2 f(\x^* + w\hat{\u}_k)) &=  \sum_{i=1}^{n} \frac{d}{dw}(\lambda_i(w)\v_i(w)\v_i(w)^{T})  \\
	\implies  \frac{d}{dw}(\nabla^2 f(\x^* + w\hat{\u}_k))\bigg\vert_{w=0} &=  \sum_{i=1}^{n} \bigg(\frac{d}{dw}(\lambda_i(w))\bigg\vert_{w=0}\v_i(0)\v_i(0)^{T} + \lambda_i(0)\frac{d}{dw}(\v_i(w))\bigg\vert_{w=0}\v_i(0)^{T} \nonumber \\ &+ \lambda_i(0)\v_i(0)\frac{d}{dw}(\v_i(w)^{T})\bigg\vert_{w=0}\bigg) \\
	\implies  \H(\hat{\u}_k) &=  \sum_{i=1}^{n} \bigg(\frac{d}{dw}(\lambda_i(w))\bigg\vert_{w=0}\v_i(0)\v_i(0)^{T} + \lambda_i(0)\frac{d}{dw}(\v_i(w))\bigg\vert_{w=0}\v_i(0)^{T} \nonumber \\ &+ \lambda_i(0)\v_i(0)\frac{d}{dw}(\v_i(w)^{T})\bigg\vert_{w=0}\bigg). \label{spectral}
	\end{align}
Note that the pair $(\lambda_i(0),\v_i(0))$ represents the $i^{th}$ eigenvalue-eigenvector pair of the unperturbed matrix $ \nabla^2 f(\x^*)$. From the Rayleigh--Schr\"{o}dinger perturbation theory (\cite{zhong2017eigenvector}), for a given first order perturbation matrix $\H(\hat{\u}_k)$ in \eqref{firstorderperturbation}, we have the following first order correction terms:
	\begin{align}
	\frac{d}{dw}(\lambda_i(w))\bigg\vert_{w=0} & = \langle \v_i(0), \H(\hat{\u}_k)\v_i(0)\rangle  \\
	\frac{d}{dw}(\v_i(w))\bigg\vert_{w=0} & = \sum_{l \neq i} \frac{\langle \v_l(0), \H(\hat{\u}_k)\v_i(0)\rangle}{\lambda_{i}(0)-\lambda_{l}(0)}\v_l(0).  \label{subsectioncorrectionterm}
	\end{align}
Observe that under \revise{Proposition \ref{eigenprop}}, we are considering the case of $m=n$, i.e., no degenerate eigenvalues in our analysis. However, we have a subsection after Lemma \ref{lem4} (generality of Lemma \ref{lem4}) that discusses the degenerate case as well.
	
	Substituting these first-order correction terms in \eqref{spectral}, we get the following result:
	\begin{align}
	\H(\hat{\u}_k) &=  \sum_{i=1}^{n} \bigg(\langle \v_i(0), \H(\hat{\u}_k)\v_i(0)\rangle \v_i(0)\v_i(0)^{T} + \lambda_i(0)(\sum_{l \neq i} \frac{\langle \v_l(0), \H(\hat{\u}_k)\v_i(0)\rangle}{\lambda_{i}(0)-\lambda_{l}(0)}\v_l(0))\v_i(0)^{T} \nonumber \\ &+ \lambda_i(0)\v_i(0)(\sum_{l \neq i} \frac{\langle \v_l(0), \H(\hat{\u}_k)\v_i(0)\rangle}{\lambda_{i}(0)-\lambda_{l}(0)}\v_l(0))^{T}\bigg). \label{errorbound}
	\end{align}
	Now, combining this result with \eqref{spectral} and substituting the subsequent matrix approximation in \eqref{integralbound} leads to the following result:
	\begin{align}
	\nabla f(\x_k) & =  \nabla f(\x^*) + \int_{p=0}^{p=1}(\nabla^2 f(\x^*) + p \norm{\u_{k}} \H(\hat{\u}_k) + \mathcal{O}(\epsilon^2)) \u_k dp \\
	& = \bigg(\nabla^2 f(\x^*) +  \frac{\norm{\u_{k}}}{2} \H(\hat{\u}_k)+ \mathcal{O}(\epsilon^2) \bigg) \u_k. \label{gradapprox}
	\end{align}
	Note that $\norm{\u_{k}} \H(\hat{\u}_k)$ and $\u_{k}$ do not depend on $p$ and hence can be pulled out of the integral.

	\subsection{Validity of the Taylor expansion in Rayleigh--Schr\"{o}dinger analysis}
	
	Recall that we used the Taylor expansion in \eqref{taylorrem} for the matrix $\G(w)$ around $w=0$. Next, we evaluated the first-order perturbation term $\H(\hat{\u}_k)$ in this expansion using the Rayleigh--Schr\"{o}dinger perturbation analysis, which is dependent on this Taylor expansion (see derivations in \cite{MBT,MBT1}). In other words, the perturbation analysis is only valid for those values of $w$ where the Taylor expansion for the matrix $\G(w)$ around $w=0$ converges. This directly reduces to the problem of finding the radius of convergence for the expansion \eqref{taylorrem}.
	
	Although evaluating the radius of convergence in the Rayleigh--Schr\"{o}dinger perturbation analysis remains an open problem in general, we can still find the radius of convergence for the expansion \eqref{taylorrem} using matrix power series.
	
	For the Taylor expansion in \eqref{taylorrem}, consider the sequence $\{r_j(\hat{\u}_k)\}$ for all $j\in \{1,2,...\}$ such that
	\begin{align}
	r_j(\hat{\u}_k) &= \norm{\bigg(\frac{d^j \G}{dw^j}\bigg\vert_{w=0}\bigg)}_2, \label{spectral11}
	\end{align}
	where $\G(w) = \nabla^2 f(\x^*+w\hat{\u}_k)$ and $w = p\norm{\u_{k}}$ with $0<p\leq 1$.
		
	Next by the Cauchy--Hadamard theorem, for any power series defined by
	\begin{align}
	h(z) &= \sum_{j=0}^{\infty} c_j (z-a)^j
	\end{align}
	where $z \in \mathbb{C}$, the radius of convergence for the series is given by
	\begin{align}
	r &= \bigg(\limsup_{j\to \infty} \sqrt[j]{\abs{c_j}}\bigg)^{-1}.
	\end{align}
	
For the case of matrix power series, the spectral radius of a matrix is used to determine the radius of convergence. From the expression of the $r_j(\hat{\u}_k)$ in \eqref{spectral11}, it is clear that the matrix $\frac{d^j \G}{dw^j}\bigg\vert_{w=0}$ is real-symmetric since $\G$ is real-symmetric. Hence, the spectral radius of this matrix is equal to its $l_2$ norm.
 	
	Using the Cauchy--Hadamard theorem on our expansion \eqref{taylorrem} for $\abs{c_j} = \frac{r_j(\hat{\u}_k)}{j!}$, we get the following radius of convergence:
	\begin{align}
	r(\hat{\u}_k) &= \bigg(\limsup_{j\to \infty} \sqrt[j]{\frac{r_j(\hat{\u}_k)}{j!}}\bigg)^{-1}.
	\end{align}
	Therefore, if $\sqrt[j]{\frac{r_j(\hat{\u}_k)}{j!}}$ is upper bounded for all $j$, then a non-zero radius of convergence is guaranteed. This implies that
	\begin{align}
	w &= p \norm{\u_{k}} <  \bigg(\limsup_{j\to \infty} \sqrt[j]{\frac{r_j(\hat{\u}_k)}{j!}}\bigg)^{-1}. \label{ROC1}
	\end{align}	
	Since $w<\epsilon$ for any $\x_{k} \in \mathcal{B}_{\epsilon}(\x^*)$, where $\x_{k} = \x^* + w \hat{\u}_k$, by setting a condition on $\epsilon$ such that  $\epsilon < \bigg(\limsup_{j\to \infty} \sqrt[j]{\frac{r_j(\hat{\u}_k)}{j!}}\bigg)^{-1}$, we can guarantee the inequality \eqref{ROC1}.
	However this result should hold for any possible unit directional vector $\hat{\u}_k$. Hence we must have
	\begin{align}
	\epsilon &< \inf_{\hat{\u}_k}\bigg(\limsup_{j\to \infty} \sqrt[j]{\frac{r_j(\hat{\u}_k)}{j!}}\bigg)^{-1} \\
	\implies \epsilon &< \inf_{{\norm{\u}=1}}\bigg(\limsup_{j\to \infty} \sqrt[j]{\frac{r_j(\u)}{j!}}\bigg)^{-1}, \label{ROC}
	\end{align}
	where \begin{align}
	r_j(\u) &= \norm{\bigg(\frac{d^j }{dw^j}\nabla^2 f(\x^*+w\u)\bigg\vert_{w=0}\bigg)}_2.
	\end{align}
	 It is to be noted that this bound on $\epsilon$ only guarantees convergence of the expansion \eqref{taylorrem} and not the convergence of terms generated by the Rayleigh--Schr\"{o}dinger perturbation analysis. Evaluating the convergence radius from the Rayleigh--Schr\"{o}dinger perturbation theory is beyond the scope of the current work. Hence this condition on $\epsilon$ is necessary but may not be sufficient.
	
	 \subsection{Note on the existence of a positive upper bound on $\epsilon$}
	 For the condition \eqref{ROC} to make sense, we must have $\inf_{{\norm{\u}=1}}\bigg(\limsup_{j\to \infty} \sqrt[j]{\frac{r_j(\u)}{j!}}\bigg)^{-1} > 0$. To this end, consider the following Taylor expansion with respect to the variable $w \geq 0$:
	 \begin{align}
	     \nabla^2 f(\x^* + w \u) = \sum\limits_{j=0}^{\infty} \frac{d^j }{dw^j}\nabla^2 f(\x^* + w \u)\bigg\vert_{w=0} \frac{w^j}{j!}, \label{ROCseries1}
	 \end{align}
	 where the above matrix-valued series converges with some strictly positive \emph{radius of convergence} (ROC) $R$ (i.e., $w \leq R$) for all $\{\u : \|\u\|_2 = 1\}$ due to the analytic nature of $f(\cdot)$. Here, we focus on convergence of the series with respect to the operator (spectral) norm and note that for any $n$-dimensional symmetric matrix $\A$ we have  the inequality $\frac{1}{n}\max_{i,l}\{\abs{[\A]_{i,l}}\} \leq \frac{1}{n}\norm{\A}_F \leq\norm{\A}_2 \leq \norm{\A}_F$. Thus, if the matrix-valued series \eqref{ROCseries1} converges for $w \leq R$ in the spectral norm then the matrix sum on the right-hand side of \eqref{ROCseries1} must also element-wise converge for the same ROC $R$. 
	 For the $(i,l)^{th}$ element of $\nabla^2 f(\x^* + w \u) $ to converge in \eqref{ROCseries1}, we must have $w \leq R \leq \bigg(\limsup_{j\to \infty} \sqrt[j]{\frac{1}{j!}\bigg\lvert\frac{d^j }{ dw^j}[\nabla^2 f(\x^* + w \u)]_{i,l}\bigg\vert_{w=0}\bigg\rvert}\bigg)^{-1} $ for any unit vector $\u$. Precisely, the ROC for \eqref{ROCseries1} is given by {$$ R = \min_{i,l}\inf_{\norm{\u} = 1}\bigg(\limsup_{j\to \infty} \sqrt[j]{\frac{1}{j!}\bigg\lvert\frac{d^j }{ dw^j}[\nabla^2 f(\x^* + w \u)]_{i,l}\bigg\vert_{w=0}\bigg\rvert}\bigg)^{-1}$$} which is strictly positive. Next, due to $ \frac{1}{n}\max_{i,l}\{\abs{[\A]_{i,l}}\}  \leq\norm{\A}_2 $, we will have the following for any $(i,l)^{th}$ element of $\A = \frac{d^j }{ dw^j}\nabla^2 f(\x^* + w \u)\bigg\vert_{w=0} $:
	  \begin{align}
   & \frac{1}{j!} \bigg\lvert\frac{d^j }{ dw^j}[\nabla^2 f(\x^* + w \u)]_{i,l}\bigg\vert_{w=0} \bigg\rvert \leq \frac{n}{j!} \norm{ \frac{d^j }{dw^j}\nabla^2 f(\x^* + w \u)\bigg\vert_{w=0}}_2 \\
  & \implies \sqrt[j]{\frac{1}{j!} \bigg\lvert\frac{d^j }{ dw^j}[\nabla^2 f(\x^* + w \u)]_{i,l}\bigg\vert_{w=0} \bigg\rvert} \leq \sqrt[j]{\frac{n}{j!} \norm{ \frac{d^j }{dw^j}\nabla^2 f(\x^* + w \u)\bigg\vert_{w=0}}_2} \\
   & \implies \limsup_{j \to \infty}\sqrt[j]{\frac{1}{j!} \bigg\lvert\frac{d^j }{ dw^j}[\nabla^2 f(\x^* + w \u)]_{i,l}\bigg\vert_{w=0} \bigg\rvert}  \leq  \limsup_{j \to \infty}n^{1/j}\sqrt[j]{\frac{1}{j!} \norm{ \frac{d^j }{dw^j}\nabla^2 f(\x^* + w \u)\bigg\vert_{w=0}}_2} \\
  &   \implies \bigg(\limsup_{j \to \infty}\sqrt[j]{\frac{1}{j!} \bigg\lvert\frac{d^j }{ dw^j}[\nabla^2 f(\x^* + w \u)]_{i,l}\bigg\vert_{w=0} \bigg\rvert}\bigg)^{-1}  \geq  \bigg(\limsup_{j \to \infty}\sqrt[j]{\frac{1}{j!} \norm{ \frac{d^j }{dw^j}\nabla^2 f(\x^* + w \u)\bigg\vert_{w=0}}_2}\bigg)^{-1} \nonumber \\&= \bigg(\limsup_{j\to \infty} \sqrt[j]{\frac{r_j(\u)}{j!}}\bigg)^{-1} \\
    & \implies R = \min_{i,l} \inf_{\norm{\u}=1} \bigg(\limsup_{j \to \infty}\sqrt[j]{\frac{1}{j!} \bigg\lvert\frac{d^j }{ dw^j}[\nabla^2 f(\x^* + w \u)]_{i,l}\bigg\vert_{w=0} \bigg\rvert}\bigg)^{-1}  \geq \inf_{\norm{\u}=1} \bigg(\limsup_{j\to \infty} \sqrt[j]{\frac{r_j(\u)}{j!}}\bigg)^{-1},
   \label{ROCseries2}
	  \end{align}
	  where we used the $\limsup$ product rule in the second last step. Now \eqref{ROCseries2} implies that the quantity $\inf_{\norm{\u}=1}\bigg(\limsup_{j\to \infty} \sqrt[j]{\frac{r_j(\u)}{j!}}\bigg)^{-1} $ is upper bounded by the radius of convergence $R$ of the series in \eqref{ROCseries1}. Next, due to the inequality $ \frac{1}{n^2} \norm{\A}_2  \leq \frac{1}{n^2} \norm{\A}_F \leq  \frac{1}{n}\max_{i,l}\{\abs{[\A]_{i,l}}\}  $, for the maximum absolute element of $\A = \frac{d^j }{ dw^j}\nabla^2 f(\x^* + w \u)\bigg\vert_{w=0} $ denoted by $\bigg\lvert\frac{d^j }{ dw^j}[\nabla^2 f(\x^* + w \u)]_{m(j),q(j)}\bigg\vert_{w=0} \bigg\rvert$ we have the following:\footnote{Notice that the position $(m(j),q(j))$ of the maximum absolute element depends on $j$.}
	  \begin{align}
   & \frac{n}{j!} \bigg\lvert\frac{d^j }{ dw^j}[\nabla^2 f(\x^* + w \u)]_{m(j),q(j)}\bigg\vert_{w=0} \bigg\rvert \geq \frac{1}{j!} \norm{ \frac{d^j }{dw^j}\nabla^2 f(\x^* + w \u)\bigg\vert_{w=0}}_2 \\
  &\implies \sqrt[j]{\frac{n}{j!} \bigg\lvert\frac{d^j }{ dw^j}[\nabla^2 f(\x^* + w \u)]_{m(j),q(j)}\bigg\vert_{w=0} \bigg\rvert} \geq \sqrt[j]{\frac{1}{j!} \norm{ \frac{d^j }{dw^j}\nabla^2 f(\x^* + w \u)\bigg\vert_{w=0}}_2} \\
   & \implies \limsup_{j \to \infty}n^{1/j}\sqrt[j]{\frac{1}{j!} \bigg\lvert\frac{d^j }{ dw^j}[\nabla^2 f(\x^* + w \u)]_{m(j),q(j)}\bigg\vert_{w=0} \bigg\rvert}  \geq  \limsup_{j \to \infty}\sqrt[j]{\frac{1}{j!} \norm{ \frac{d^j }{dw^j}\nabla^2 f(\x^* + w \u)\bigg\vert_{w=0}}_2} \\
    & \implies \bigg(\limsup_{j \to \infty}\sqrt[j]{\frac{1}{j!} \bigg\lvert\frac{d^j }{ dw^j}[\nabla^2 f(\x^* + w \u)]_{m(j),q(j)}\bigg\vert_{w=0} \bigg\rvert}\bigg)^{-1}  \leq  \bigg(\limsup_{j \to \infty}\sqrt[j]{\frac{1}{j!} \norm{ \frac{d^j }{dw^j}\nabla^2 f(\x^* + w \u)\bigg\vert_{w=0}}_2}\bigg)^{-1} \nonumber \\ & \hspace{8.5cm}= \bigg(\limsup_{j\to \infty} \sqrt[j]{\frac{r_j(\u)}{j!}}\bigg)^{-1}\\
   & \implies R \leq \inf_{\norm{\u}=1} \bigg(\limsup_{j \to \infty}\sqrt[j]{\frac{1}{j!} \bigg\lvert\frac{d^j }{ dw^j}[\nabla^2 f(\x^* + w \u)]_{m(j),q(j)}\bigg\vert_{w=0} \bigg\rvert}\bigg)^{-1}  \leq  \inf_{\norm{\u}=1} \bigg(\limsup_{j\to \infty} \sqrt[j]{\frac{r_j(\u)}{j!}}\bigg)^{-1},
   \label{ROCseries3}
	  \end{align}
	  where the L.H.S. of the last inequality holds by $\min_{i,l} \inf_{\norm{\u}=1} \bigg(\limsup_{j \to \infty}\sqrt[j]{\frac{1}{j!} \bigg\lvert\frac{d^j }{ dw^j}[\nabla^2 f(\x^* + w \u)]_{i,l}\bigg\vert_{w=0} \bigg\rvert}\bigg)^{-1} \leq \inf_{\norm{\u}=1} \bigg(\limsup_{j \to \infty}\sqrt[j]{\frac{1}{j!} \bigg\lvert\frac{d^j }{ dw^j}[\nabla^2 f(\x^* + w \u)]_{m(j),q(j)}\bigg\vert_{w=0} \bigg\rvert}\bigg)^{-1} $.
	Finally, combining \eqref{ROCseries2} and \eqref{ROCseries3} we get:
	\begin{align}
	   R \leq \inf_{\norm{\u}= 1} \bigg(\limsup_{j\to \infty} \sqrt[j]{\frac{r_j(\u)}{j!}}\bigg)^{-1} &\leq R \\
	   \implies \inf_{\norm{\u}= 1} \bigg(\limsup_{j\to \infty} \sqrt[j]{\frac{r_j(\u)}{j!}}\bigg)^{-1} &= R.
	\end{align}
	 \endproof

\section{Radial vector $\u_{k}$ in terms of the initialization $\u_{0}$} \label{Appendix B}
	 \subsection{\textbf{Proof of Lemma \ref{lem5}}}
	 \proof{}
	  Combining the equation $\u_{k} = \x_{k} - \x^*$ this with gradient descent update yields
	\begin{align}
	\u_{k+1} - \u_{k} = -\alpha \nabla f(\x_{k}).
	\end{align}
	Next, substituting \eqref{gradapprox} here, we get the following recursion:
	\begin{align}
	\u_{k+1} - \u_{k} & =  - \alpha \bigg(\nabla^2 f(\x^*) +  \frac{\norm{\u_{k}}}{2} \H(\hat{\u}_k) + \mathcal{O}(\epsilon^2) \bigg) \u_{k}  \\
	\u_{k+1} & =  \bigg(\mathbf{I}- \alpha \bigg(\nabla^2 f(\x^*) +  \frac{\norm{\u_{k}}}{2} \H(\hat{\u}_k)+\mathcal{O}(\epsilon^2) \bigg)\bigg) \u_{k}.
	\end{align}
	Finally substituting \eqref{errorbound} here and applying the spectral theorem to the matrices $\mathbf{I}$ and $\nabla^2 f(\x^*)$ yields
	\begin{align}
	\u_{k+1} & =  \bigg(\sum_{i=1}^{n} \v_i(0)\v_i(0)^{T}- \alpha \bigg(\sum_{i=1}^{n} \lambda_i(0)\v_i(0)\v_i(0)^{T} +  \frac{\norm{\u_{k}}}{2}\bigg( \sum_{i=1}^{n} \bigg(\langle \v_i(0), \H(\hat{\u}_k)\v_i(0)\rangle \v_i(0)\v_i(0)^{T} \nonumber \\ & + \lambda_i(0)(\sum_{l \neq i} \frac{\langle \v_l(0), \H(\hat{\u}_k)\v_i(0)\rangle}{\lambda_{i}(0)-\lambda_{l}(0)}\v_l(0))\v_i(0)^{T} + \lambda_i(0)\v_i(0)(\sum_{l \neq i} \frac{\langle \v_l(0), \H(\hat{\u}_k)\v_i(0)\rangle}{\lambda_{i}(0)-\lambda_{l}(0)}\v_l(0))^{T}\bigg) \bigg) \bigg)+\mathcal{O}(\epsilon^2)\bigg) \u_{k}  \\
	\u_{k+1} & =  \bigg[\sum_{i=1}^{n} \bigg(1-\alpha \lambda_{i}(0)-\alpha \frac{\norm{\u_{k}}}{2}\langle \v_i(0), \H(\hat{\u}_k)\v_i(0)\rangle\bigg)\v_i(0)\v_i(0)^{T}\nonumber \\ & - \alpha  \frac{\norm{\u_{k}}}{2}\sum_{i=1}^{n}  \sum_{l \neq i}  \frac{\langle \v_l(0), \H(\hat{\u}_k)\v_i(0)\rangle\lambda_i(0)}{\lambda_{i}(0)-\lambda_{l}(0)}\bigg(\v_l(0)\v_i(0)^{T}   +\v_i(0)\v_l(0)^{T}\bigg)  \bigg]\u_{k} +\mathcal{O}(\epsilon^2)\u_{k}.  \label{interbound}
	\end{align}
	Next, we start analyzing the coefficients of spectral components $\v_i(0)\v_l(0)^{T}$ for any $(i,l)$ pair.
	
	\subsubsection{Coefficient bounds:}
	We start with \eqref{interbound} and analyze it in terms of the stable subspace $\mathcal{E}_{S}$ and unstable subspace $\mathcal{E}_{US}$ of the Hessian $\nabla^2 f(\x^*)$. To this end we rewrite \eqref{interbound} and split its first term into the stable and unstable spectral components. The stable spectral components result from the positive eigenvalues of $\nabla^2 f(\x^*)$ whereas the unstable spectral components result from its negative eigenvalues.
	
	\begin{align}
	\hspace{-2cm} \u_{k+1} &  =  \bigg[\sum_{i=1}^{n} \bigg(1-\alpha \lambda_{i}(0)-\alpha \frac{\norm{\u_{k}}}{2}\langle \v_i(0), \H(\hat{\u}_k)\v_i(0)\rangle\bigg)\v_i(0)\v_i(0)^{T}\nonumber \\ & - \alpha  \frac{\norm{\u_{k}}}{2}\sum_{i=1}^{n}  \sum_{l \neq i}  \frac{\langle \v_l(0), \H(\hat{\u}_k)\v_i(0)\rangle\lambda_i(0)}{\lambda_{i}(0)-\lambda_{l}(0)}\bigg(\v_l(0)\v_i(0)^{T}   +\v_i(0)\v_l(0)^{T}\bigg)  \bigg]\u_{k} +\mathcal{O}(\epsilon^2)\u_{k} \\
	 & =  \bigg[\sum_{i \in \mathcal{N}_{S}} \bigg(1-\alpha \lambda_{i}(0)-\alpha \frac{\norm{\u_{k}}}{2}\langle \v_i(0), \H(\hat{\u}_k)\v_i(0)\rangle\bigg)\v_i(0)\v_i(0)^{T} \nonumber \\ & + \sum_{j \in \mathcal{N}_{US}} \bigg(1-\alpha \lambda_{j}(0)-\alpha \frac{\norm{\u_{k}}}{2}\langle \v_j(0), \H(\hat{\u}_k)\v_j(0)\rangle\bigg)\v_j(0)\v_j(0)^{T}\nonumber \\ & - \alpha  \frac{\norm{\u_{k}}}{2}\sum_{i=1}^{n}  \sum_{l \neq i}  \frac{\langle \v_l(0), \H(\hat{\u}_k)\v_i(0)\rangle\lambda_i(0)}{\lambda_{i}(0)-\lambda_{l}(0)}\bigg(\v_l(0)\v_i(0)^{T}   +\v_i(0)\v_l(0)^{T}\bigg)  \bigg]\u_{k} +\mathcal{O}(\epsilon^2)\u_{k} \\
	 & =  \bigg[\sum_{i \in \mathcal{N}_{S}}  c^{s}_i(k)\v_i(0)\v_i(0)^{T} + \sum_{j \in \mathcal{N}_{US}}  c^{us}_j(k)\v_j(0)\v_j(0)^{T}\nonumber \\ & + \sum_{i=1}^{n}  \sum_{l \neq i} \bigg(d_{l,i}(k)\v_l(0)\v_i(0)^{T}   +d_{i,l}(k)\v_i(0)\v_l(0)^{T}\bigg)  \bigg]\u_{k} +\mathcal{O}(\epsilon^2)\u_{k}, \label{recursion}
	\end{align}
where the coefficient terms $c^{s}_i(k)$, $c^{us}_j(k)$ and $d_{l,i}(k)$ in \eqref{recursion} are as follows:
	\begin{align}
	c^{s}_i(k) &= \bigg(1-\alpha \lambda_{i}(0)-\alpha \frac{\norm{\u_{k}}}{2}\langle \v_i(0), \H(\hat{\u}_k)\v_i(0)\rangle\bigg)   \label{coeffbounds1}\\
	c^{us}_j(k) &= \bigg(1-\alpha \lambda_{j}(0)-\alpha \frac{\norm{\u_{k}}}{2}\langle \v_j(0), \H(\hat{\u}_k)\v_j(0)\rangle\bigg)  \label{coeffbounds2} \\
	d_{i,l}(k) &= d_{l,i}(k) =  \frac{\langle \v_l(0), \H(\hat{\u}_k)\v_i(0)\rangle\lambda_i(0)\alpha \norm{\u_{k}}}{2(\lambda_{l}(0)-\lambda_{i}(0))}. \label{coeffbounds3}
	\end{align}	
	Now, from \eqref{hessianapprox} and the Lipschitz continuity of the Hessian (Assumption \textbf{A3}), we get the following bound:
	\begin{align}
	\nabla^2 f(\x^* + p{\u_k}) &= \nabla^2 f(\x^*) + p \norm{\u_{k}} \H(\hat{\u}_k) + \mathcal{O}(\epsilon^2).
	\end{align}
	Recall that the term $\mathcal{O}(\epsilon^2)$ comes from \eqref{taylorremfirstorder}. Therefore, to further simplify the above equation, we replace $\mathcal{O}(\epsilon^2)$ with $\int_{0}^{w} u \frac{d^2 \G}{du^2} du$ from \eqref{taylorremfirstorder} where $w =p \norm{\u_{k}}$. Then taking the norm of both sides, followed by triangle inequality and using Assumption \textbf{A3} yields
	\begin{align}
		\nabla^2 f(\x^* + p{\u_k}) &= \nabla^2 f(\x^*) + p \norm{\u_{k}} \H(\hat{\u}_k) + \int_{0}^{w} u \frac{d^2 \G}{du^2} du  \\
	\norm{\H(\hat{\u}_k)}_2 &=   \frac{1}{p \norm{\u_{k}}}\norm{\nabla^2 f(\x^* + p{\u}_k) - \nabla^2 f(\x^*)- \int_{0}^{w} u \frac{d^2 \G}{du^2}du}_2\\
	 &\leq  \frac{M}{p \norm{\u_{k}}}\norm{
		\x^* + p{\u}_k - \x^*} + \frac{\norm{ \int_{0}^{w} u \frac{d^2 \G}{du^2}du}_2}{p \norm{\u_{k}}} \\
 &\leq   M + \frac{ \bigg(\int_{0}^{w}\norm{  \frac{d^2 \G}{du^2}}_2^2 du \bigg)^{\frac{1}{2}} \bigg(\int_{0}^{w}u^2 du\bigg)^{\frac{1}{2}}}{w} \leq M + \frac{B_2w}{\sqrt{3}} \leq M + \mathcal{O}(\epsilon). \label{reviewreporteqnew}
	\end{align}
	Note that in the last step, we used the Cauchy Schwarz inequality followed by the same bound $\norm{ \frac{d^2 \G}{du^2}}_2 \leq B_2$ as in the steps following \eqref{taylorremfirstorder}. For the case when $p \norm{\u_{k}} \to 0$, the bound on $ \norm{\H(\hat{\u}_k)}_2$ can be evaluated by using the substitution $w =p\norm{\u_{k}}$:
	\begin{align}
	\norm{\H(\hat{\u}_k)}_2  &\leq \lim_{p\norm{\u_{k}} \to 0} \frac{M}{p \norm{\u_{k}}}\norm{
		\x^* + p{\u}_k - \x^*} + \lim_{p\norm{\u_{k}} \to 0} \frac{ \bigg(\int_{0}^{w}\norm{  \frac{d^2 \G}{du^2}}_2^2 du \bigg)^{\frac{1}{2}} \bigg(\int_{0}^{w}u^2 du\bigg)^{\frac{1}{2}}}{p \norm{\u_{k}}} \\
	  &\leq \lim_{w \to 0} \frac{M}{w}w + \lim_{w\to 0} \frac{ \bigg(\int_{0}^{w}\norm{  \frac{d^2 \G}{du^2}}_2^2 du \bigg)^{\frac{1}{2}} \bigg(\int_{0}^{w}u^2 du\bigg)^{\frac{1}{2}}}{w} \\
	 &\leq M + \lim_{w\to 0} \frac{\bigg(\int_{0}^{w}\norm{  \frac{d^2 \G}{du^2}}_2^2 du \bigg)^{\frac{1}{2}}w^{1/2} }{\sqrt{3}} = M + \lim_{w\to 0}\bigg(\int_{0}^{w}\norm{  \frac{d^2 \G}{du^2}}_2^2 du \bigg)^{\frac{1}{2}}\lim_{w\to 0} \frac{w^{1/2} }{\sqrt{3}} = M.
	\end{align}
	Note that in the last step, we used $\lim_{w\to 0}\bigg(\int_{0}^{w}\norm{  \frac{d^2 \G}{du^2}}_2^2 du \bigg)^{\frac{1}{2}} = \bigg(\int_{0}^{1}\lim_{w\to 0}\mathds{1}_{[0,w]}\norm{  \frac{d^2 \G}{du^2}}_2^2 du \bigg)^{\frac{1}{2}}=0$ by the dominated convergence theorem where $\mathds{1}_{[0,w]}$ is the indicator function on $[0,w]$.
	
		Hence for any eigenvectors $\v_i(0), \v_j(0)$ we have that
	\begin{align}
	-M- \mathcal{O}(\epsilon)& \leq \langle \v_i(0), \H(\hat{\u}_k)\v_j(0)\rangle \leq M + \mathcal{O}(\epsilon).
	\end{align}
		Using Assumptions \textbf{A2} and \textbf{A4}, for the stable subspace $\mathcal{E}_{S}$, we have the following bound on $\lambda_{i}(0)$:
	\begin{align}
	\beta \leq \lambda_{i}(0) \leq L.
	\end{align}
	Similarly for the unstable subspace $\mathcal{E}_{US}$, we have the following bound on $\lambda_{j}(0)$ from Assumptions \textbf{A2} and \revise{\textbf{A4}}:
	\begin{align}
	-L \leq \lambda_{j}(0) \leq -\beta.
	\end{align}
	Now substituting these bounds into \eqref{coeffbounds1}, \eqref{coeffbounds2}, \eqref{coeffbounds3} and using the fact that $\norm{\u_{k}} < \epsilon$, we get the following bounds on the coefficients:
	\begin{align}
	\bigg(1 - \alpha L -\frac{\alpha \epsilon M}{2} -\mathcal{O}(\epsilon^2)  \bigg) \leq c^{s}_i(k) & \leq \bigg(1 - \alpha \beta +\frac{\alpha \epsilon M}{2} + \mathcal{O}(\epsilon^2) \bigg)  \label{bound1} \\
	\bigg(1 +\alpha \beta -\frac{\alpha \epsilon M}{2}  - \mathcal{O}(\epsilon^2) \bigg) \leq c^{us}_j(k) & \leq \bigg(1 +\alpha L +\frac{\alpha \epsilon M}{2}+\mathcal{O}(\epsilon^2)   \bigg) \label{bound2} \\
	-\frac{\alpha \epsilon ML}{2 \delta} - \mathcal{O}(\epsilon^2)   \leq d_{i,l}(k) & \leq \frac{\alpha \epsilon ML}{2 \delta}+ \mathcal{O}(\epsilon^2).  \label{bound3}
	\end{align}
	
	After establishing the bounds on the coefficients $c^{s}_i(k), c^{us}_j(k), d_{i,l}(k)  $, we further analyze the recursive vector update equation \eqref{recursion} and induct it from $k=0$ to $k= K-1$ so as to obtain $\u_K$ in terms of $\u_0$:
	\begin{align}
	\u_{k+1} & =  \bigg[\sum_{i \in \mathcal{N}_{S}}  c^{s}_i(k)\v_i(0)\v_i(0)^{T} + \sum_{j \in \mathcal{N}_{US}}  c^{us}_j(k)\v_j(0)\v_j(0)^{T}\nonumber \\ & + \sum_{i=1}^{n}  \sum_{l \neq i} \bigg(d_{l,i}(k)\v_l(0)\v_i(0)^{T}   +d_{i,l}(k)\v_i(0)\v_l(0)^{T}\bigg)  \bigg]\u_{k} +\mathcal{O}(\epsilon^2)\u_{k} \\
	\implies \u_{K} & = \prod_{k=0}^{K-1} \bigg[\mathcal{O}(\epsilon^2) +\sum_{i \in \mathcal{N}_{S}}  c^{s}_i(k)\v_i(0)\v_i(0)^{T} + \sum_{j \in \mathcal{N}_{US}}  c^{us}_j(k)\v_j(0)\v_j(0)^{T}\nonumber \\ & + \sum_{i=1}^{n}  \sum_{l \neq i} \bigg(d_{l,i}(k)\v_l(0)\v_i(0)^{T}   +d_{i,l}(k)\v_i(0)\v_l(0)^{T}\bigg)  \bigg]\u_{0}. \label{matrixprod}
	\end{align}	
	Observe that in the above expression, the vector $\u_{K}$ results from a product of $K$ matrices. Each of these matrices comes from a linear combination of $n^2$ different matrices given by $\v_i(0)\v_i(0)^{T}$ for $i \in \mathcal{N}_{S} $, $\v_j(0)\v_j(0)^{T}$ for $j \in \mathcal{N}_{US}$, the cross terms $\v_l(0)\v_i(0)^{T},\v_i(0)\v_l(0)^{T}$ with $i \neq l$ and in addition to this a matrix term of order $\mathcal{O}(\epsilon^2)$.
	
	Next, using the orthogonality of eigenvectors we obtain $\v_i(0)^T\v_j(0) = 0$ for $i \neq j$ and $\v_i(0)^T\v_j(0) = 1$ for $i = j$. Therefore by induction it can be readily inferred that the $K$ matrix product is a linear combination of the same $n^2$ matrices plus all the matrix error terms of the order $\mathcal{O}(\epsilon^2)$ and above.
	Hence we rewrite \eqref{matrixprod} as follows:
	\begin{align}
	\u_{K} = \prod_{k=0}^{K-1} \bigg[\A_k + \B_k +\mathcal{O}(\epsilon^2)   \bigg]\u_{0}, \label{substitution}
	\end{align}
	where $\A_k = \sum_{i \in \mathcal{N}_{S}}  c^{s}_i(k)\v_i(0)\v_i(0)^{T} + \sum_{j \in \mathcal{N}_{US}}  c^{us}_j(k)\v_j(0)\v_j(0)^{T}$ and $\B_k = \sum_{i=1}^{n}  \sum_{l \neq i} \bigg(d_{l,i}(k)\v_l(0)\v_i(0)^{T}   +d_{i,l}(k)\v_i(0)\v_l(0)^{T}\bigg)$. From \eqref{bound3} the term $\B_k$ is of order $\mathcal{O}(\epsilon)$. Therefore, this equation can be written more compactly as
	\begin{align}
	\u_{K} = \prod_{k=0}^{K-1} \bigg[\A_k + \epsilon \P_k   \bigg]\u_{0},
	\end{align} where $\epsilon \P_k = \B_k+\mathcal{O}(\epsilon^2)$.
	
	Next we analyze the matrix product $ \prod_{k=0}^{K-1} \bigg[\A_k + \epsilon \P_k   \bigg]$. Taking the norm of this product, followed by the supremum over $k$ and using the triangle inequality yields
	\begin{align}
	\norm{\prod_{k=0}^{K-1} \bigg[\A_k + \epsilon \P_k   \bigg]}_2 & \leq  \prod_{k=0}^{K-1} \norm{\bigg[\A_k + \epsilon \P_k   \bigg]}_2 \\
 & \leq  \prod_{k=0}^{K-1} \sup_{0\leq k \leq K-1}\norm{\bigg[\A_k + \epsilon \P_k   \bigg]}_2 \\
	 & \leq  \prod_{k=0}^{K-1}\bigg[\sup_{0\leq k \leq K-1} \norm{\A_k}_2 + \epsilon \sup_{0\leq k \leq K-1}\norm{\P_k}_2 \bigg] \\
		 & \leq  \prod_{k=0}^{K-1}\bigg[\norm{\A}_2 + \epsilon \norm{\P}_2 \bigg] = \bigg(\norm{\A}_2 + \epsilon \norm{\P}_2\bigg)^K, \label{binomialbound}
	\end{align}
	where in the last step we have used the substitutions $ \sup_{0\leq k \leq K-1}\norm{\A_k}_2 = \norm{\A}_2$ and $\sup_{0\leq k \leq K-1}\norm{\P_k}_2 = \norm{\P}_2$ for some arbitrary matrices $\A$ and $\P$.
	
	Now observe that the product term on the right-hand side of \eqref{binomialbound} has a binomial expansion which can be written compactly as
	\begin{align}
	\norm{\prod_{k=0}^{K-1} \bigg[\A_k + \epsilon \P_k   \bigg]}_2 & \leq   \sum_{r=0}^{K} \binom{K}{r}( \epsilon \norm{\P}_2)^r \norm{\A}_2^{K-r}
	= \norm{\A}_2^K\bigg( 1 + \epsilon \frac{\norm{\P}_2}{\norm{\A}_2}\bigg)^K. \label{orderbound}
	\end{align}
	Next, consider the term $\bigg( 1 + \epsilon \frac{\norm{\P}_2}{\norm{\A}_2}\bigg)^K$ on the right-hand side of above bound. For the function $g_{\omega}(x) = (1+x)^{\omega}$ such that $\omega \in \mathbb{R}$, its expansion and the remainder term are given by
	\begin{align}
(1+x)^{\omega} &= \sum_{k=0}^{\infty} \binom{\omega}{k}x^k \\
R_j(x) & = \int_{0}^{x} \frac{(x-z)^j}{j!} (j+1)! \binom{\omega}{j+1} (1+z)^{\omega -j-1}dz,
	\end{align} 	
	where we have that $\limsup_{j\to \infty}R_j(x) = 0$ for $\abs{x} < 1$.
	
	Therefore using this remainder expression for the term $\bigg( 1 + \epsilon \frac{\norm{\P}_2}{\norm{\A}_2}\bigg)^K$ with $x=\epsilon  \frac{\norm{\P}_2}{\norm{\A}_2}$ we will have
	\begin{align}
	R_1(x) & = \int_{0}^{x} \frac{(x-z)}{1!} 2! \binom{K}{2} (1+z)^{K -2}dz \\
	& = K(K-1) \bigg( \frac{(1+x)^K}{K(K-1)} - \frac{1+x}{K-1} + \frac{1}{K} \bigg). \label{binomrem}
	\end{align}
	For $\abs{x}<1$ and $\abs{Kx}\ll1$, we can use the approximation $(1+x)^K \approx 1+ Kx + \binom{K}{2}x^2$. Substituting this approximation in \eqref{binomrem}, we get $R_1(x)$ as
	\begin{align}
	R_1(x) & \approx K(K-1) \bigg( \frac{1+Kx + \binom{K}{2}x^2}{K(K-1)} - \frac{1+x}{K-1} + \frac{1}{K} \bigg) = \frac{K(K-1)}{2}x^2 \\
	R_1\bigg(\epsilon  \frac{\norm{\P}_2}{\norm{\A}_2}\bigg) & \approx  \frac{K(K-1)}{2}\bigg(\epsilon  \frac{\norm{\P}_2}{\norm{\A}_2}\bigg)^2 = \mathcal{O}\bigg((K \epsilon)^2\bigg).
	\end{align}
    Hence for $\epsilon\frac{\norm{\P}_2}{\norm{\A}_2} < 1$ and $K \epsilon \ll1$, we can substitute this bound in \eqref{orderbound} as follows:
    \begin{align}
    \norm{\prod_{k=0}^{K-1} \bigg[\A_k + \epsilon \P_k   \bigg]}_2 & \leq  \norm{\A}_2^K\bigg( 1 + \epsilon \frac{\norm{\P}_2}{\norm{\A}_2}\bigg)^K \\ &= \norm{\A}_2^K \bigg( 1+ K\epsilon \frac{\norm{\P}_2}{\norm{\A}_2} + R_1\bigg(\epsilon  \frac{\norm{\P}_2}{\norm{\A}_2}\bigg) \bigg) \\
    & \approx  \norm{\A}_2^K \bigg( 1+ K\epsilon \frac{\norm{\P}_2}{\norm{\A}_2} + \mathcal{O}\bigg((K\epsilon)^2\bigg) \bigg).  \label{approximation1}
    \end{align}
    This approximate upper bound implies that the upper bound on the norm of matrix product $\prod_{k=0}^{K-1} \bigg[\A_k + \epsilon \P_k   \bigg]$ can be approximately expanded up to an $ \epsilon $ precision term accompanied with a remainder term of $\mathcal{O}\bigg(\norm{\A}_2^K (K\epsilon)^2\bigg)$ as long as $K \epsilon \ll1$.

    Next we obtain a lower bound on the inverse of the norm of matrix product $\prod_{k=0}^{K-1} \bigg[\A_k + \epsilon \P_k   \bigg]^{-1}$.
    Taking the inverse of the norm of this product, using the identities $\norm{\Z}_2 \geq \norm{\Z^{-1}}_2^{-1} $, $ \norm{(\mathbf{I}+\Z)^{-1}}_2^{-1} \geq (1- \norm{\Z}_2)$, followed by taking the infimum over $k$ yields
    \begin{align}
        \norm{\prod_{k=0}^{K-1} \bigg[\A_k + \epsilon \P_k   \bigg]^{-1}}_2^{-1}
 & \geq  \prod_{k=0}^{K-1} \norm{\A_k^{-1}}_2^{-1} \bigg(1-\epsilon \norm{ \A_k^{-1}  \P_k }_2  \bigg)\\
 & \geq  \prod_{k=0}^{K-1} \inf_{0\leq k \leq K-1}\norm{\A_k^{-1}}_2^{-1} \bigg(1-\epsilon \sup_{0\leq k \leq K-1}\norm{ \A_k^{-1}  \P_k }_2  \bigg) \\
    & \geq  \bigg(\inf_{0\leq k \leq K-1}\norm{\A_k^{-1}}_2^{-1}\bigg)^K \bigg(1-\epsilon \sup_{0\leq k \leq K-1}\norm{ \A_k^{-1}  \P_k }_2  \bigg)^K .\label{lowerboundmat}
    \end{align}
    Now repeating the previous analysis of the upper bound here will give the conclusion that the lower bound on inverse of the norm of matrix product $\prod_{k=0}^{K-1} \bigg[\A_k + \epsilon \P_k   \bigg]^{-1}$ can be approximately computed up to $K\epsilon$ precision provided $\epsilon \sup_{0\leq k \leq K-1}\norm{ \A_k^{-1}  \P_k }_2 < 1$ and $K \epsilon \ll 1$ if the step size $\alpha< \frac{1}{L}$. The reasoning for having $\alpha< \frac{1}{L}$ will be discussed in the subsequent section when we derive some feasible range for $\epsilon$ as well as the case where $\alpha \approx \frac{1}{L}$. In particular, the inequality \eqref{lowerboundmat} can be simplified further as
      \begin{align}
      \norm{\prod_{k=0}^{K-1} \bigg[\A_k + \epsilon \P_k   \bigg]^{-1}}_2^{-1} \geq \norm{\A^{-1}}_2^{-K} \bigg( 1- K\epsilon \frac{\norm{\P}_2}{\norm{\A^{-1}}^{-1}_2} - \mathcal{O}\bigg((K\epsilon)^2\bigg) \bigg)	 \label{approximation1a}
      \end{align}
      for $K\epsilon \ll 1$ and $\epsilon \frac{\norm{\P}_2}{\norm{\A^{-1}}^{-1}_2} <1 $ where we have that $ \sup_{0\leq k \leq K-1}\norm{\A_k^{-1}}_2 = \norm{\A^{-1}}_2$ and $\sup_{0\leq k \leq K-1}\norm{\P_k}_2 = \norm{\P}_2$ for the matrices $\A$ and $\P$ used previously.

 Now, if $\nu_n\leq \dots \leq \nu_1$ are the absolute value of the eigenvalues of the matrix product $ \prod_{k=0}^{K-1} \bigg[\A_k + \epsilon \P_k   \bigg]$, then using \eqref{approximation1} and \eqref{approximation1a}, we have the condition
 \begin{align}
 \norm{\A^{-1}}_2^{-K} \bigg( 1- K\epsilon \frac{\norm{\P}_2}{\norm{\A^{-1}}^{-1}_2} - \mathcal{O}\bigg((K\epsilon)^2\bigg) \bigg) \leq \nu_n \leq \dots \leq \nu_1  \leq  \norm{\A}_2^K \bigg( 1+ K\epsilon \frac{\norm{\P}_2}{\norm{\A}_2} + \mathcal{O}\bigg((K\epsilon)^2\bigg) \bigg).
\end{align}

    Therefore we can conclude that the matrix product \eqref{matrixprod} can be approximately computed up to $K\epsilon$ precision provided $K\epsilon \ll 1$, $\epsilon \frac{\norm{\P}_2}{\norm{\A}_2} <1 $ and $\epsilon \frac{\norm{\P}_2}{\norm{\A^{-1}}^{-1}_2} <1 $. At this point, we are interested in analyzing the matrix product in \eqref{matrixprod} only for iterations $K = \mathcal{O}(\frac{1}{\epsilon})$. This is done so as to derive exit times and initial conditions for trajectories that can escape a strict saddle point in linear time. It is also remarked that we could have retained the higher-order terms $\mathcal{O}\bigg(\norm{\A}_2^K(K\epsilon)^r\bigg)$ in the above matrix product \eqref{lowerboundmat} if we wanted to analyze trajectories with polynomial or even exponential rates of escape.
    	 \endproof
  \subsection{\textbf{Proof of Lemma \ref{lem6}}}  	
\proof{}
	For values of $K = \mathcal{O}(\frac{1}{\epsilon})$ we explicitly compute the matrix product in \eqref{matrixprod} up to $K\epsilon$ precision and drop all the higher order terms ($\epsilon^2$ and above) that collectively act as a single remainder term of an approximate order $\mathcal{O}\bigg(\norm{\A}_2^K(K\epsilon)^2\bigg)$. From \eqref{bound3} we know that only the coefficients $d_{i,l}(k)$ are of order $\mathcal{O}(\epsilon)$, hence we now expand \eqref{matrixprod} only up to first order in $d_{i,l}(k)$ to obtain the following approximation:
	\begin{align}
	\u_{K} & \approx \tilde{\u}_K=\bigg[\sum_{i \in \mathcal{N}_{S}}  \bigg(\prod_{k=0}^{K-1}c^{s}_i(k)\bigg)\v_i(0)\v_i(0)^{T} + \sum_{j \in \mathcal{N}_{US}}  \prod_{k=0}^{K-1}\bigg(c^{us}_j(k)\bigg)\v_j(0)\v_j(0)^{T}\nonumber \\ & + \sum_{i \in \mathcal{N}_{S} }  \sum_{l \in \mathcal{N}_{S}} \sum_{r=0}^{K-1} \bigg(\prod_{k=0}^{r-1} c^{s}_i(k)\bigg) d_{i,l}(r) \bigg( \prod_{k=r+1}^{K-1} c^{s}_l(k)\bigg) \bigg(\v_l(0)\v_i(0)^{T}   +\v_i(0)\v_l(0)^{T} \bigg)
	\nonumber \\ & + \sum_{i \in \mathcal{N}_{S} }  \sum_{l \in \mathcal{N}_{US}} \sum_{r=0}^{K-1} \bigg(\prod_{k=0}^{r-1} c^{s}_i(k)\bigg) d_{i,l}(r) \bigg( \prod_{k=r+1}^{K-1} c^{us}_l(k)\bigg) \bigg(\v_l(0)\v_i(0)^{T}   +\v_i(0)\v_l(0)^{T} \bigg)
	\nonumber \\ & + \sum_{i \in \mathcal{N}_{US} }  \sum_{l \in \mathcal{N}_{S}} \sum_{r=0}^{K-1} \bigg(\prod_{k=0}^{r-1} c^{us}_i(k)\bigg) d_{i,l}(r) \bigg( \prod_{k=r+1}^{K-1} c^{s}_l(k)\bigg) \bigg(\v_l(0)\v_i(0)^{T}   +\v_i(0)\v_l(0)^{T} \bigg)
	\nonumber \\ & + \sum_{i \in \mathcal{N}_{US} }  \sum_{l \in \mathcal{N}_{US}} \sum_{r=0}^{K-1} \bigg(\prod_{k=0}^{r-1} c^{us}_i(k)\bigg) d_{i,l}(r) \bigg( \prod_{k=r+1}^{K-1} c^{us}_l(k)\bigg) \bigg(\v_l(0)\v_i(0)^{T}   +\v_i(0)\v_l(0)^{T} \bigg)
	\bigg]\u_{0}, \label{u_approx}
	\end{align}
		where we have that $\tilde{\u}_K$ as the $\epsilon$ approximate trajectory.	
		
	Next we express $\u_{0}$ as the sum of projections onto the stable subspace and unstable subspace of $\nabla^2 f(\x^{*})$ as follows:
	\begin{align}
	\u_{0} &= \epsilon\sum_{i \in \mathcal{N}_{S} } {\theta}^{s}_{i} \v_i(0) + \epsilon\sum_{j \in \mathcal{N}_{US} } {\theta}^{us}_{j} \v_j(0) \label{spectrum}\\
	\sum_{i \in \mathcal{N}_{S} }& ({\theta}^{s}_{i})^2  + \sum_{j \in \mathcal{N}_{US} } ({\theta}^{us}_{j})^2 = 1,
	\end{align}
	where $\epsilon\theta^s_i = \langle \u_{0}, \v_i(0)\rangle$, $\epsilon\theta^{us}_j = \langle \u_{0}, \v_j(0)\rangle$ with $\v_i(0) \in \mathcal{E}_S$ and $\v_j(0) \in \mathcal{E}_{US}$ respectively. Observe that \eqref{spectrum} has an $\epsilon$ multiplier because $\norm{\u_0} = \epsilon$. This is due to the fact that $\u_{0} + \x^* = \x_{0}$ and $\x_{0} \in \mathcal{\bar{B}}_{\epsilon}(\x^*) \backslash \mathcal{B}_{\epsilon}(\x^*)$.
	
 Now for all $i$ and $j$, the Hessian $\nabla^2 f(\x^*)$ can have eigenvectors $\v_i(0)$ and $\v_j(0)$ as well as $-\v_i(0)$ and $-\v_j(0)$. Therefore for the sake of analysis, the signs with these eigenvectors are chosen such that the respective coefficients $\theta^s_i$ and ${\theta}^{us}_{j}$ are positive for all $i$ and $j$. It is easy to show that such a choice always exists for all $i$ and $j$ because if $\langle \u_{0}, \v_i(0)\rangle > 0$ then $\langle \u_{0}, -\v_i(0)\rangle < 0$ and vice versa for any $i$ (and analogously for the index $j$).

	Finally substituting $\u_{0}$ in \eqref{u_approx}, we get the following result for $\u_{K}$:
	\begin{align}
	\u_{K}\approx \tilde{\u}_K =\epsilon\sum_{i \in \mathcal{N}_{S}}  \bigg(&\prod_{k=0}^{K-1}c^{s}_i(k){\theta}^{s}_i + \sum_{l \in \mathcal{N}_{S}} \sum_{r=0}^{K-1} \prod_{k=0}^{r-1} c^{s}_i(k) d_{i,l}(r)  \prod_{k=r+1}^{K-1} c^{s}_l(k){\theta}^{s}_l \nonumber \\ & + \sum_{l \in \mathcal{N}_{US}} \sum_{r=0}^{K-1} \prod_{k=0}^{r-1} c^{s}_i(k) d_{i,l}(r)  \prod_{k=r+1}^{K-1} c^{us}_l(k){\theta}^{us}_l  \bigg) \v_i(0) + \nonumber \\
	\epsilon \sum_{j \in \mathcal{N}_{US}}  \bigg(&\prod_{k=0}^{K-1}c^{us}_j(k){\theta}^{us}_j + \sum_{l \in \mathcal{N}_{S}} \sum_{r=0}^{K-1} \prod_{k=0}^{r-1} c^{us}_j(k) d_{j,l}(r)  \prod_{k=r+1}^{K-1} c^{s}_l(k){\theta}^{s}_l \nonumber \\ & + \sum_{l \in \mathcal{N}_{US}} \sum_{r=0}^{K-1} \prod_{k=0}^{r-1} c^{us}_j(k) d_{j,l}(r)  \prod_{k=r+1}^{K-1} c^{us}_l(k){\theta}^{us}_l  \bigg)  \v_j(0). \label{exactrajectory}
	\end{align}
	
	\subsubsection{Bounds on $\epsilon$:}
	Recall that from \eqref{approximation1} we established that the first-order approximation of the matrix product \eqref{matrixprod} is only valid for $\epsilon\frac{\norm{\P}_2}{\norm{\A}_2} < 1$ and $K \epsilon \ll1$. Next, from \eqref{binomialbound} we have that $ \sup_{0\leq k \leq K-1}\norm{\A_k}_2 = \norm{\A}_2$ and $\sup_{0\leq k \leq K-1}\norm{\P_k}_2 = \norm{\P}_2$.

    From \eqref{substitution} we have the following:
    \begin{align}
    \A_k &= \sum_{i \in \mathcal{N}_{S}}  c^{s}_i(k)\v_i(0)\v_i(0)^{T} + \sum_{j \in \mathcal{N}_{US}}  c^{us}_j(k)\v_j(0)\v_j(0)^{T}, \ \text{and} \label{expansionA}\\
    \B_k &= \sum_{i=1}^{n}  \sum_{l \neq i} \bigg(d_{l,i}(k)\v_l(0)\v_i(0)^{T}   +d_{i,l}(k)\v_i(0)\v_l(0)^{T}\bigg),
    \end{align}
     with $\epsilon \P_k = \B_k+\mathcal{O}(\epsilon^2)$.

     Observe that $\A_k$ is a matrix in its spectral decomposed form where the coefficients $c^{s}_i(k)$ and $c^{us}_i(k)$ correspond to the eigenvalues of $\A_k$. Therefore applying the bounds \eqref{bound1} and \eqref{bound2} we have the following result:
     \begin{align}
     \norm{\A}_2 &= \sup_{0\leq k \leq K-1}\norm{\A_k}_2 \\ & = \sup_{0\leq k \leq K-1}\bigg\{\max_{i \in \mathcal{N}_{S}, j \in \mathcal{N}_{US}}\{c^{s}_i(k), c^{us}_j(k)\}  \bigg\} \\
     & = \bigg(1 +\alpha L +\frac{\alpha \epsilon M}{2}+\mathcal{O}(\epsilon^2)   \bigg). \label{bound4}
     \end{align}
     Next, taking the norm of both sides of ${\epsilon \P_k}  = \B_k + \mathcal{O}(\epsilon^2)$, taking supremum over $k$ followed by the triangle inequality and then using \eqref{bound3} we get the following upper bound:
     \begin{align}
     \sup_{0\leq k \leq K-1}\norm{\epsilon \P_k}_2 & = \sup_{0\leq k \leq K-1}\norm{\B_k + \mathcal{O}(\epsilon^2)}_2 \\
     & \leq \sup_{0\leq k \leq K-1}\norm{\B_k}_2 + \mathcal{O}(\epsilon^2) \\
     & \leq \sum_{i=1}^{n}  \sum_{l \neq i}\bigg(\sup_{0\leq k \leq K-1}\norm{d_{l,i}(k)\v_l(0)\v_i(0)^{T}}_2   +\sup_{0\leq k \leq K-1}\norm{d_{i,l}(k)\v_i(0)\v_l(0)^{T}}_2\bigg) + \mathcal{O}(\epsilon^2) \\
     & \leq \sum_{i=1}^{n}  \sum_{l \neq i}\bigg(\sup_{0\leq k \leq K-1}\norm{d_{l,i}(k)\v_l(0)\v_i(0)^{T}}_F   +\sup_{0\leq k \leq K-1}\norm{d_{i,l}(k)\v_i(0)\v_l(0)^{T}}_F\bigg) + \mathcal{O}(\epsilon^2) \\
      &= \sum_{i=1}^{n}  \sum_{l \neq i}\bigg(\sup_{0\leq k \leq K-1}\abs{d_{l,i}(k)}   +\sup_{0\leq k \leq K-1}\abs{d_{i,l}(k)}\bigg) + \mathcal{O}(\epsilon^2) \\
      & \leq \frac{\alpha \epsilon M L n^2}{\delta} + \mathcal{O}(\epsilon^2), \label{bound5}
     \end{align}
     where in the last couple of steps we used the following properties of any matrix $\Z$: $     \norm{\Z}_2 \leq \norm{\Z}_F $, and $\norm{\Z}_F = \sqrt{\text{tr}(\Z \Z^T)}$.

     Now we require that $\epsilon\frac{\norm{\P}_2}{\norm{\A}_2} < 1$. Using \eqref{bound4}, this condition becomes
     \begin{align}
     \sup_{0\leq k \leq K-1}\norm{\epsilon \P_k}_2=\norm{\epsilon \P}_2 &< \norm{\A}_2 = \bigg(1 +\alpha L +\frac{\alpha \epsilon M}{2}+\mathcal{O}(\epsilon^2)   \bigg).
     \end{align}
     Therefore to obtain a bound on $\epsilon$ we can utilize \eqref{bound5} and set this condition as follows:
     \begin{align}
     \sup_{0\leq k \leq K-1}\norm{\epsilon \P_k}_2=\norm{\epsilon \P}_2 &\leq\frac{\alpha \epsilon M L n^2}{\delta} + \mathcal{O}(\epsilon^2)< \norm{\A}_2 = \bigg(1 +\alpha L +\frac{\alpha \epsilon M}{2}+\mathcal{O}(\epsilon^2)   \bigg) \\
     \frac{\alpha \epsilon M L n^2}{\delta} -\frac{\alpha \epsilon M}{2} & < 1+ \alpha L  + \mathcal{O}(\epsilon^2) \\
     \epsilon &< \frac{2 \delta(1+\alpha L)}{\alpha M ( 2 L n^2   -\delta)}  + \mathcal{O}(\epsilon^2). \label{condition1}
     \end{align}
     Note that this condition on $\epsilon$ is sufficient but may not be necessary since we are using an upper bound on $\sup_{0\leq k \leq K-1}\norm{\epsilon \P_k}_2$ from \eqref{bound5} as a lower bound for $\norm{\A}_2$. Hence, the inequality may shrink the feasible set for $\epsilon$ making it a sufficient condition but not necessary.

     Having established a range for $\epsilon$ from the upper bound \eqref{approximation1}, we utilize the lower bound \eqref{lowerboundmat} to get the complete feasible range for $\epsilon$.
     From the bound \eqref{lowerboundmat} we need that $\epsilon \sup_{0\leq k \leq K-1}\norm{ \A_k^{-1}  \P_k }_2 < 1$. Now for this particular condition to work, $\A_k$ should not have eigenvalues close to $0$ or of order $\mathcal{O}(\epsilon)$. Recall that from \eqref{expansionA}, $\A_k$ has its eigenvalues as $c_i^s(k)$ and $c_j^{us}(k)$ which are bounded by the inequalities in \eqref{bound1}, \eqref{bound2}. For $\alpha \approx \frac{1}{L}$, the lower bound in \eqref{bound1} becomes $\mathcal{O}(\epsilon)$. Hence we analyze the two cases corresponding to different ranges of $\alpha$ separately.
     \subsubsection{Case 1---$\alpha \in \bigg(0,\frac{1}{L}-\mathcal{O}(\epsilon)\bigg]$:}
     For this case, we can use the condition $\epsilon \sup_{0\leq k \leq K-1}\norm{ \A_k^{-1}  \P_k }_2 < 1$ in \eqref{lowerboundmat}.
       To obtain a certain feasible range on $\epsilon$, this condition can be set as follows:
     \begin{align}
    \epsilon \sup_{0\leq k \leq K-1}\norm{ \A_k^{-1}  \P_k }_2< \sup_{0\leq k \leq K-1}\norm{  \A_k^{-1}}_2     \sup_{0\leq k \leq K-1}\norm{\epsilon \P_k }_2 &< 1 \\
      \sup_{0\leq k \leq K-1}\bigg\{\max_{i \in \mathcal{N}_{S}, j \in \mathcal{N}_{US}}\bigg\{\frac{1}{c^{s}_i(k)}, \frac{1}{c^{us}_j(k)}\bigg\}  \bigg\} \sup_{0\leq k \leq K-1}\norm{\epsilon \P_k }_2&<1 \\
    \bigg(1-\alpha L -\frac{\alpha \epsilon M}{2}-\mathcal{O}(\epsilon^2)   \bigg)^{-1} \bigg( \frac{\alpha \epsilon M L n^2}{\delta} + \mathcal{O}(\epsilon^2)\bigg) & < 1 \\
     \frac{2 \delta(1-\alpha L)}{\alpha M ( 2 L n^2   +\delta)}  + \mathcal{O}(\epsilon^2) &>\epsilon.\label{condition2}
     \end{align}
      Note that this condition on $\epsilon$ is sufficient but may not be necessary.\\
      Moreover, combining the conditions \eqref{condition1} and \eqref{condition2} with \eqref{ROC} we get the following necessary bound:
     \begin{align}
     \epsilon < \min \bigg\{\inf_{{\norm{\u}=1}}\bigg(\limsup_{j\to \infty} \sqrt[j]{\frac{r_j(\u)}{j!}}\bigg)^{-1},\frac{2 \delta(1-\alpha L)}{\alpha M ( 2 L n^2   +\delta)}  + \mathcal{O}(\epsilon^2)\bigg\}.
     \end{align}
     Finally it is also required to have $K\epsilon \ll 1$ or $K \ll \frac{1}{\epsilon}$. Therefore this condition implies
     \begin{align}
     K = \mathcal{O}\bigg(\frac{1}{\epsilon}\bigg).
     \end{align}
      \subsubsection{Case 2---$\alpha \in \bigg(\frac{1}{L}-\mathcal{O}(\epsilon),\frac{1}{L}\bigg]:$ } For this case, observe that the lower bound in \eqref{lowerboundmat} is of order $\mathcal{O}(\epsilon^K)$. Further simplifying this lower bound and taking the infimum term inside, we obtain the following:
     \begin{align}
     	     \norm{\prod_{k=0}^{K-1} \bigg[\A_k + \epsilon \P_k   \bigg]}_2 & \geq  \bigg(\inf_{0\leq k \leq K-1}\norm{\A_k^{-1}}_2^{-1}\bigg)^K \bigg(1-\epsilon \sup_{0\leq k \leq K-1}\norm{ \A_k^{-1}  \P_k }_2  \bigg)^K \\
     	     & \geq  \bigg(\inf_{0\leq k \leq K-1}\norm{\A_k^{-1}}_2^{-1}\bigg)^K \bigg(1- \sup_{0\leq k \leq K-1}\norm{  \A_k^{-1}}_2     \sup_{0\leq k \leq K-1}\norm{\epsilon \P_k }_2 \bigg)^K \\
     	      & \geq  \bigg(\inf_{0\leq k \leq K-1}\norm{\A_k^{-1}}_2^{-1}- \frac{\sup_{0\leq k \leq K-1}\norm{  \A_k^{-1}}_2}{\sup_{0\leq k \leq K-1}\norm{  \A_k^{-1}}_2}     \sup_{0\leq k \leq K-1}\norm{\epsilon \P_k }_2 \bigg)^K \\
     	   & \geq  \bigg(\bigg|  \bigg(1-\alpha L -\frac{\alpha \epsilon M}{2}-\mathcal{O}(\epsilon^2)   \bigg)\bigg|- \bigg( \frac{\alpha \epsilon M L n^2}{\delta} + \mathcal{O}(\epsilon^2)\bigg) \bigg)^K.
     \end{align}
     Now for $\alpha = \frac{1}{L}$, the above lower bound will be $(C\epsilon)^K$ where $C$ is some constant. Therefore, for this lower bound to converge to $0$ for large $K$ we must necessarily have $C \epsilon < 1$ which implies
     \begin{align}
     \bigg| \frac{ \epsilon M}{2L} - \frac{\epsilon M n^2}{\delta} + \mathcal{O}(\epsilon^2) \bigg| & <1 \\
     \frac{1}{\frac{Mn^2}{\delta}- \frac{M}{2L} }+ \mathcal{O}(\epsilon^2) & >   \epsilon\\
     \frac{2L\delta}{M(2Ln^2 -\delta )}+ \mathcal{O}(\epsilon^2) &> \epsilon.
     \end{align}
     Finally, combining this condition on $\epsilon$ with \eqref{condition1} and \eqref{ROC} for $\alpha = \frac{1}{L}$, we get that
     \begin{align}
      \epsilon &< \min \bigg\{\inf_{{\norm{\u}=1}}\bigg(\limsup_{j\to \infty} \sqrt[j]{\frac{r_j(\u)}{j!}}\bigg)^{-1},\frac{4L\delta}{ M ( 2 L n^2   -\delta)}  + \mathcal{O}(\epsilon^2),\frac{2L\delta}{M(2Ln^2 -\delta )}+ \mathcal{O}(\epsilon^2) \bigg\}  \\
      \epsilon &< \min \bigg\{\inf_{{\norm{\u}=1}}\bigg(\limsup_{j\to \infty} \sqrt[j]{\frac{r_j(\u)}{j!}}\bigg)^{-1},\frac{2L\delta}{M(2Ln^2 -\delta )} + \mathcal{O}(\epsilon^2)\bigg\}.
     \end{align}
     The condition $K = \mathcal{O}\bigg(\frac{1}{\epsilon}\bigg)$ is still required to hold.
     \endproof

     \section{Lower bounds on the distance between $\x_{K}$ and $\x^*$}
     \label{Appendix C}
     \subsection{\textbf{Proof of Lemma \ref{lem7}}}
     \proof{}
     An approximate equation for $\u_{K}$ in terms of $\u_{0}$ is given by \eqref{exactrajectory}. This approximation holds for all values of $K$ from $1$ to $K_{exit}$, where $K_{exit}$ denotes the iteration number of escape from $\mathcal{B}_{\epsilon}(\x^*)$. Formally $K_{exit}$ can be expressed as
     \begin{align}
     K_{exit} & = \inf_{K\geq 1} \bigg\{K \hspace{0.2cm}\bigg| \hspace{0.2cm} \norm{\tilde{\u}_K}^2 > \epsilon^2 \bigg\},   \label{exittime}
     \end{align} where the squared norm is used for the sake of simplifying subsequent analysis involving lower bounds. However, the sequence $\{\tilde{\u}_K\}_{K=0}^{K_{exit}}$ cannot be determined solely from the initialization $\u_{0}$. To uniquely determine any $\tilde{\u}_K$, we still need to know the coefficient terms $c^{s}_i(k)$, $c^{us}_j(k)$ and $d_{l,i}(k)$ for all values of $k$ from $0$ to $K-1$. The only information available in this regard is the bound on these coefficients from \eqref{bound1}, \eqref{bound2} and \eqref{bound3}. Therefore it becomes impossible to predetermine the entire sequence $\{\tilde{\u}_K \}_{K=0}^{K_{exit}}$ just based on the knowledge of $\u_{0}$ .

     To circumvent this problem, we introduce a set $S_{\epsilon}$ which is the set of all possible $\epsilon$ precision trajectories generated by the approximate equation \eqref{exactrajectory}. Recall that while deriving the approximation \eqref{exactrajectory}, we expanded terms appearing in the product \eqref{matrixprod} only up to order $\mathcal{O}(\epsilon)$; hence we can call these approximate sequences as $\epsilon$ precision trajectories with respect to $\x^*$. For a fixed initialization of $\u_{0}$, the set $S_{\epsilon}$ is given by
     \begin{align}
     S_{\epsilon} & = \bigg\{  \{\tilde{\u}_K^{\tau}\}_{K=1}^{K_{exit}^{\tau}}  \hspace{0.2cm}\bigg| \hspace{0.2cm}\u_{0} \bigg \},
     \end{align}
     where each possible $\epsilon$ precision trajectory is parameterized by some $\tau \in \mathbb{R}$, $K_{exit}^{\tau}$ is the escape iteration for the $\tau$-parameterized $\epsilon$-precision trajectory and $\tilde{\u}_K^{\tau}$ satisfies \eqref{exactrajectory} for every $\tau$. Note that $\tau$ varies with variations in the sequence $\bigg \{\{c^{s}_i(k), c^{us}_j(k), d_{l,i}(k)\}_{k=0}^{K-1} \bigg \}_{K=1}^{K_{exit}}$ which are in turn controlled by variations in the coefficient terms from the bounds \eqref{bound1}, \eqref{bound2} and \eqref{bound3}. Since the set $S_{\epsilon}$ contains all possible $\epsilon$ precision trajectories, the actual $\epsilon$-precise trajectory that the radial vector $\u_{K}$ takes inside the ball $\mathcal{B}_{\epsilon}(\x^*)$ will also belong to the set $S_{\epsilon}$. Let this actual $\epsilon$ precision trajectory be parameterized by some $\tau = \omega$. Therefore we have that
     \begin{align}
\{\tilde{\u}_K^{\omega}\}_{K=1}^{K_{exit}^{\omega}} \in S_{\epsilon}.
\end{align} Moreover, $\tilde{\u}_K^{\omega}$ satisfies the approximate equation \eqref{exactrajectory}. Next using \eqref{exittime}, we can write the escape iteration for the $\tau$-parameterized $\epsilon$ precision trajectory as
 \begin{align}
 K_{exit}^{\tau} & = \inf_{K\geq 1} \bigg\{K \hspace{0.2cm}\bigg| \hspace{0.2cm} \norm{\tilde{\u}_K^{\tau}}^2 > \epsilon^2 \bigg\}. \label{exittimetau}
 \end{align}
 We now define a quantity $K^{\iota}$ such that
 \begin{align}
 K^{\iota} & = \inf_{K\geq 1} \bigg\{K \hspace{0.2cm}\bigg| \hspace{0.2cm} \inf_{\tau}\bigg\{\norm{\tilde{\u}_K^{\tau}}^2 \bigg \}> \epsilon^2 \bigg\}.  \label{upperexittime}
 \end{align}
\subsubsection{Claim of the lemma:}
 \begin{align}
  K^{\iota} &\geq \sup_{\tau} \bigg\{K_{exit}^{\tau} \bigg\} \label{exittimebound} = \sup_{\tau} \inf_{K\geq 1} \bigg\{K \hspace{0.2cm}\bigg| \hspace{0.2cm} \norm{\tilde{\u}_K^{\tau}}^2 > \epsilon^2 \bigg\}.
 \end{align}
 \textit{Proof by contradiction:}
 Let us assume that for some $\tau = a$ the escape iteration $K_{exit}^{a}$ is such that $K_{exit}^{a} > K^{\iota} $. From the definition of $K^{\iota}$ in \eqref{upperexittime}, $K^{\iota}$ is the smallest iteration such that $\inf_{\tau}\bigg\{\norm{\tilde{\u}_{K^{\iota}}^{\tau}}^2 \bigg \}> \epsilon^2 $. This implies $ \norm{\tilde{\u}_{K^{\iota}}^a}^2 > \epsilon^2$. However, this is not possible since it contradicts the definition of infimum from \eqref{exittimetau} for $\tau = a$. Therefore we must have $K_{exit}^{a} \leq K^{\iota} $ and this should hold for any $a$. Hence, we must have $ K^{\iota} \geq \sup_{\tau} \bigg\{K_{exit}^{\tau} \bigg\}$.

 Since the actual $\epsilon$-precise trajectory given by $\{\tilde{\u}_K^{\omega}\}_{K=1}^{K_{exit}^{\omega}}$ belongs to the $\tau$-parameterized set $S_{\epsilon}$, hence $K_{exit}^{\omega} \leq K^{\iota}$. Therefore it is sufficient to develop an upper bound on $K^{\iota}$ in order to draw conclusions about $K_{exit}^{\omega}$. In the subsequent section, we analyze the lower bound on $\norm{\tilde{\u}_K}^2$ to obtain this $K^{\iota}$.
     \endproof

     \subsection{\textbf{Proof of Theorem \ref{theorem1}}}
     \proof{}

     Taking the norm squared on both sides of \eqref{exactrajectory} we get the following:
     \begin{align}
     \norm{\tilde{\u}_K}^2 = \epsilon^2\sum_{i \in \mathcal{N}_{S}}  \bigg(&\underbrace{\prod_{k=0}^{K-1}c^{s}_i(k){\theta}^{s}_i}_{T_1} + \underbrace{\sum_{l \in \mathcal{N}_{S}} \sum_{r=0}^{K-1} \prod_{k=0}^{r-1} c^{s}_i(k) d_{i,l}(r)  \prod_{k=r+1}^{K-1} c^{s}_l(k){\theta}^{s}_l}_{T_2} \nonumber \\ & + \underbrace{\sum_{l \in \mathcal{N}_{US}} \sum_{r=0}^{K-1} \prod_{k=0}^{r-1} c^{s}_i(k) d_{i,l}(r)  \prod_{k=r+1}^{K-1} c^{us}_l(k){\theta}^{us}_l}_{T_3}  \bigg)^2  + \nonumber \\
     \epsilon^2 \sum_{j \in \mathcal{N}_{US}}  \bigg(&\underbrace{\prod_{k=0}^{K-1}c^{us}_j(k){\theta}^{us}_j}_{T_4} + \underbrace{\sum_{l \in \mathcal{N}_{S}} \sum_{r=0}^{K-1} \prod_{k=0}^{r-1} c^{us}_j(k) d_{j,l}(r)  \prod_{k=r+1}^{K-1} c^{s}_l(k){\theta}^{s}_l}_{T_5} \nonumber \\ & + \underbrace{\sum_{l \in \mathcal{N}_{US}} \sum_{r=0}^{K-1} \prod_{k=0}^{r-1} c^{us}_j(k) d_{j,l}(r)  \prod_{k=r+1}^{K-1} c^{us}_l(k){\theta}^{us}_l}_{T_6}  \bigg)^2 \\
      =  \epsilon^2\bigg(&\sum_{i \in \mathcal{N}_{S}} (T_1 + T_2 + T_3)^2 + \sum_{j \in \mathcal{N}_{US}}(T_4 + T_5 + T_6)^2\bigg).
     \end{align}
     Now this equation is satisfied by $\tilde{\u}_K^{\tau}$ for every $\tau$. Hence for any given $\tau$ we can write
     \begin{align}
     \norm{\tilde{\u}_K^{\tau}}^2 =  \epsilon^2\bigg(&\sum_{i \in \mathcal{N}_{S}} (T_1(\tau) + T_2(\tau) + T_3(\tau))^2 + \sum_{j \in \mathcal{N}_{US}}(T_4(\tau) + T_5(\tau) + T_6(\tau))^2\bigg), \label{trajectoryterms}
     \end{align}
     where $\tau$ varies with variations in the sequence $\bigg \{\{c^{s}_i(k), c^{us}_j(k), d_{l,i}(k)\}_{k=0}^{K-1} \bigg \}_{K=1}^{K_{exit}}$.

     Using \eqref{bound1}, \eqref{bound2} and \eqref{bound3} we get the bounds on these coefficient product terms from $T_1(\tau)$ to $T_6(\tau)$. Starting with the term $T_1(\tau)$ we have that
     \begin{align}
      \inf_{\tau}T_1(\tau)  & = \prod_{k=0}^{K-1}\inf_{\tau} \bigg \{c^{s}_i(k)\bigg \}{\theta}^{s}_i =    \bigg(1 - \alpha L - \frac{\alpha \epsilon M}{2}-\mathcal{O}(\epsilon^2) \bigg)^K {\theta}^{s}_i, \ \text{and}
     \end{align}
     \begin{align}
     \sup_{\tau}T_1(\tau) &= \prod_{k=0}^{K-1}\sup_{\tau}  \bigg \{c^{s}_i(k) \bigg \}{\theta}^{s}_i = \bigg(1 - \alpha \beta + \frac{\alpha \epsilon M}{2} + \mathcal{O}(\epsilon^2) \bigg)^K {\theta}^{s}_i
     \end{align}
     for positive $c^{s}_i(k)$. Next for the term $T_2(\tau)$, first consider the lower bound
    \begin{align}
   \inf_{\tau} T_2(\tau)  & \geq \sum_{l \in \mathcal{N}_{S}} \sum_{r=0}^{K-1}  \inf_{\tau} \bigg \{{d_{i,l}(r)}\prod_{k=0}^{r-1} c^{s}_i(k)  \prod_{k=r+1}^{K-1} c^{s}_l(k)\bigg \}{\theta}^{s}_l  \\
    & \geq \sum_{l \in \mathcal{N}_{S}} \sum_{r=0}^{K-1} - \sup_{\tau} \bigg \{\abs{d_{i,l}(r)}\bigg \}\sup_{\tau} \bigg \{\prod_{k=0}^{r-1} c^{s}_i(k)  \prod_{k=r+1}^{K-1} c^{s}_l(k)\bigg \}{\theta}^{s}_l \\
    & \geq \sum_{l \in \mathcal{N}_{S}} \sum_{r=0}^{K-1} - \bigg(\frac{\alpha \epsilon M L}{2 \delta} + \mathcal{O}(\epsilon^2) \bigg)\sup_{\tau} \bigg \{\prod_{k=0}^{r-1} c^{s}_i(k)  \prod_{k=r+1}^{K-1}  c^{s}_l(k)\bigg \}{\theta}^{s}_l \\
    & =
    \sum_{l \in \mathcal{N}_{S}} \sum_{r=0}^{K-1} - \bigg(\frac{\alpha \epsilon M L}{2 \delta} + \mathcal{O}(\epsilon^2) \bigg)\prod_{k=0}^{r-1} \sup_{\tau} \bigg \{c^{s}_i(k) \bigg \} \prod_{k=r+1}^{K-1} \sup_{\tau} \bigg \{ c^{s}_l(k)\bigg \}{\theta}^{s}_l \\
    & =  	-K\bigg(1 - \alpha \beta + \frac{\alpha \epsilon M}{2} + \mathcal{O}(\epsilon^2)\bigg)^{K-1} \bigg(\frac{\alpha \epsilon M L}{2 \delta} + \mathcal{O}(\epsilon^2) \bigg) \sum_{l \in \mathcal{N}_{S}} {\theta}^{s}_l,
    \end{align}
    where we have $c_i^s(k)\geq 0$ for all $i$ and $k$. The upper bound on $T_2(\tau)$ is as follows:
      \begin{align}
    \sup_{\tau} T_2(\tau)  & \leq \sum_{l \in \mathcal{N}_{S}} \sum_{r=0}^{K-1}  \sup_{\tau} \bigg \{d_{i,l}(r)\prod_{k=0}^{r-1} c^{s}_i(k)  \prod_{k=r+1}^{K-1} c^{s}_l(k)\bigg \}{\theta}^{s}_l  \\
    & \leq \sum_{l \in \mathcal{N}_{S}} \sum_{r=0}^{K-1}  \sup_{\tau} \bigg \{\abs{d_{i,l}(r)}\bigg \}\sup_{\tau} \bigg \{\prod_{k=0}^{r-1} c^{s}_i(k)  \prod_{k=r+1}^{K-1} c^{s}_l(k)\bigg \}{\theta}^{s}_l \\
    & = \sum_{l \in \mathcal{N}_{S}} \sum_{r=0}^{K-1}  \bigg(\frac{\alpha \epsilon M L}{2 \delta} + \mathcal{O}(\epsilon^2) \bigg)\sup_{\tau} \bigg \{\prod_{k=0}^{r-1} c^{s}_i(k)  \prod_{k=r+1}^{K-1}  c^{s}_l(k)\bigg \}{\theta}^{s}_l \\
    & =
 K\bigg(1 - \alpha \beta + \frac{\alpha \epsilon M}{2} + \mathcal{O}(\epsilon^2)\bigg)^{K-1} \bigg(\frac{\alpha \epsilon M L}{2 \delta} + \mathcal{O}(\epsilon^2) \bigg) \sum_{l \in \mathcal{N}_{S}} {\theta}^{s}_l.
    \end{align}
     For the term $T_3(\tau)$, first consider the lower bound
     \begin{align}
\inf_{\tau} T_3(\tau)     & \geq \sum_{l \in \mathcal{N}_{US}} \sum_{r=0}^{K-1} \inf_{\tau} \bigg \{ \prod_{k=0}^{r-1}  c^{s}_i(k)  d_{i,l}(r)   \prod_{k=r+1}^{K-1}  c^{us}_l(k)\bigg\}{\theta}^{us}_l \\
    & \geq \sum_{l \in \mathcal{N}_{US}} \sum_{r=0}^{K-1}  -\sup_{\tau} \bigg\{\abs{d_{i,l}(r)}\bigg\} \sup \bigg\{\prod_{k=0}^{r-1} c^{s}_i(k)   \prod_{k=r+1}^{K-1} c^{us}_l(k)\bigg\}{\theta}^{us}_l     \\
      & = \sum_{l \in \mathcal{N}_{US}} \sum_{r=0}^{K-1}- \bigg(\frac{\alpha \epsilon M L}{2 \delta} + \mathcal{O}(\epsilon^2) \bigg) \prod_{k=0}^{r-1} \sup_{\tau} \bigg\{c^{s}_i(k)\bigg\}   \prod_{k=r+1}^{K-1} \sup_{\tau} \bigg\{c^{us}_l(k)\bigg\}{\theta}^{us}_l     \\
      & = \sum_{l \in \mathcal{N}_{US}} \sum_{r=0}^{K-1}- \bigg(\frac{\alpha \epsilon M L}{2 \delta} + \mathcal{O}(\epsilon^2) \bigg)  \bigg(1 -\alpha \beta +\frac{\alpha \epsilon M}{2}+\mathcal{O}(\epsilon^2)   \bigg)^{r}  \bigg(1 +\alpha L +\frac{\alpha \epsilon M}{2}+\mathcal{O}(\epsilon^2)   \bigg)^{K-r-1}{\theta}^{us}_l     \\
           &  =  -\bigg(\frac{\alpha \epsilon M L}{2 \delta} + \mathcal{O}(\epsilon^2) \bigg)\frac{\bigg(1 + \alpha L + \frac{\alpha \epsilon M}{2} + \mathcal{O}(\epsilon^2)\bigg)^K - \bigg(1 - \alpha \beta - \frac{\alpha \epsilon M}{2}+ \mathcal{O}(\epsilon^2)\bigg)^K}{ ( \alpha L+ \alpha \beta + \mathcal{O}(\epsilon^2))} \sum_{l \in \mathcal{N}_{US}} {\theta}^{us}_l     \\
    & >     -\bigg(\frac{ \alpha \epsilon M L}{2 \delta} + \mathcal{O}(\epsilon^2) \bigg)\frac{\bigg(1 + \alpha L + \frac{\alpha \epsilon M}{2} + \mathcal{O}(\epsilon^2)\bigg)^K }{ ( \alpha L+\alpha\beta+ \mathcal{O}(\epsilon^2))} \sum_{l \in \mathcal{N}_{US}} {\theta}^{us}_l.
     \end{align}
     Note that here in the last step we used a loose lower bound by dropping the negative term from the numerator for the sake of simplifying the subsequent analysis.
     Similarly, an upper bound on $T_3(\tau)$ can be obtained, which is as follows:
     \begin{align}
     \sup_{\tau}T_3(\tau) & <  \bigg(\frac{ \alpha\epsilon M L}{2 \delta} + \mathcal{O}(\epsilon^2) \bigg)\frac{\bigg(1 + \alpha L + \frac{\alpha \epsilon M}{2} + \mathcal{O}(\epsilon^2)\bigg)^K }{ ( \alpha L+\alpha\beta + \mathcal{O}(\epsilon^2))} \sum_{l \in \mathcal{N}_{US}} {\theta}^{us}_l.
     \end{align}
     Now that we have derived the bounds for the terms $T_1(\tau), T_2(\tau), T_3(\tau)$, the bounds for remaining terms $T_4(\tau), T_5(\tau), T_6(\tau)$ can be derived along similar lines. Since the algebra is somewhat tedious, we leave these derivations to the reader and directly present the bounds.

     The term $T_4(\tau)$ is bounded as
     \begin{align}
    \inf_{\tau} T_4(\tau) & =     \bigg(1 + \alpha \beta - \frac{\alpha \epsilon M}{2} -\mathcal{O}(\epsilon^2)\bigg)^K {\theta}^{us}_j, \ \text{and}
     \end{align}
     \begin{align}
 \sup_{\tau}T_4(\tau) & = \bigg(1 +\alpha L + \frac{\alpha \epsilon M}{2}+\mathcal{O}(\epsilon^2)\bigg)^K {\theta}^{us}_j.
     \end{align}
   The lower and upper bound on term $T_5(\tau)$ are as follows:
   \begin{align}
   \inf_{\tau}T_5(\tau) &  >	 -\bigg(\frac{\alpha\epsilon M L}{2 \delta} +\mathcal{O}(\epsilon^2)\bigg) \frac{\bigg(1 + \alpha L + \frac{\alpha \epsilon M}{2}+\mathcal{O}(\epsilon^2)\bigg)^K}{(\alpha L+\alpha \beta + \mathcal{O}(\epsilon^2))} \sum_{l \in \mathcal{N}_{S}} {\theta}^{s}_l, \ \text{and} \\
    \sup_{\tau}T_5(\tau) &   < \bigg(\frac{\alpha\epsilon M L}{2 \delta} +\mathcal{O}(\epsilon^2)\bigg) \frac{\bigg(1 + \alpha L + \frac{\alpha \epsilon M}{2}+\mathcal{O}(\epsilon^2)\bigg)^K}{(\alpha L + \alpha\beta + \mathcal{O}(\epsilon^2))} \sum_{l \in \mathcal{N}_{S}} {\theta}^{s}_l.
   \end{align}
    The lower and upper bound on term $T_6(\tau)$ are as follows:
     \begin{align}
    \inf_{\tau}T_6(\tau) & \geq -K\bigg(1 + \alpha L + \frac{\alpha \epsilon M}{2} + \mathcal{O}(\epsilon^2)\bigg)^{K-1} \bigg(\frac{\alpha \epsilon M L}{2 \delta} + \mathcal{O}(\epsilon^2) \bigg) \sum_{l \in \mathcal{N}_{US}} {\theta}^{us}_l, \ \text{and}  \\
     \sup_{\tau}T_6(\tau) &  \leq  K\bigg(1 + \alpha L + \frac{\alpha \epsilon M}{2} + \mathcal{O}(\epsilon^2)\bigg)^{K-1} \bigg(\frac{\alpha \epsilon M L}{2 \delta} + \mathcal{O}(\epsilon^2) \bigg) \sum_{l \in \mathcal{N}_{US}} {\theta}^{us}_l.
     \end{align}
     Using these results and dropping higher order terms ($\mathcal{O}(\epsilon^2)$ and above), we can get the lower bound on $\norm{\tilde{\u}_K^{\tau}}^2$. From \eqref{trajectoryterms}, observe that
     \begin{align}
\norm{\tilde{\u}_K}^2 = & \epsilon^2\bigg(\sum_{i \in \mathcal{N}_{S}} (T_1 + T_2 + T_3)^2 + \sum_{j \in \mathcal{N}_{US}}(T_4 + T_5 + T_6)^2\bigg).
     \end{align}
     Let $Y_1(\tau) = \sum_{i \in \mathcal{N}_{S}} (T_1(\tau) + T_2(\tau) + T_3(\tau))^2$ and $Y_2(\tau) = \sum_{j \in \mathcal{N}_{US}} (T_4(\tau) + T_5(\tau) + T_6(\tau))^2$. Using \eqref{trajectoryterms}, we can see that
     \begin{align}
\norm{\tilde{\u}_K^{\tau}}^2 = & \epsilon^2\bigg(Y_1(\tau) + Y_2(\tau)\bigg).
     \end{align} Now using the bounds for $T_1(\tau), T_2(\tau), T_3(\tau)$ we have the following lower bound on $Y_1(\tau)$:
     \begin{align}
   \hspace{-0.2cm}  \inf_{\tau}Y_1(\tau) & \geq \sum_{i \in \mathcal{N}_{S}} \bigg(\inf_{\tau} \bigg \{
     T_1^2(\tau) + T_2^2(\tau) + T_3^2(\tau) + 2 T_1(\tau)T_2(\tau) + 2T_2(\tau)T_3(\tau) + 2T_3(\tau) T_1(\tau)  \bigg \} \bigg)\\
      & \hspace{-0.5cm} \geq \sum_{i \in \mathcal{N}_{S}} \bigg(\inf_{\tau}
     \underbrace{T_1^2(\tau)}_{>0} + \inf_{\tau} \underbrace{T_2^2(\tau)}_{\geq 0} + \inf_{\tau} \underbrace{T_3^2(\tau)}_{\geq 0} + 2  \underbrace{\inf_{\tau}\bigg\{ T_1(\tau)T_2(\tau)\bigg\}}_{<0} + 2  \underbrace{\inf_{\tau}\bigg\{ T_2(\tau)T_3(\tau)\bigg\}}_{<0}  + 2  \underbrace{\inf_{\tau}\bigg\{ T_3(\tau)T_1(\tau)\bigg\}}_{<0}\bigg) \\
           & \hspace{-0.5cm} >\bigg( \sum_{i \in \mathcal{N}_{S}} \inf_{\tau}
     \underbrace{T_1^2(\tau)}_{>0} + 0 + 0 - 2  \underbrace{\sup_{\tau}\abs{T_1(\tau)}}_{> 0} \underbrace{\sup_{\tau}\abs{T_2(\tau)}}_{>0}- 2  \underbrace{\sup_{\tau}\abs{T_2(\tau)}}_{> 0} \underbrace{\sup_{\tau}\abs{T_3(\tau)}}_{>0}- 2  \underbrace{\sup_{\tau}\abs{T_3(\tau)}}_{> 0} \underbrace{\sup_{\tau}\abs{T_1(\tau)}}_{>0}  \bigg) \\
 & \hspace{-0.5cm} = \sum_{i \in \mathcal{N}_{S}} \bigg(   \bigg(1 - \alpha L - \frac{\alpha \epsilon M}{2}-\mathcal{O}(\epsilon^2) \bigg)^{2K} ({\theta}^{s}_i)^2  - 2 K{\theta}^{s}_i\bigg(1 - \alpha \beta + \frac{\alpha \epsilon M}{2} + \mathcal{O}(\epsilon^2)\bigg)^{2K-1} \bigg(\frac{\alpha \epsilon M L}{2 \delta} +\mathcal{O}(\epsilon^2) \bigg) \sum_{l \in \mathcal{N}_{S}} {\theta}^{s}_l \nonumber \\ & \hspace{-0.5cm} -2K \bigg(1 - \alpha \beta + \frac{\alpha \epsilon M}{2} + \mathcal{O}(\epsilon^2) \bigg)^{K-1}  \bigg(\frac{ \alpha \epsilon M L}{2 \delta} +\mathcal{O}(\epsilon^2)\bigg)^2\frac{\bigg(1 + \alpha L + \frac{\alpha \epsilon M}{2} + \mathcal{O}(\epsilon^2)\bigg)^K }{ ( \alpha L+\alpha \beta + \mathcal{O}(\epsilon^2))} \sum_{l \in \mathcal{N}_{US}} {\theta}^{us}_l \sum_{i \in \mathcal{N}_{S}} {\theta}^{s}_i  \bigg)\nonumber \\ & \hspace{-0.5cm} -2 \bigg(1 - \alpha \beta + \frac{\alpha \epsilon M}{2} + \mathcal{O}(\epsilon^2) \bigg)^K  \bigg(\frac{ \alpha \epsilon M L}{2 \delta} +\mathcal{O}(\epsilon^2)\bigg)\frac{\bigg(1 + \alpha L + \frac{\alpha \epsilon M}{2} + \mathcal{O}(\epsilon^2)\bigg)^K }{ ( \alpha L+\alpha \beta + \mathcal{O}(\epsilon^2))} \sum_{l \in \mathcal{N}_{US}} {\theta}^{us}_l \sum_{i \in \mathcal{N}_{S}} {\theta}^{s}_i  \\
 & \hspace{-0.5cm} =  \sum_{i \in \mathcal{N}_{S}} \bigg(   \bigg(1 - \alpha L - \frac{\alpha \epsilon M}{2}-\mathcal{O}(\epsilon^2) \bigg)^{2K} ({\theta}^{s}_i)^2  - 2 K{\theta}^{s}_i\bigg(1 - \alpha \beta + \frac{\alpha \epsilon M}{2} + \mathcal{O}(\epsilon^2)\bigg)^{2K-1} \bigg(\frac{\alpha \epsilon M L}{2 \delta} +\mathcal{O}(\epsilon^2) \bigg) \sum_{l \in \mathcal{N}_{S}} {\theta}^{s}_l \bigg)\nonumber \\ & \hspace{-3cm}-2 \bigg(1+\frac{K\sum_{i \in \mathcal{N}_{S}}\bigg(\frac{ \alpha \epsilon M L}{2 \delta}+\mathcal{O}(\epsilon^2) \bigg)}{\bigg(1 - \alpha \beta + \frac{\alpha \epsilon M}{2} + \mathcal{O}(\epsilon^2) \bigg)}\bigg)\bigg(1 - \alpha \beta + \frac{\alpha \epsilon M}{2} + \mathcal{O}(\epsilon^2) \bigg)^K  \bigg(\frac{ \alpha \epsilon M L}{2 \delta}+\mathcal{O}(\epsilon^2) \bigg)\frac{\bigg(1 + \alpha L + \frac{\alpha \epsilon M}{2} + \mathcal{O}(\epsilon^2)\bigg)^K }{ ( \alpha L+\alpha \beta + \mathcal{O}(\epsilon^2))} \sum_{l \in \mathcal{N}_{US}} {\theta}^{us}_l \sum_{i \in \mathcal{N}_{S}} {\theta}^{s}_i \\
  &\hspace{-0.5cm} =  \sum_{i \in \mathcal{N}_{S}} \bigg(   \bigg(1 - \alpha L - \frac{\alpha \epsilon M}{2}-\mathcal{O}(\epsilon^2) \bigg)^{2K} ({\theta}^{s}_i)^2  - 2 K{\theta}^{s}_i\bigg(1 - \alpha \beta + \frac{\alpha \epsilon M}{2} + \mathcal{O}(\epsilon^2)\bigg)^{2K-1} \bigg(\frac{\alpha \epsilon M L}{2 \delta} +\mathcal{O}(\epsilon^2) \bigg) \sum_{l \in \mathcal{N}_{S}} {\theta}^{s}_l \bigg)\nonumber \\ &\hspace{-0.5cm} -2 \bigg(1+\mathcal{O}(K\epsilon)\bigg)\bigg(1 - \alpha \beta + \frac{\alpha \epsilon M}{2} + \mathcal{O}(\epsilon^2) \bigg)^K  \bigg(\frac{ \alpha \epsilon M L}{2 \delta}+\mathcal{O}(\epsilon^2) \bigg)\frac{\bigg(1 + \alpha L + \frac{\alpha \epsilon M}{2} + \mathcal{O}(\epsilon^2)\bigg)^K }{ ( \alpha L+\alpha \beta + \mathcal{O}(\epsilon^2))} \sum_{l \in \mathcal{N}_{US}} {\theta}^{us}_l \sum_{i \in \mathcal{N}_{S}} {\theta}^{s}_i,
     \end{align}
     where in the last step we replaced the term $K\sum_{i \in \mathcal{N}_{S}}\frac{\bigg(\frac{ \alpha \epsilon M L}{2 \delta}+\mathcal{O}(\epsilon^2) \bigg)}{\bigg(1 - \alpha \beta + \frac{\alpha \epsilon M}{2} + \mathcal{O}(\epsilon^2) \bigg)}$ with $\mathcal{O}(K\epsilon)$ for $K\epsilon \ll 1$ and $\bigg(1 - \alpha \beta + \frac{\alpha \epsilon M}{2} + \mathcal{O}(\epsilon^2) \bigg) \gg \epsilon$. This is because the numerator $ \bigg(\frac{ \alpha \epsilon M L}{2 \delta}+\mathcal{O}(\epsilon^2) \bigg)$ is of $\mathcal{O}(\epsilon)$; hence, we require the denominator $\bigg(1 - \alpha \beta + \frac{\alpha \epsilon M}{2} + \mathcal{O}(\epsilon^2) \bigg)$ to be of constant order, i.e., independent of $\epsilon$. Similarly, using the bounds for $T_4(\tau), T_5(\tau), T_6(\tau)$ we have the following lower bound for $Y_2(\tau)$:

       \begin{align}
 \hspace{-0.2cm}\inf_{\tau}Y_2(\tau) & \geq \sum_{j \in \mathcal{N}_{US}} \bigg(\inf_{\tau} \bigg \{
      T_4^2(\tau) + T_5^2(\tau) + T_6^2(\tau) + 2 T_4(\tau)T_6(\tau) + 2T_6(\tau)T_5(\tau) + 2T_5(\tau) T_4(\tau)  \bigg \} \bigg)\\
      & \geq \sum_{j \in \mathcal{N}_{US}} \bigg(\inf_{\tau}
      \underbrace{T_4^2(\tau)}_{>0} + \inf_{\tau} \underbrace{T_5^2(\tau)}_{\geq 0} + \inf_{\tau} \underbrace{T_6^2(\tau)}_{\geq 0} + 2  \underbrace{\inf_{\tau}\bigg\{T_4(\tau)T_6(\tau)\bigg\}}_{<0} + 2  \underbrace{\inf_{\tau}\bigg\{T_6(\tau)T_5(\tau)\bigg\}}_{<0}  + 2  \underbrace{\inf_{\tau}\bigg\{T_5(\tau)T_4(\tau)\bigg\}}_{<0} \bigg) \\
        & >\bigg( \sum_{j \in \mathcal{N}_{US}} \inf_{\tau}
      \underbrace{T_4^2(\tau)}_{>0} + 0 + 0 - 2  \underbrace{\sup_{\tau}\abs{T_4(\tau)}}_{> 0} \underbrace{\sup_{\tau}\abs{T_6(\tau)}}_{>0}- 2  \underbrace{\sup_{\tau}\abs{T_6(\tau)}}_{> 0} \underbrace{\sup_{\tau}\abs{T_5(\tau)}}_{>0}- 2  \underbrace{\sup_{\tau}\abs{T_5(\tau)}}_{> 0} \underbrace{\sup_{\tau}\abs{T_4(\tau)}}_{>0}  \bigg) \\
      & \hspace{-1.5cm}= \sum_{j \in \mathcal{N}_{US}} \bigg(   \bigg(1 + \alpha \beta - \frac{\alpha \epsilon M}{2}-\mathcal{O}(\epsilon^2) \bigg)^{2K} ({\theta}^{us}_j)^2  - 2 K{\theta}^{us}_j\bigg(1 + \alpha L + \frac{\alpha \epsilon M}{2} + \mathcal{O}(\epsilon^2)\bigg)^{2K-1} \bigg(\frac{\alpha \epsilon M L}{2 \delta} +\mathcal{O}(\epsilon^2) \bigg) \sum_{l \in \mathcal{N}_{US}} {\theta}^{us}_l \bigg) \nonumber \\ & -2   \bigg(1+\mathcal{O}(K\epsilon)\bigg)\bigg(\frac{ \alpha \epsilon M L}{2 \delta}+\mathcal{O}(\epsilon^2) \bigg)\frac{\bigg(1 + \alpha L + \frac{\alpha \epsilon M}{2} + \mathcal{O}(\epsilon^2)\bigg)^{2K} }{ ( \alpha L+\alpha \beta + \mathcal{O}(\epsilon^2))} \sum_{l \in \mathcal{N}_{S}} {\theta}^{s}_l \sum_{j \in \mathcal{N}_{US}}{\theta}^{us}_j.
      \end{align}
      Finally combining these two bounds yields the following lower bound on $\inf_{\tau} \norm{\tilde{\u}^{\tau}_K}^2 $:
     \begin{align}
     \inf_{\tau}\norm{\tilde{\u}_K^{\tau}}^2 = &   \epsilon^2\bigg(\inf_{\tau}Y_1(\tau) + \inf_{\tau}Y_2(\tau)\bigg) \\
        & \hspace{-1.5cm} >\epsilon^2\bigg[\bigg(1 - \alpha L - \frac{\alpha \epsilon M}{2} -\mathcal{O}(\epsilon^2)\bigg)^{2K} \sum_{i \in \mathcal{N}_{S}} ({\theta}^{s}_i)^2  - 2K\bigg(1 - \alpha \beta + \frac{\alpha \epsilon M}{2}+\mathcal{O}(\epsilon^2)\bigg)^{2K-1}  \bigg(\frac{ \alpha \epsilon M L}{2 \delta}+\mathcal{O}(\epsilon^2) \bigg)  (\sum_{i \in \mathcal{N}_{S}} {\theta}^{s}_i)^2 + \nonumber \\ & \hspace{-1.5cm}\bigg(1 + \alpha \beta - \frac{\alpha \epsilon M}{2}-\mathcal{O}(\epsilon^2)\bigg)^{2K}\sum_{j \in \mathcal{N}_{US}}( {\theta}^{us}_j)^2   -  2K\bigg(1 + \alpha L + \frac{\alpha \epsilon M}{2}+\mathcal{O}(\epsilon^2)\bigg)^{2K-1}  \bigg(\frac{ \alpha \epsilon M L}{2 \delta}+\mathcal{O}(\epsilon^2) \bigg) (\sum_{j \in \mathcal{N}_{US}} {\theta}^{us}_j)^2  -\nonumber \\ &  \hspace{-1.5cm} 2 \bigg(1+\mathcal{O}(K\epsilon)\bigg)\frac{\bigg(\frac{ \alpha \epsilon M L}{2 \delta}+\mathcal{O}(\epsilon^2) \bigg) }{( \alpha L + \alpha\beta + \mathcal{O}(\epsilon^2))}\bigg(1 + \alpha L + \frac{\alpha \epsilon M}{2}+\mathcal{O}(\epsilon^2)\bigg)^K\bigg(1 - \alpha \beta + \frac{\alpha \epsilon M}{2}+\mathcal{O}(\epsilon^2)\bigg)^{K}\sum_{j \in \mathcal{N}_{US}}{\theta}^{us}_j \sum_{i \in \mathcal{N}_{S}}{\theta}^{s}_i  - \nonumber \\ & \hspace{-1.5cm} 2 \bigg(1+\mathcal{O}(K\epsilon)\bigg)\frac{\bigg(\frac{ \alpha \epsilon M L}{2 \delta}+\mathcal{O}(\epsilon^2) \bigg)}{( \alpha L + \alpha\beta + \mathcal{O}(\epsilon^2))}\bigg(1 + \alpha L + \frac{\alpha \epsilon M}{2}+\mathcal{O}(\epsilon^2)\bigg)^{2K}\sum_{j \in \mathcal{N}_{US}}{\theta}^{us}_j \sum_{i \in \mathcal{N}_{S}}{\theta}^{s}_i
     \bigg]
     \end{align}
     \begin{align}
  & \hspace{-0.5cm} >  \epsilon^2\bigg[\bigg(1 - \alpha L - \frac{\alpha \epsilon M}{2} -\mathcal{O}(\epsilon^2)\bigg)^{2K} \sum_{i \in \mathcal{N}_{S}} ({\theta}^{s}_i)^2  - 2nK\bigg(1 - \alpha \beta + \frac{\alpha \epsilon M}{2}+\mathcal{O}(\epsilon^2)\bigg)^{2K-1}  \bigg(\frac{ \alpha \epsilon M L}{2 \delta}+\mathcal{O}(\epsilon^2) \bigg)  \sum_{i \in \mathcal{N}_{S}} ({\theta}^{s}_i)^2 + \nonumber \\ & \hspace{-0.5cm}\bigg(1 + \alpha \beta - \frac{\alpha \epsilon M}{2}-\mathcal{O}(\epsilon^2)\bigg)^{2K}\sum_{j \in \mathcal{N}_{US}}( {\theta}^{us}_j)^2   -  2nK\bigg(1 + \alpha L + \frac{\alpha \epsilon M}{2}+\mathcal{O}(\epsilon^2)\bigg)^{2K-1}  \bigg(\frac{ \alpha \epsilon M L}{2 \delta}+\mathcal{O}(\epsilon^2) \bigg) \sum_{j \in \mathcal{N}_{US}} ({\theta}^{us}_j)^2  - \nonumber \\ & \hspace{-0.5cm}  n\bigg(1+\mathcal{O}(K\epsilon)\bigg)\frac{\bigg(\frac{ \alpha \epsilon M L}{2 \delta}+\mathcal{O}(\epsilon^2) \bigg) }{( \alpha L + \alpha\beta + \mathcal{O}(\epsilon^2))}\bigg(1 + \alpha L + \frac{\alpha \epsilon M}{2}+\mathcal{O}(\epsilon^2)\bigg)^K\bigg(1 - \alpha \beta + \frac{\alpha \epsilon M}{2}+\mathcal{O}(\epsilon^2)\bigg)^{K}\bigg(\sum_{i \in \mathcal{N}_{S}}({\theta}^{s}_i)^2 + \sum_{j \in \mathcal{N}_{US}}({\theta}^{us}_j)^2 \bigg) - \nonumber \\ & \hspace{-0.5cm} n\bigg(1+\mathcal{O}(K\epsilon)\bigg)\frac{\bigg(\frac{ \alpha \epsilon M L}{2 \delta}+\mathcal{O}(\epsilon^2) \bigg)}{( \alpha L + \alpha\beta + \mathcal{O}(\epsilon^2))}\bigg(1 + \alpha L + \frac{\alpha \epsilon M}{2}+\mathcal{O}(\epsilon^2)\bigg)^{2K}\bigg(\sum_{i \in \mathcal{N}_{S}}({\theta}^{s}_i)^2 + \sum_{j \in \mathcal{N}_{US}}({\theta}^{us}_j)^2 \bigg)
  \bigg].   \label{ineqthm1}
     \end{align}
     Note that in the last step we have used the following inequalities:
     \begin{align}
     n\sum_{i \in \mathcal{N}_{S}}({\theta}^{s}_i)^2 & \geq (\sum_{i \in \mathcal{N}_{S}}{\theta}^{s}_i)^2,  \nonumber \\
     n\sum_{j \in \mathcal{N}_{US}}({\theta}^{us}_j)^2 & \geq (\sum_{j \in \mathcal{N}_{US}}{\theta}^{us}_j)^2, \ \text{and} \nonumber \\
     2\sum_{j \in \mathcal{N}_{US}}{\theta}^{us}_j \sum_{i \in \mathcal{N}_{S}}{\theta}^{s}_i &\leq n\sum_{j \in \mathcal{N}_{US}}({\theta}^{us}_j)^2+ n\sum_{i \in \mathcal{N}_{S}}({\theta}^{s}_i)^2,  \nonumber
     \end{align}
     where $n$ is the dimension of the domain of the function $f(\cdot)$. The above condition can be more compactly written as
     \begin{align}
       \epsilon^2 \geq \inf_{\tau}\norm{\tilde{\u}_K^{\tau}}^2 > &  \epsilon^2 \Psi(K), \label{thm1}
     \end{align}
     where we have that
     \begin{align}
\Psi(K) =& \bigg(c_1^{2K} -2Kc_2^{2K-1}b_1 - b_2c_3^Kc_2^K - b_2c_3^{2K}\bigg)\sum_{i \in \mathcal{N}_{S}}({\theta}^{s}_i)^2 + \bigg( c_4^{2K} - 2Kc_3^{2K-1}b_1- b_2c_3^Kc_2^K -b_2c_3^{2K}\bigg)\sum_{j \in \mathcal{N}_{US}}({\theta}^{us}_j)^2,
     \end{align}
and $c_1 = \bigg(1 - \alpha L - \frac{\alpha \epsilon M}{2}- \mathcal{O}(\epsilon^2) \bigg)$, $c_2 =\bigg(1 - \alpha \beta + \frac{\alpha \epsilon M}{2}+\mathcal{O}(\epsilon^2)\bigg) $, $c_3 = \bigg(1 + \alpha L + \frac{\alpha \epsilon M}{2}+\mathcal{O}(\epsilon^2)\bigg)$, $c_4 =\bigg(1 + \alpha \beta - \frac{\alpha \epsilon M}{2}-\mathcal{O}(\epsilon^2)\bigg)$, $b_1 = \bigg(\frac{\alpha\epsilon M L n}{2 \delta}+\mathcal{O}(\epsilon^2)\bigg)$ and $b_2 = \frac{\bigg(\frac{\alpha\epsilon M L n}{2 \delta}+\mathcal{O}(\epsilon^2)\bigg)\bigg(1+\mathcal{O}(K\epsilon)\bigg)}{ \bigg( \alpha L + \alpha\beta + \mathcal{O}(\epsilon^2)\bigg)}$.

The condition in \eqref{thm1} holds for all such $K$ where $\inf_{\tau}\norm{\tilde{\u}_K^{\tau}}^2 \leq \epsilon^2$. Therefore to obtain $K^{\iota}$ defined in \eqref{upperexittime}, we need to solve for $K$ where $ \epsilon^2 \leq \epsilon^2 \Psi(K) $ or equivalently $1 \leq \Psi(K)$ where the condition $\inf_{\tau}\norm{\tilde{\u}_K^{\tau}}^2 \leq \epsilon^2$ gets inverted using inequality \eqref{thm1}.
\subsubsection{Claim for the value of $K$ in Theorem \ref{theorem1}:}
Since the infimum in \eqref{thm1} is taken over all $\tau$, the condition in \eqref{thm1} holds true for all $K$ in the range $1 \leq K < \sup_{\tau} \bigg\{K_{exit}^{\tau} \bigg\}  $.\\
\textit{Proof of the claim:}
 Recall that from the definition of $K^{\iota}$ from \eqref{upperexittime}, $ K^{\iota}$ satisfies the following condition:
\begin{align}
   \inf_{\tau}\norm{\tilde{\u}_{K^{\iota}-1}^{\tau}}^2 \leq \epsilon^2 < \inf_{\tau}\norm{\tilde{\u}_{K^{\iota}}^{\tau}}^2,   \label{check1}
\end{align}
where the lower bound implies that the infimum over all $\tau$-parameterized approximate trajectories has not yet escaped the ball $\mathcal{B}_{\epsilon}(\x^*)$. Let there exist some $\bar{K}$ where $\bar{K} \in \bigg\{1,2,...,\sup_{\tau} \bigg\{K_{exit}^{\tau} \bigg\}\bigg\}$ such that the condition in \eqref{thm1} holds for all $K \in  \bigg\{1,...,\bar{K}-1\bigg \}$ and fails to hold for all $K \in  \bigg\{\bar{K},...,\sup_{\tau} \bigg\{K_{exit}^{\tau} \bigg\}\bigg\}$, i.e. we have the condition $\epsilon^2  \leq   \epsilon^2 \Psi({K}) < \inf_{\tau}\norm{\tilde{\u}_{{K}}^{\tau}}^2 $ for $K \geq \bar{K}$. This implies that
\begin{align}
  \epsilon^2 \Psi(\bar{K}-1)<\inf_{\tau}\norm{\tilde{\u}_{\bar{K}-1}^{\tau}}^2 \leq \epsilon^2  \leq &  \epsilon^2 \Psi(\bar{K}) < \inf_{\tau}\norm{\tilde{\u}_{\bar{K}}^{\tau}}^2. \label{check2}
\end{align}
From conditions \eqref{check1} and \eqref{check2} we get that $\bar{K} = K^{\iota}$.
Since $\bar{K} \in \bigg\{1,2,...,\sup_{\tau} \bigg\{K_{exit}^{\tau} \bigg\}\bigg\}$ we have that $ K^{\iota}=\bar{K} \leq \sup_{\tau} \bigg\{K_{exit}^{\tau} \bigg\} \leq K^{\iota}$. Hence we must have that $\bar{K} = \sup_{\tau} \bigg\{K_{exit}^{\tau} \bigg\}$.
     \endproof

\section{Proof of Theorem \ref{theorem2} (Exit time for the infimum of $\epsilon$-precision trajectories)} \label{Appendix D}

\proof{}

 Further simplifying the inequality in \eqref{ineqthm1} by dropping order $\mathcal{O}(\epsilon^2)$ and $\mathcal{O}(K\epsilon)$ terms (for $K\epsilon \ll 1$) appearing on its right hand side and using \eqref{thm1}, we get the following approximate lower bound:
\begin{align}
\hspace{-2cm}1 \gtrapprox  &  \bigg(\bigg[ \bigg(1 - \alpha L - \frac{\alpha \epsilon M}{2}\bigg)^{2K}  - 2K\bigg(1 - \alpha \beta + \frac{\alpha \epsilon M}{2}\bigg)^{2K-1}  \frac{\alpha \epsilon M Ln}{2 \delta} \bigg] \sum_{i \in \mathcal{N}_{S}} ({\theta}^{s}_i)^2 + \nonumber \\ & \bigg[\bigg(1 + \alpha \beta - \frac{\alpha \epsilon M}{2}\bigg)^{2K} -  2K\bigg(1 + \alpha L + \frac{\alpha \epsilon M}{2}\bigg)^{2K-1}  \frac{\alpha \epsilon M Ln}{2 \delta}\bigg]  \sum_{j \in \mathcal{N}_{US}} ({\theta}^{us}_j)^2  \nonumber \\ &-\frac{\alpha\epsilon M L n}{2 \delta ( \alpha L + \alpha\beta )}\bigg(1 + \alpha L + \frac{\alpha \epsilon M}{2}\bigg)^K\bigg(1 - \alpha \beta + \frac{\alpha \epsilon M}{2}\bigg)^{K}\bigg(\sum_{i \in \mathcal{N}_{S}}({\theta}^{s}_i)^2 + \sum_{j \in \mathcal{N}_{US}}({\theta}^{us}_j)^2 \bigg)  \nonumber \\ & -\frac{\alpha\epsilon M Ln }{2 \delta ( \alpha L + \alpha\beta)}\bigg(1 + \alpha L + \frac{\alpha \epsilon M}{2}\bigg)^{2K}\bigg)\bigg(\sum_{i \in \mathcal{N}_{S}}({\theta}^{s}_i)^2 + \sum_{j \in \mathcal{N}_{US}}({\theta}^{us}_j)^2 \bigg)
 \\
\hspace{-2cm}1 \gtrapprox  &  \bigg(\bigg[ \bigg(1 - \alpha L - \frac{\alpha \epsilon M}{2}\bigg)^{2K}  - 2K\bigg(1 - \alpha \beta + \frac{\alpha \epsilon M}{2}\bigg)^{2K-1}  \frac{\alpha \epsilon M Ln}{2 \delta} \bigg] \sum_{i \in \mathcal{N}_{S}} ({\theta}^{s}_i)^2 + \nonumber \\ & \bigg[\bigg(1 + \alpha \beta - \frac{\alpha \epsilon M}{2}\bigg)^{2K} -  2K\bigg(1 + \alpha L + \frac{\alpha \epsilon M}{2}\bigg)^{2K-1}  \frac{\alpha \epsilon M Ln}{2 \delta}\bigg]  \sum_{j \in \mathcal{N}_{US}} ({\theta}^{us}_j)^2   - \epsilon M Ln \frac{\bigg(1 + \alpha L + \frac{\alpha \epsilon M}{2}\bigg)^{2K} }{ \delta (L+\beta )}
\bigg),
\end{align}
where in the last step we used the relation $\bigg(\sum_{i \in \mathcal{N}_{S}}({\theta}^{s}_i)^2 + \sum_{j \in \mathcal{N}_{US}}({\theta}^{us}_j)^2 \bigg) =1 $ and the inequality $\bigg(1 - \alpha \beta + \frac{\alpha \epsilon M}{2}\bigg)< \bigg(1 + \alpha L + \frac{\alpha \epsilon M}{2}\bigg)$.
Now, if we substitute the step size $\alpha = \frac{1}{L}$, we get the following approximate inequality:
\begin{align}
\hspace{-2cm}1 \gtrapprox  &  \bigg(\bigg[ \bigg( - \frac{ \epsilon M}{2L}\bigg)^{2K}  - 2K\bigg(1 - \frac{\beta}{L} + \frac{\epsilon M}{2L}\bigg)^{2K-1}  \frac{ \epsilon M n}{2 \delta} \bigg] \sum_{i \in \mathcal{N}_{S}} ({\theta}^{s}_i)^2 + \nonumber \\ & \bigg[\bigg(1 + \frac{\beta}{L} - \frac{\epsilon M}{2L}\bigg)^{2K} -  2K\bigg(2 + \frac{ \epsilon M}{2L}\bigg)^{2K-1}  \frac{\epsilon M n}{2 \delta}\bigg]  \sum_{j \in \mathcal{N}_{US}} ({\theta}^{us}_j)^2   - \epsilon M Ln \frac{\bigg(2 + \frac{ \epsilon M}{2L}\bigg)^{2K} }{ \delta (L+ \beta )}
\bigg)  \label{originallowerbound}\\
\hspace{-2cm}1 \gtrapprox  &  \bigg(\bigg[  - 2K\bigg(1 - \frac{\beta}{L} + \frac{\epsilon M}{2L}\bigg)^{2K-1}  \frac{ \epsilon M n}{2 \delta} \bigg] \sum_{i \in \mathcal{N}_{S}} ({\theta}^{s}_i)^2 + \nonumber \\ & \bigg[\bigg(1 + \frac{\beta}{L} - \frac{\epsilon M}{2L}\bigg)^{2K} -  2K\bigg(2 + \frac{ \epsilon M}{2L}\bigg)^{2K-1}  \frac{\epsilon M n}{2 \delta}\bigg]  \sum_{j \in \mathcal{N}_{US}} ({\theta}^{us}_j)^2 - \epsilon M Ln \frac{\bigg(2 + \frac{ \epsilon M}{2L}\bigg)^{2K} }{ \delta ( L+\beta )}
\bigg),  \label{lbtranscend}
\end{align}
where in the last step we dropped the $\bigg( - \frac{ \epsilon M}{2L}\bigg)^{2K} $ term from right hand side.

In order to obtain $K^{\iota}$ and hence the exit time $K_{exit}$, we need to solve for values of $K$ where the approximate inequality in \eqref{lbtranscend} becomes an equality. Hence, we look into the two possible cases for this value $K$, i.e., large $K$ and small $K$. Note that in the next subsections we only consider those cases where our unstable projection $ \sum_{j \in \mathcal{N}_{US}} ({\theta}^{us}_j)^2$ is not too close to $0$. We now obtain the exit time $K_{exit}$ for the two cases.
\subsubsection{Case 1---Large $K$:}
If $K$ is large with $K = \mathcal{O}\bigg(\frac{1}{\epsilon}\bigg)$ then we can use the Lambert $W$ function \cite{corless1996lambertw} to solve the above transcendental inequality \eqref{lbtranscend}. Specifically for obtaining linear escape rates i.e., $K = \mathcal{O}\bigg(\log\bigg(\frac{1}{\epsilon}\bigg)\bigg)$, we set $\frac{1}{\bigg(2 + \frac{ \epsilon M}{2L}\bigg)^{2K}} = \rho \epsilon^c  $ for some $\rho > 0$, $c>0$, $\bigg(1 - \frac{\beta}{L} + \frac{\epsilon M}{2L}\bigg)^{2K} = \eta \epsilon^d$ for some $\eta > 0$, $d>0$ where $\bigg(1 - \frac{\beta}{L} + \frac{\epsilon M}{2L}\bigg)<1$ and divide both sides of \eqref{lbtranscend} by the term $\bigg(2 + \frac{ \epsilon M}{2L}\bigg)^{2K}$ to get the following approximate inequality:
\begin{align}
\hspace{-1cm}\frac{1}{\bigg(2 + \frac{ \epsilon M}{2L}\bigg)^{2K}} \gtrapprox  &  \bigg(\bigg[ \underbrace{- 2K\frac{\bigg(1 - \frac{\beta}{L} + \frac{\epsilon M}{2L}\bigg)^{2K-1}}{\bigg(2 + \frac{ \epsilon M}{2L}\bigg)^{2K}} \frac{ \epsilon M n}{2 \delta} \bigg] \sum_{i \in \mathcal{N}_{S}} ({\theta}^{s}_i)^2}_{F_1} + \nonumber \\ & \bigg[\frac{\bigg(1 + \frac{\beta}{L} - \frac{\epsilon M}{2L}\bigg)^{2K}}{\bigg(2 + \frac{ \epsilon M}{2L}\bigg)^{2K}} -  2K\bigg(2 + \frac{ \epsilon M}{2L}\bigg)^{-1}  \frac{\epsilon M n}{2 \delta}\bigg]  \sum_{j \in \mathcal{N}_{US}} ({\theta}^{us}_j)^2 -  \frac{\epsilon M Ln}{ \delta (L+ \beta )}
\bigg) . \label{largekbound}
\end{align}
Dropping the first term $F_1$ on right hand side for large $K$ (this term has order $\mathcal{O}\bigg(\epsilon^{(1+c+d)}\log \bigg(\frac{1}{\epsilon}\bigg) \bigg)$ with $c>0$, $d>0$) and making the substitution of $\rho \epsilon^c$ on the left hand side, we get the following bound:
\begin{align}
\hspace{-0.3cm}\rho \epsilon^c \hspace{0.08cm}\gtrapprox    \bigg(\frac{1 + \frac{\beta}{L} - \frac{\epsilon M}{2L}}{2 + \frac{ \epsilon M}{2L}}\bigg)^{2K} \sum_{j \in \mathcal{N}_{US}} ({\theta}^{us}_j)^2 - & 2K\bigg(2 + \frac{ \epsilon M}{2L}\bigg)^{-1}  \frac{\epsilon M n}{2 \delta} \sum_{j \in \mathcal{N}_{US}} ({\theta}^{us}_j)^2 -  \frac{\epsilon M Ln}{ \delta (L+\beta)} \\
\hspace{-0.3cm}   \bigg(\frac{1 + \frac{\beta}{L} - \frac{\epsilon M}{2L}}{2 + \frac{ \epsilon M}{2L}}\bigg)^{2K}   \lessapprox \hspace{0.1cm}&   2K\bigg(2 + \frac{ \epsilon M}{2L}\bigg)^{-1}  \frac{\epsilon M n}{2 \delta} + \frac{\epsilon\bigg(\rho \epsilon^{(c-1)} + \frac{ M Ln}{ \delta (L+ \beta )} \bigg)}{\sum_{j \in \mathcal{N}_{US}} ({\theta}^{us}_j)^2 }. \label{well-cond}
\end{align}
When the problem is well conditioned, i.e., $\bigg(1 - \frac{\beta}{L} + \frac{\epsilon M}{2L}\bigg)<1$ or equivalently $\frac{\beta}{L}> \frac{\epsilon M}{2L}$, then we are guaranteed fast escape under good initial unstable projections. Now, solving for the values of $K$ where the inequality \eqref{well-cond} becomes equality, we make use of the general transcendental equation $ q^x = ax + b$ whose solution is given by
\begin{align}
x & = - \frac{W(-\frac{\log q}{a}q^{-\frac{b}{a}})}{\log q} - \frac{b}{a},
\end{align}
where $W(\cdot)$ is the Lambert $W$ function. On comparing the coefficients, we have $x = 2K $ and the constants as follows:
\begin{align}
a  = \bigg(2 + \frac{ \epsilon M}{2L}\bigg)^{-1}  \frac{\epsilon M n}{2 \delta} ,
b & = \frac{\epsilon\bigg(\rho\epsilon^{(c-1)} + \frac{ M Ln}{ \delta (L+ \beta)} \bigg)}{\sum_{j \in \mathcal{N}_{US}} ({\theta}^{us}_j)^2 } ,
q  =  \bigg(\frac{1 + \frac{\beta}{L} - \frac{\epsilon M}{2L}}{2 + \frac{ \epsilon M}{2L}}\bigg).
\end{align}

For large values of any argument $y$, the Lambert $W$ function is bounded by $W(y) \leq \log(y)$. If the quantity $\sum_{j \in \mathcal{N}_{US}} ({\theta}^{us}_j)^2 $ is not too close to $0$ and is lower bounded, i.e., $ \sum_{j \in \mathcal{N}_{US}} ({\theta}^{us}_j)^2 \geq  \Delta$ then we have an initial projection onto the unstable subspace of the saddle point. Using the Lambert $W$ function bound and substituting the coefficients, we have following bound on $K$:
\begin{align}
\hspace{-2cm} 2K = &  \frac{1}{\log\bigg(\frac{2 + \frac{ \epsilon M}{2L}}{1 + \frac{\beta}{L} - \frac{\epsilon M}{2L}}\bigg)} W \bigg(\bigg(2 + \frac{ \epsilon M}{2L}\bigg) \frac{2 \delta}{\epsilon M n} \log\bigg(\frac{2 + \frac{ \epsilon M}{2L}}{1 + \frac{\beta}{L} - \frac{\epsilon M}{2L}}\bigg)\bigg(\frac{2 + \frac{ \epsilon M}{2L}}{1 + \frac{\beta}{L} - \frac{\epsilon M}{2L}}\bigg)^{\frac{2 \delta \bigg(2 + \frac{ \epsilon M}{2L}\bigg) \bigg(\rho\epsilon^{(c-1)} + \frac{ M Ln}{ \delta (L+ \beta)} \bigg)}{M n \sum_{j \in \mathcal{N}_{US}} ({\theta}^{us}_j)^2 }}\bigg) - \nonumber \\ & \hspace{2cm} \frac{2 \delta \bigg(2 + \frac{ \epsilon M}{2L}\bigg) \bigg(\rho\epsilon^{(c-1)} + \frac{ M Ln}{ \delta (L+ \beta)} \bigg)}{M n \sum_{j \in \mathcal{N}_{US}} ({\theta}^{us}_j)^2 }  \\
\hspace{-2cm} 2K \leq &  \frac{1}{\log\bigg(\frac{2 + \frac{ \epsilon M}{2L}}{1 + \frac{\beta}{L} - \frac{\epsilon M}{2L}}\bigg)} \log \bigg(\bigg(2 + \frac{ \epsilon M}{2L}\bigg) \frac{2 \delta}{\epsilon M n} \log\bigg(\frac{2 + \frac{ \epsilon M}{2L}}{1 + \frac{\beta}{L} - \frac{\epsilon M}{2L}}\bigg)\bigg(\frac{2 + \frac{ \epsilon M}{2L}}{1 + \frac{\beta}{L} - \frac{\epsilon M}{2L}}\bigg)^{\frac{2 \delta \bigg(2 + \frac{ \epsilon M}{2L}\bigg) \bigg(\rho\epsilon^{(c-1)} + \frac{ M Ln}{ \delta (L+ \beta)} \bigg)}{M n \sum_{j \in \mathcal{N}_{US}} ({\theta}^{us}_j)^2 }}\bigg) - \nonumber \\ & \hspace{2cm} \frac{2 \delta \bigg(2 + \frac{ \epsilon M}{2L}\bigg) \bigg(\rho\epsilon^{(c-1)} + \frac{ M Ln}{ \delta (L+ \beta)} \bigg)}{M n \sum_{j \in \mathcal{N}_{US}} ({\theta}^{us}_j)^2 } \\
\hspace{-2cm} 2K \leq &  \frac{\log \bigg(\bigg(2 + \frac{ \epsilon M}{2L}\bigg) \log\bigg(\frac{2 + \frac{ \epsilon M}{2L}}{1 + \frac{\beta}{L} - \frac{\epsilon M}{2L}}\bigg)\frac{2 \delta}{\epsilon M n}\bigg)}{\log\bigg(\frac{2 + \frac{ \epsilon M}{2L}}{1 + \frac{\beta}{L} - \frac{\epsilon M}{2L}}\bigg)} + \frac{2 \delta \bigg(2 + \frac{ \epsilon M}{2L}\bigg) \bigg(\rho\epsilon^{(c-1)} + \frac{ M Ln}{ \delta (L+ \beta)} \bigg)}{M n \sum_{j \in \mathcal{N}_{US}} ({\theta}^{us}_j)^2 } -  \frac{2 \delta \bigg(2 + \frac{ \epsilon M}{2L}\bigg) \bigg(\rho\epsilon^{(c-1)} + \frac{ M Ln}{ \delta (L+ \beta)} \bigg)}{M n \sum_{j \in \mathcal{N}_{US}} ({\theta}^{us}_j)^2 }\\
K \leq & \frac{\log \bigg(\bigg(2 + \frac{ \epsilon M}{2L}\bigg) \log\bigg(\frac{2 + \frac{ \epsilon M}{2L}}{1 + \frac{\beta}{L} - \frac{\epsilon M}{2L}}\bigg)\frac{2 \delta}{\epsilon M n}\bigg)}{2\log\bigg(\frac{2 + \frac{ \epsilon M}{2L}}{1 + \frac{\beta}{L} - \frac{\epsilon M}{2L}}\bigg)} = \mathcal{O}\bigg(\log\bigg(\frac{1}{\epsilon}\bigg)\bigg). \label{escapetime1}
\end{align}
Notice that the $K$ solved here is an approximate solution to \eqref{well-cond} where the inequality in \eqref{well-cond} gets inverted. Since the condition \eqref{lbtranscend} gets reversed at $K=K^{\iota}$, we therefore get the condition $K^{\iota} \lessapprox \frac{\log \bigg(\bigg(2 + \frac{ \epsilon M}{2L}\bigg) \log\bigg(\frac{2 + \frac{ \epsilon M}{2L}}{1 + \frac{\beta}{L} - \frac{\epsilon M}{2L}}\bigg)\frac{2 \delta}{\epsilon M n}\bigg)}{2\log\bigg(\frac{2 + \frac{ \epsilon M}{2L}}{1 + \frac{\beta}{L} - \frac{\epsilon M}{2L}}\bigg)}$ and using the fact that $ K_{exit} < K^{\iota}$ gives the desired conclusion of $K_{exit} \leq K^{\iota} = \mathcal{O}\bigg(\log\bigg(\frac{1}{\epsilon}\bigg)\bigg)$. \revise{The bound $\epsilon < \frac{2 \beta}{M}$ follows from the fact that $\frac{\beta}{L} > \frac{\epsilon M}{2L} $.}

Hence, we have escape rates of order $\mathcal{O}\bigg(\log\bigg(\frac{1}{\epsilon}\bigg)\bigg)$ for the case when our problem is well conditioned and does not have a very small unstable projection. It is remarked that this is only an upper bound on $K$ and the iterate is likely to escape way before this time. Also, this result supports our analysis of the trajectory function for values of $K = \mathcal{O}\bigg(\frac{1}{\epsilon}\bigg)$.

It is worth mentioning that dropping of the first term $F_1$ with order $\mathcal{O}\bigg(\epsilon^{(1+c+d)}\log \bigg(\frac{1}{\epsilon}\bigg) \bigg)$ from the right hand side of inequality \eqref{largekbound} is justified since from the particular upper bound of $K^{\iota}$ from \eqref{escapetime1} it can be inferred that $c>1$.

From the substitution $\frac{1}{\bigg(2 + \frac{ \epsilon M}{2L}\bigg)^{2K}} = \rho \epsilon^c$ where $2K = \frac{\log \bigg(\bigg(2 + \frac{ \epsilon M}{2L}\bigg) \log\bigg(\frac{2 + \frac{ \epsilon M}{2L}}{1 + \frac{\beta}{L} - \frac{\epsilon M}{2L}}\bigg)\frac{2 \delta}{\epsilon M n}\bigg)}{\log\bigg(\frac{2 + \frac{ \epsilon M}{2L}}{1 + \frac{\beta}{L} - \frac{\epsilon M}{2L}}\bigg)}$ we have that
\begin{align}
\log \bigg(\frac{1}{\rho\epsilon^c } \bigg) & = 2K\log\bigg(2 + \frac{ \epsilon M}{2L}\bigg)   \\
c\log \bigg(\frac{1}{\sqrt[c]{\rho}\epsilon } \bigg) & =  \frac{\log \bigg(\bigg(2 + \frac{ \epsilon M}{2L}\bigg) \log\bigg(\frac{2 + \frac{ \epsilon M}{2L}}{1 + \frac{\beta}{L} - \frac{\epsilon M}{2L}}\bigg)\frac{2 \delta}{\epsilon M n}\bigg)}{\log\bigg(\frac{2 + \frac{ \epsilon M}{2L}}{1 + \frac{\beta}{L} - \frac{\epsilon M}{2L}}\bigg)} \log\bigg(2 + \frac{ \epsilon M}{2L}\bigg)   \\
c & = \frac{ \log\bigg(2 + \frac{ \epsilon M}{2L}\bigg)}{ \log\bigg(2 + \frac{ \epsilon M}{2L}\bigg) - \log\bigg(1 + \frac{\beta}{L} - \frac{\epsilon M}{2L}  \bigg)} > 1,
\end{align}
where we have $\log \bigg(\frac{1}{\sqrt[c]{\rho}\epsilon } \bigg) =  \log \bigg(\bigg(2 + \frac{ \epsilon M}{2L}\bigg) \log\bigg(\frac{2 + \frac{ \epsilon M}{2L}}{1 + \frac{\beta}{L} - \frac{\epsilon M}{2L}}\bigg)\frac{2 \delta}{\epsilon M n}\bigg)$. Now with $c>1$, we will have the following condition for any $d>0$:
\begin{align}
\lim_{\epsilon \to 0^+}\frac{\epsilon^{(1+c+d)}\log \bigg(\frac{1}{\epsilon}\bigg)}{\epsilon^2} = 0.
\end{align}
Hence, for sufficiently small $\epsilon$, the term $F_1$ can be of at most order $\mathcal{O}(\epsilon^2)$.

\textbf{Comments on the projection $\sum_{j \in \mathcal{N}_{US}} ({\theta}^{us}_j)^2 $ :} Recall that from \eqref{well-cond}, we solved for values of $K$ where this inequality becomes an equality. However, this solution for such $K$ may not necessarily exist. For instance, the left hand side of \eqref{well-cond} given by $ \bigg(\frac{1 + \frac{\beta}{L} - \frac{\epsilon M}{2L}}{2 + \frac{ \epsilon M}{2L}}\bigg)^{2K}$ is a decreasing function of $K$ whereas the right hand side of this inequality given by $ 2K\bigg(2 + \frac{ \epsilon M}{2L}\bigg)^{-1}  \frac{\epsilon M n}{2 \delta} + \frac{\epsilon\bigg(\rho \epsilon^{(c-1)} + \frac{ M Ln}{ \delta (L+ \beta )} \bigg)}{\sum_{j \in \mathcal{N}_{US}} ({\theta}^{us}_j)^2 }$ is an increasing function of $K$. Hence for a solution $K$ to exist where these two quantities become equal, we must necessarily have that
\begin{align}
\bigg(\frac{1 + \frac{\beta}{L} - \frac{\epsilon M}{2L}}{2 + \frac{ \epsilon M}{2L}}\bigg)^{2K} \bigg\vert_{K=0}   > &   2K\bigg(2 + \frac{ \epsilon M}{2L}\bigg)^{-1}  \frac{\epsilon M n}{2 \delta}\bigg\vert_{K=0} + \frac{\epsilon\bigg(\rho\epsilon^{(c-1)} + \frac{ M Ln}{ \delta (L+ \beta )} \bigg)}{\sum_{j \in \mathcal{N}_{US}} ({\theta}^{us}_j)^2 } \\
\sum_{j \in \mathcal{N}_{US}} ({\theta}^{us}_j)^2  > & \epsilon\bigg(\rho\epsilon^{(c-1)} + \frac{ M Ln}{ \delta (L+ \beta )} \bigg) > \epsilon \frac{ M Ln}{ \delta (L+ \beta )},
\end{align}
where we can set $\Delta > \epsilon \frac{ M Ln}{ \delta (L+ \beta )} $ and therefore require the condition $\sum_{j \in \mathcal{N}_{US}} ({\theta}^{us}_j)^2 \geq  \Delta$. Note that this is only a necessary condition for the existence of $K$ from \eqref{escapetime1} but is not sufficient.

\subsubsection{Case 2---Small K:}
Recall that while developing the inequality \eqref{largekbound} from \eqref{lbtranscend}, we used the fact that $K$ is sufficiently large. However, for very small values of $K$, i.e., $K < \mathcal{O}\bigg(\log\bigg(\frac{1}{\epsilon}\bigg)\bigg)$, the transformation of the inequality \eqref{lbtranscend} into \eqref{largekbound} may not necessarily hold true. In that case, a different approach is required to solve for $K$. Since the new solutions for $K$ will be very small values, we can skip the analysis for small $K$ case and extrapolate it to the previous result of $K \leq K^{\iota} = \mathcal{O}\bigg(\log\bigg(\frac{1}{\epsilon}\bigg)\bigg) $ which is a linear exit time solution. We now complete the proof of Theorem \ref{theorem2} by establishing one last result.
\subsubsection{Claim: The set of $\epsilon$-precision trajectories with linear exit times from the ball $\mathcal{B}_{\epsilon}(\x^*)$ is non-empty}\label{nonvacuousclaim}
\proof
Observe that from \eqref{lbtranscend}, we need to find $K$ where this approximate inequality becomes an equality. Let the initial condition be such that $\sum_{j \in \mathcal{N}_{US}} ({\theta}^{us}_j)^2=1$; then \eqref{lbtranscend} can be given by
\begin{align}
\hspace{-2cm}1 \gtrapprox  &   \bigg(1 + \frac{\beta}{L} - \frac{\epsilon M}{2L}\bigg)^{2K} -  \bigg[2K\bigg(2 + \frac{ \epsilon M}{2L}\bigg)^{2K-1}  \frac{\epsilon M n}{2 \delta}  + \epsilon M Ln \frac{\bigg(2 + \frac{ \epsilon M}{2L}\bigg)^{2K} }{ \delta ( L+\beta )}\bigg] \\
\hspace{-2cm}\bigg(2 + \frac{ \epsilon M}{2L}\bigg)^{-2K} \gtrapprox  &  \bigg(\frac{1 + \frac{\beta}{L} - \frac{\epsilon M}{2L}}{2 + \frac{ \epsilon M}{2L}}\bigg)^{2K} -  \underbrace{\bigg[2K\bigg(2 + \frac{ \epsilon M}{2L}\bigg)^{-1}  \frac{\epsilon M n}{2 \delta}  +  \frac{\epsilon M Ln }{ \delta ( L+\beta )}\bigg]}_{L_1}. \label{nonemptycondition}
\end{align}
It is easy to infer that the right-hand side of \eqref{nonemptycondition} is negative for $K=\frac{\log \bigg(\bigg(2 + \frac{ \epsilon M}{2L}\bigg) \log\bigg(\frac{2 + \frac{ \epsilon M}{2L}}{1 + \frac{\beta}{L} - \frac{\epsilon M}{2L}}\bigg)\frac{2 \delta}{\epsilon M n}\bigg)}{2\log\bigg(\frac{2 + \frac{ \epsilon M}{2L}}{1 + \frac{\beta}{L} - \frac{\epsilon M}{2L}}\bigg)}$ where this value of $K$ comes from \eqref{escapetime1}. Hence, the approximate inequality in \eqref{nonemptycondition} holds for this value of $K$. However, for small positive values of $K$, one can check that the right-hand side of \eqref{nonemptycondition} is greater than its left-hand side, provided $\epsilon$ is sufficiently small and the problem is well-conditioned. This is because the term $L_1$ on the right-hand side of \eqref{nonemptycondition} is of order $\mathcal{O}(\epsilon)$ for small positive values of $K$ whereas we have that $\bigg(2 + \frac{ \epsilon M}{2L}\bigg)^{-2K} < \bigg(\frac{1 + \frac{\beta}{L} - \frac{\epsilon M}{2L}}{2 + \frac{ \epsilon M}{2L}}\bigg)^{2K}$ for any positive $K$.

Therefore, the approximate inequality in \eqref{nonemptycondition} becomes an equality for some $K= \mathcal{O}(\log(\epsilon^{-1}))$ and we have that $K^{\iota}= \mathcal{O}(\log(\epsilon^{-1}))$. As a result, the exit time $K_{exit}$ is linear for the initial condition $\sum_{j \in \mathcal{N}_{US}} ({\theta}^{us}_j)^2=1$ since $K_{exit}< K^{\iota}$. It should be noted that the proof of a linear exit time for the general initial condition $\Delta \leq \sum_{j \in \mathcal{N}_{US}} ({\theta}^{us}_j)^2 < 1$ can be developed along similar lines though it may require more effort.
\endproof

\section{Counterexample to the monotonicity property}
\label{Appendix E}
Consider a trajectory of the gradient descent method that satisfies the following boundary condition for some $\rho \in (0,1)$:
\begin{align}
 \frac{M\epsilon^2}{2 \beta (1- \rho)} = \frac{M\norm{\x_{0} - \x^*}^2}{2 \beta (1- \rho)} >  \langle\v_n, \x_{0} - \x^* \rangle \geq \underbrace{\langle\v_n, \x_{1} - \x^* \rangle \geq \frac{M\norm{\x_{1} - \x^*}^2}{2 \beta (1- \rho)}}_{I_1}, \label{crudeanalysiscond1}
\end{align}
where $ \norm{\x_{0} - \x^*} = \epsilon$ by definition. This trajectory violates the strict monotonicity of $ \langle\v_n, \x_{k} - \x^* \rangle$ for $k=0$. But it is straightforward to see that the condition $I_1$ along with \eqref{crudeanalysis3}--\eqref{crudeanalysis4.2}, in which $\epsilon$ is replaced by $\norm{\x_{1} - \x^*}$, ensures geometric growth of $\langle\v_n, \x_{k} - \x^* \rangle$ for all $k \geq 1$. This in turn guarantees linear exit time from $\mathcal{B}_{\epsilon}(\x^*)$ for the trajectory starting at $\x_0$, even though the strict monotonicity property is violated.  
Now our goal is to prove the existence of at-least one such trajectory that satisfies \eqref{crudeanalysiscond1} so as to construct the counterexample. In order to have the condition $ \langle\v_n, \x_{0} - \x^* \rangle \geq  \langle\v_n, \x_{1} - \x^* \rangle$ from \eqref{crudeanalysiscond1} with $\norm{\x_{1}- \x^* } \leq \epsilon$, we require
\begin{align}
    \langle\v_n, \x_{0} - \x^* \rangle &\geq  \langle\v_n, \x_{1} - \x^* \rangle \\ \Longleftrightarrow \langle\v_n, \x_{0} - \x^* \rangle &\geq  \langle\v_n, \x_{0} - \x^* \rangle - \alpha \langle \v_n, \nabla f(\x_0) \rangle \\
    \Longleftrightarrow \alpha \langle \v_n, \nabla f(\x_0) \rangle & \geq 0 \\
    \Longleftrightarrow \bigg\langle \v_n, \int_{p=0}^{p=1}\nabla^2 f(\x^* + p(\x_0 -\x^*)) (\x_0 -\x^*) dp \bigg\rangle & \geq 0 \\
    \Longleftrightarrow \langle \v_n, \nabla^2 f(\x^* ) (\x_0 -\x^*)  \rangle + \langle \v_n, P(\x_0) (\x_0 -\x^*)  \rangle &\geq 0, \label{crudeanalysis5}
\end{align}
where $ P(\x_0) = \int_{p=0}^{p=1}\nabla^2 f(\x^* + p(\x_0 -\x^*))  dp - \nabla^2 f(\x^* )$ and $ \norm{  P(\x_0) }_2 \leq \frac{M\epsilon}{2}$. Next, without loss of generality, write $\x_0 - \x^* = \epsilon\sum\limits_{j=1}^{n}a_j \v_j $ with $a_j \in [0,1]$ for all $j$, $ \sum_{j} a_j^2 =1 $, and $a_n = \frac{\langle\v_n, \x_{0} - \x^* \rangle}{\epsilon} = \frac{M \epsilon\sigma}{2 \beta (1-\rho)}$ for some positive $\sigma$ (note that $\sigma$ cannot be $1$ since we require the condition $ \langle\v_n, \x_{0} - \x^* \rangle < \frac{M \epsilon^2}{2 \beta (1-\rho)}$ from the left-hand-side of \eqref{crudeanalysiscond1}). Substituting the expression for $\x_0 - \x^*$ in \eqref{crudeanalysis5} followed by substituting $a_n$ and using the fact that $\lambda_n \geq -L$ from Assumption \textbf{A2} yields
\begin{align}
  \langle \v_n, \nabla f(\x_0) \rangle &= \langle \v_n, (\nabla^2 f(\x^*) + P(\x_0))(\x_0 -\x^*) \rangle \\ &=\epsilon  a_n \lambda_n + \langle \v_n, P(\x_0) (\x_0 -\x^*)  \rangle \geq -L \frac{M \epsilon^2\sigma}{2 \beta (1-\rho)}  + \epsilon\langle \v_n, P(\x_0) \sum\limits_{j=1}^{n}a_j \v_j  \rangle \geq 0.  \label{crudeanalysiscond2}
\end{align}
 Now there will exist some twice continuously differentiable function $f(\cdot)$ for which $\norm{P(\x_0)}_2 = \frac{M \epsilon}{2}$ for a given $\x_0$. Writing $ P(\x_0)$ in terms of the $\v_j$'s using the rank-one decomposition we get $ P(\x_0) = \sum\limits_{i=1}^n \sum\limits_{j=1}^n c_{ij} \v_i \v_j^T$ where $c_{ij} = c_{ji}$ since $ P(\x_0)$ is symmetric and we have the constraint $ \frac{M \epsilon}{2} \leq  \sqrt{\sum\limits_{i=1}^n \sum\limits_{j=1}^n c_{ij}^2} \leq \frac{M n \epsilon }{2}$. Hence one can fix $c_{in} = \frac{M \epsilon}{2} a_i $ for some twice continuously differentiable $f(\cdot)$ and substitute the resulting $ P(\x_0)$ into \eqref{crudeanalysiscond2} to get
 \begin{align}
     -L \frac{M \epsilon^2\sigma}{2 \beta (1-\rho)}  + \epsilon\langle \v_n, P(\x_0) \sum\limits_{j=1}^{n}a_j \v_j  \rangle &\geq 0 \\
     \Longleftrightarrow \epsilon\bigg\langle \sum\limits_{i=1}^n \sum\limits_{j=1}^n c_{ij} \v_i \v_j^T \v_n,  \sum\limits_{j=1}^{n}a_j \v_j \bigg \rangle &\geq L \frac{M \epsilon^2\sigma}{2 \beta (1-\rho)} \\
     \Longleftrightarrow \epsilon\bigg\langle \sum\limits_{i=1}^n  c_{in} \v_i ,  \sum\limits_{j=1}^{n}a_j \v_j \bigg \rangle &\geq L \frac{M \epsilon^2\sigma}{2 \beta (1-\rho)} \\
      \Longleftrightarrow \frac{M \epsilon^2}{2 } &\geq L \frac{M \epsilon^2\sigma}{2 \beta (1-\rho)} \\
      \Longleftrightarrow \frac{\beta (1 - \rho)}{L} & \geq \sigma. \label{crudeanalysiscond4}
 \end{align}

 Also, from \eqref{crudeanalysis1} we have that $\x_1 -\x^* = (\mathbf{I}- \alpha \nabla^2 f(\x^*))(\x_0 -\x^*) - \alpha r(\x_0)$ where $ \norm{r(\x_0) } \leq \frac{M \epsilon^2}{2} $. Hence, $ \norm{ \x_1 -\x^*} \leq \epsilon \sqrt{\sum\limits_{j=1}^n (1-\alpha \lambda_j)^2 a_j^2} + \frac{M \epsilon^2}{2} < \epsilon (\sqrt{ (1-\alpha \beta)^2 +4a_n^2} + \frac{M \epsilon}{2}) =  \epsilon (\sqrt{ (1-\alpha \beta)^2 +\mathcal{O}(\epsilon^2)} + \frac{M \epsilon}{2}) $. Next, using \eqref{crudeanalysiscond2} and \eqref{crudeanalysiscond4} we get that $ \langle\v_n, \x_{1} - \x^* \rangle = \langle\v_n, \x_{0} - \x^* \rangle - \alpha \langle\v_n, \nabla f(\x_{0}) \rangle =  \frac{M \epsilon^2\sigma}{2 \beta (1-\rho)}- \mathcal{O}(\epsilon^3) = \frac{M \epsilon^2}{2 L}- \mathcal{O}(\epsilon^3)$ by choosing $ \sigma = \frac{\beta (1 - \rho)}{L} - \mathcal{O}(\epsilon) $ when evaluating $\langle\v_n, \nabla f(\x_{0}) \rangle $. For the inequality $I_1$ to hold, we require $ \langle\v_n, \x_{1} - \x^* \rangle \geq \frac{M\norm{\x_{1} - \x^*}^2}{2 \beta (1- \rho)}$. This can be achieved by requiring the condition
 \begin{align}
     \langle\v_n, \x_{1} - \x^* \rangle =  \frac{M \epsilon^2\sigma}{2 \beta (1-\rho)}- \mathcal{O}(\epsilon^3) &\geq  \frac{M\epsilon^2}{2 \beta (1- \rho)} \bigg(\sqrt{ (1-\alpha \beta)^2 +\mathcal{O}(\epsilon^2)} + \frac{M \epsilon}{2}\bigg)^2 \geq \frac{M\norm{\x_{1} - \x^*}^2}{2 \beta (1- \rho)}\\
     \implies\frac{\beta (1 - \rho)}{L} - \mathcal{O}(\epsilon) =\sigma & \geq \bigg(\sqrt{ (1-\alpha \beta)^2 +\mathcal{O}(\epsilon^2)} + \frac{M \epsilon}{2}\bigg)^2 + \mathcal{O}(\epsilon). \label{crudeanalysiscond5}
 \end{align}
Now both \eqref{crudeanalysiscond4} and \eqref{crudeanalysiscond5} will be satisfied for $\alpha = \frac{1}{L}$ provided $\frac{\beta}{L}$ is close to $1$, $\epsilon$ is sufficiently small and $ \rho$ is not too large. Hence we have obtained a value of $\sigma$ and in turn $a_n = \frac{\langle \v_n, \x_0-\x^*\rangle}{\epsilon} = \frac{M\epsilon \sigma}{2 \beta (1-\rho)}$, i.e., the initial boundary condition, for which \eqref{crudeanalysiscond1} is satisfied on some twice continuously differentiable function $f(\cdot)$.
\end{appendices}

\ifx\undefined\BySame
\newcommand{\BySame}{\leavevmode\rule[.5ex]{3em}{.5pt}\ }
\fi
\ifx\undefined\textsc
\newcommand{\textsc}[1]{{\sc #1}}
\newcommand{\emph}[1]{{\em #1\/}}
\let\tmpsmall\small
\renewcommand{\small}{\tmpsmall\sc}
\fi

\end{document}